\documentclass[a4paper]{article} 

\usepackage[T1]{fontenc}
\usepackage[english]{babel}
\usepackage{hyperref}
\usepackage{amsmath,amsfonts,amscd}
\usepackage{mathtools}
\usepackage{calrsfs}
\usepackage{graphicx}
\usepackage{multirow}
\usepackage{multicol}
\usepackage{color}
\usepackage{bm}
\usepackage{subcaption}
\usepackage{setspace}
\usepackage{cuted}
\usepackage[sort&compress,numbers]{natbib}
\usepackage{mathptmx}      
\usepackage{ulem}

\providecommand{\keywords}[1]
{
  \small	
  \textbf{\textit{Keywords---}} #1
}

\usepackage{amsthm}
\newtheorem{remark}{Remark}
\newtheorem*{acknowledgements}{Acknowledgements}
\newcommand{\diff}{\mathrm{d}}
\newcommand{\divv}{\mathrm{div}}
\newcommand{\Nabla}{\nabla}
\newcommand{\Ha}{\mathcal{H}}

\DeclareMathAlphabet{\pazocal}{OMS}{zplm}{m}{n}
\DeclareMathAlphabet\mathbfcal{OMS}{cmsy}{b}{n}

\AtBeginDocument{\let\latexlabel\label}

\begin{document}


%


\title{A Damage Phase-Field Model for Fractional Viscoelastic Materials in Finite Strain}

\author{T. C.  da Costa Haveroth$^{1}$,  G. A.  Haveroth$^{2}$,  M. L.  Bittencourt$^{3}$,  J. L.  Boldrini$^{1}$ \\
 \small $^{1}$ Institute of Mathematics,  Statistics and Scientific Computing,  University of Campinas,  SP, 13083-859,  Brazil \\
 \small $^{2}$ Institute of Mathematics and Computer Science,  University of São Paulo at São Carlos,  SP, 13560-970,  Brazil   \\
 \small $^{3}$ School of Mechanical Engineering, University of Campinas, SP, 13083-970, Brazil
}

\maketitle

\begin{abstract}
This paper proposes a thermodynamically consistent phase-field damage model for viscoelastic materials following the strategy developed by Boldrini et al. \cite{boldrini2016non}.
Suitable free-energy and pseudo-potentials of dissipation are developed to build a model leading to a stress-strain relation, under the assumption of finite {strain}, in terms of fractional derivatives.
A novel degradation function, which properly couples stress response and damage evolution for viscoelastic materials, is proposed.
We obtain a set of differential equations that accounts for the evolution of motion, damage, and temperature.
In the present work, for simplicity, this model is numerically solved for isothermal cases by using a semi-implicit/explicit scheme.
Several numerical tests, including fitting with experimental data, show that the developed model accounts appropriately for damage in viscoelastic materials for small and finite strains.
Non-isothermal numerical simulations will be considered in future works.
\end{abstract}


\section{Introduction}
\label{intro}

Interest in damage modeling for {viscoelastic materials has increased greatly in recent years.
Due to their application in the areas of engineering, biology and structural analysis, the appropriate characterization of the viscoelastic behavior} is mandatory to predict component failure, making this theme a very challenging and contemporary research topic.

The earliest contributions for modeling fractures in viscoelastic material date back to the mid-1960s, considering the studies of Knaus et al. \cite{knauss1963rupture,knauss1966time,knauss1969stable,wnuk1970delayed}, Williams \cite{williams1964structural,williams1965initiation} and Schapery \cite{schapery1964application}.
In these studies, crack description is included with the prescription of a critical strain that is established empirically.
Since then, robustness of the models has progressed and many works combining theoretical and computational aspects have been proposed.

Some of the traditional models use the cohesive zone method \cite{tijssens2000modeling,tijssens2000simulation}, which, although widely adopted, presents some difficulties related to the insertion of a sharp interface.
Models based on X-FEM \cite{moes1999finite,yu2011modeling,ozupek2016computational} and peridynamics \cite{madenci2017ordinary} have also sparked researchers' interest, but these {strategies} require considerable reformulation of computational methods or present difficulties to account for nonlinear viscoelasticity \cite{thamburaja2019fracture}.
Thus, continuum approaches have emerged as an alternative to overcome some of these disadvantages \cite{duddu2013nonlocal,nguyen2016large,thamburaja2019fracture}.
Particularly, phase-field models are an interesting concept to deal with material damage due to the ability to describe state changes in a continuum way. 
In other words, they replace the sharp interface by a gradual, but fast, description of the state change induced by the crack propagation; they  may also easily couple thermal and deformation processes by taking into account the influence on stored energy of the material \cite{wang2010phase,shanthraj2016phase}. 
Moreover, the diffuse approximation of discontinuities diminishes the burden of remeshing during crack propagation \cite{nguyen2016phase}. 

One important aspect to be considered is the thermodynamic consistency of the phase-field models.
In this regard, many authors have presented interesting contributions. 
Miehe et al. \cite{miehe2015phase} outlined a framework for phase-field models of crack propagation in elastic solids. 
Fabrizio and co-workers \cite{AmendolaFabrizio,FGM2006} also presented an isothermal model for describing damage and fatigue for non-isothermal cases.
Boldrini et al. \cite{boldrini2016non} presented a general thermodynamically consistent phase-field model for damage and fatigue where the behavior of particular material classes are considered by their corresponding free-energy potentials and pseudo-potentials of dissipation. 
Haveroth et al. \cite{haveroth2020non} included the effects of plasticity in the model presented in \cite{boldrini2016non} and compared several simulated results with experimental data.

The application of phase-field to describe damage in viscoelastic materials was considered by Sch{\"a}nzel \cite{schanzel2015phase} and
Shen et al. \cite{shen2019fracture}, who used traditional rheological models of springs and dashpots to describe viscoelasticity.
These works proposed effective models to predict the material response under loading conditions, although the thermodynamic consistency of these formulations are unclear. 
Furthermore, models based on chains of springs and dashpots frequently require the identification of many material parameters for the constitutive equation.

Despite these difficulties, traditional models based on mechanical
analogies have been used since the mid-19th century \cite{maxwell1867iv,green1946new,casula1992generalized,findley2013creep}, both to describe linear and non-linear viscoelastic behavior \cite{truesdell1955simplest,green1957mechanics,green1959mechanics,christensen1980nonlinear,simo1987fully,holzapfel1996new,le1993three}.
These models are widely used as they are particularly useful for predicting the material response in a purely phenomenological way.
Generalized rheological chains can be used to model a large number of viscoelastic materials; however, as mentioned previously, it can lead to complications for the inverse identification problem since various springs and dashpots may be involved.
Models based on fractional derivatives have emerged as an interesting alternative to describe viscoelastic behavior. 
According to Welch et al. \cite{welch1999}, the use of fractional derivative operators typically demand fewer rheological elements, providing more flexibility to the models. 

Although the relationship between viscoelasticity and fractional derivatives started only after 1930, 
nowadays viscoelasticity analysis is one of the fields with the most extensive applications of fractional derivatives \cite{diethelm1999solution,lion2004payne,schimidt2006,mainardispa2011,caputo2016applications,lazopoulos2016fractional}. {Recent contributions include the works of Jaishankar and McKinley \cite{jaishankar2013power,jaishankar2014fractional}. 
These authors used simple fractional constitutive relationships to characterize the power-law rheological behavior exhibited by viscoelastic interfaces \cite{jaishankar2013power}. 
They also proposed a viscoelastic fractional formulation to describe linear and nonlinear viscoelastic properties of complex liquids and soft solids \cite{jaishankar2014fractional}.
Xu and Jiang \cite{xu2017creep} also used fractional viscoelastic models to characterize creep behavior. 
They provided fitting with experimental data showing the effectiveness of the proposed modeling.

Concerning the thermodynamic aspects of viscoelasticity theories, we mention the classical studies of Coleman \cite{coleman1964thermodynamics,coleman1964thermodynamics2}, which in turn inspired the works of Christensen \cite{christensen1968obtaining, christensen1980nonlinear} and Laws \cite{laws1967thermodynamics}. 
These authors applied the classical Colleman-Noll approach including the called memory effects in the free-energy in order to account for the hereditary effects of viscoelasticity.
Since then, many other researches have presented contributions on this subject \cite{fabrizio1992mathematical,biot1955variational, reese1998theory, schapery1997nonlinear}.
Regarding the fractional viscoelastic theory, 
Lion \cite{lion1997thermodynamics} presented a full derivation of the fractional Zener model from the point of view of thermodynamics, where a free-energy potential is derived and the corresponding potential of dissipation is obtained.
Fabrizio \cite{fabrizio2014fractional} developed a model with fractional derivatives and compared it with the classical Volterra theory. 
This author also proposed a free-energy associated with fractional viscoelasticity.
Alfano and Musto \cite{alfano2017thermodynamic} revisited a fractional model proposed in \cite{musto2015fractional} and presented a thermodynamic derivation that resulted in a linear viscoelastic model.

Recently, fractional viscoelastic models have been coupled to damage models  \cite{sumelka2017hyperelastic,sumelka2020modelling,alfano2017thermodynamic}. 
Krasnobrizha et al. \cite{krasnobrizha2016hysteresis} presented an elastoplastic damage model with fractional derivatives that distinguished the dissipation due to the material damage, plasticity and viscoelasticity. 
This collaborative model is validated for a woven composite with thermoset and thermoplastic matrices.
Tang et al. \cite{tang2018new} proposed a damage model in viscoplastic materials to describe creep in rocks. 
In this work, the fractional derivative is used to describe viscoelastic behavior coupled with a continuum damage approach. 
Good curve fittings of experimental data were obtained, showing that the model can reflect ongoing damage during rock creep.
Caputo and Fabrizio \cite{caputo2015damage} coupled phase-field and fractional derivatives to describe damage in viscoelastic materials. 
They considered the fractional order of the derivative as a phase-field variable which represents the damage evolution, but once more, the thermodynamical consistency of this formulation is unclear.

Although significant progress has been made, many of the constitutive models for damage in viscoelastic materials do not account simultaneously for crack initiation, modeling unloading processes, nonlinear viscoelasticity, or even thermal effects.
Moreover, apart from the works of Tijssens et al. \cite{tijssens2000modeling,tijssens2000simulation}, Schanzel \cite{schanzel2015phase} and Thamburaja et al. \cite{thamburaja2019fracture}, the proposed models are limited to small strain. 

Motivated by this situation, we propose a thermodynamically consistent framework coupling the benefits of phase-fields and fractional derivatives to describe damage in viscoelastic materials.
It results in a model, written in the Lagrangian configuration, that describes the diffuse crack interface by a scalar variable which evolves according to a set of governing equations derived from thermodynamic considerations and leads to automatic crack initiation, that cannot be predicted by discrete fracture models. 
Furthermore, the coupling with the fractional model considers the viscoelastic effects by using less material parameters than
those required in the traditional rheological models.
The model presented here allows finite strain and is based on a mathematical phase-field framework that is similar to the one proposed by Boldrini et al. \cite{boldrini2016non} for linear elastic brittle material, which guarantees thermodynamic consistency and can also include non-isothermal aspects without violating such consistency.

For this, we firstly develop a general phase-field model to account for damage by defining a novel free-energy potential {which includes memory effects}.
The particularization of the model for viscoelastic materials is done by using a suitable free-energy potential and pseudo-potentials of dissipation. 
In particular, the choice of viscoelastic free-energy potential leads to a fractional order differential stress/strain relation. 
This equation includes a degradation function \cite{kuhn2015degradation}, 
which role is to ensure that the part of the driving force associated to the hyperelastic interaction appears in the equation for damage evolution. 
We propose a new degradation function suitable for describing damage in viscoelastic materials.

The model is given by a nonlinear system of {fractional} partial differential equations for the evolution of motion, damage and temperature in materials with viscolastic behavior that is solved by using a {semi-implicit/explicit} finite element scheme. 
Numerical examples include a one-dimensional application of the model to describe the dynamic response of a viscoelastic rod and check the influence of some terms for the stress equation.} 
Afterwards, the two-dimensional extension is used to simulate tensile tests that include loading and unloading processes. 
We also perform an experimental curve fitting for tensile tests with samples of high density polyethylene (HDPE) for small and large strains,
for which the model presents good curve fitting properties in loading process and also good ability to predict the behavior of the specimen for unloading processes. 

\section{Development of the Model}
\label{development_of_the_model}

Consider a body {$\pazocal{B}\subset \mathbb{R}^3$ in the reference configuration with Lagrangian coordinates denoted by $\bm{p}$ and an arbitrary regular subdomain $\pazocal{D}\subset\pazocal{B}$} with boundary $\partial \pazocal{D}$. 
The fundamental state of the body is described by the density of mass $\rho$ which satisfies the principle of mass conservation,  
dynamic variables $\bm{u}$ and $\bm{v}$ representing, respectively, displacement and velocity vector fields, 
and the specific density of internal energy $e$.
The governing equation for $\bm{v}$ is obtained from the principle of virtual power (PVP). 
The first principle of thermodynamics is used for $e$. 

Suppose that $\pazocal{B}$ can develop damage due to strain process. 
We assume that damage evolution can be described by a scalar phase field.
In this case, the phase-field variable $\varphi$ corresponds to the volumetric fraction of damaged material and lies in the interval $[0,1]$; $\varphi=0$ is associated with the undamaged material and $\varphi=1$ with the fractured material. 
In the context of this work, damage is considered a dynamic variable with a corresponding equation obtained from the PVP.

Application of the PVP will require the definition of virtual velocities ${\delta}\bm{v}$ and ${\delta}c$, that are, respectively, admissible macroscopic virtual velocity (the time rate of change of displacement) and admissible microscopic virtual velocity (the time rate of change of dynamic phase-field $\varphi$).

\subsection{Basic Mechanical Aspects}

Following similar arguments developed by {Fr\'emond \cite{fremond2013non} and Boldrini et al. \cite{boldrini2016non}}, the basic governing equations considered here emerge from the mechanical principles which are summarized below.

\begin{enumerate}
	\item The {principle of mass conservation} states that the total mass in a closed system is unaltered by any physical and chemical actions, that is $\dot{\rho} = 0$,
where the dot notation $\dot{(\cdot)}=\frac{\partial}{\partial t}(\cdot)$ corresponds to time derivative.

	\item The {principle of virtual power} (PVP) states the equilibrium of the virtual powers of inertia $\pazocal{P}_a$, internal $\pazocal{P}_i$ and external $\pazocal{P}_e$ loads for any virtual actions (${\delta}\bm{v},{\delta}c)$ and subdomain $\pazocal{D}$ as
\begin{eqnarray}
	\label{pvp}
	\pazocal{P}_a 
	= \pazocal{P}_i 
	+ \pazocal{P}_e,
\end{eqnarray}
\noindent
where
\begin{eqnarray}
	\label{internal_power}
	\pazocal{P}_i = 
	{- \int_{\pazocal{D}} \boldsymbol{P} : \Nabla \left({\delta}\bm{v}\right)\ \diff \pazocal{D} }
	{- \int_{\pazocal{D}}\left[ k \delta {c} + \bm{h} \cdot \Nabla (\delta c)\right] \ \diff \pazocal{D} },\nonumber\\
\end{eqnarray}
	\begin{eqnarray}
  	\label{external_power}
	\pazocal{P}_e  
	= {\int_{\pazocal{D}} \rho \bm{f} . {\delta}\bm{v}\ \diff \pazocal{D}}
	{+ \int_{\partial \pazocal{D}} \left[ \bm{t} . {\delta}\bm{v} + t_h \delta c \right]\ \diff (\partial\pazocal{D})},
\end{eqnarray}
\noindent
and
\begin{eqnarray}
	\label{inertia_power}
	\pazocal{P}_a  
	= \int_{\pazocal{D}} \rho \dot{\bm{v}} . \bm{\delta}\bm{v} \ \diff \pazocal{D}.
\end{eqnarray}
\noindent
{In Eq. \eqref{internal_power},} 
$\bm{P}$ is the first Piola-Kirchhoff stress tensor;
we recall that $\bm{P} = \bm{F} \bm{S}$, where $\bm{F}$ is the gradient deformation tensor, and $\bm{S}$ is the symmetric second Piola--Kirchhoff stress tensor; also,
$k$ is a volume density of energy by unit of $\varphi$ and $\bm{h}$ is an energy flux vector by unit of $\varphi$ \cite{fremond2013non}.
We denote by  $\Nabla (\cdot)$ the gradient operator in the Lagrangian configuration.
The first term on the right-hand side of Eq. \eqref{internal_power} is the classical stress power, while the second term is the power of the interior generalized forces related to the material damage.
In Eq. \eqref{external_power}, $\bm{f}$ is the body force vector  field per unit of mass, $\bm{t}$ is the macroscopic stress vector field and $t_h$ is the superficial density of energy supplied to the material by the flux $\bm{h}$. 
The first integral in Eq. \eqref{external_power} is related to the virtual power of actions at a distance, while the last two terms in the second integral are associated to the virtual powers of the surface loads. 
It is assumed that there are no exterior microscopic actions affecting the damage of the material (e.g., aging or corrosion).

By replacing Eqs. \eqref{internal_power}, \eqref{external_power} and \eqref{inertia_power} into Eq. \eqref{pvp}, using $\bm{\delta}\bm{v}=\bm{0}$ and the fact that $\delta{c}$ is arbitrary, we obtain
\begin{equation*}
	\refstepcounter{equation} \latexlabel{eq_inertia_power_02a}
	\refstepcounter{equation} \latexlabel{eq_inertia_power_02b}
	\divv(\bm{h}) - k = 0 \; \;  \textmd{in} \; \pazocal{D} \quad \textmd{and} \quad
	\bm{h} \cdot \bm{n}_0 = t_h \; \;  \textmd{in } \partial\pazocal{D},
	\tag{\ref{eq_inertia_power_02a}-\ref{eq_inertia_power_02b}}
\end{equation*}
\noindent
where $\bm{n}_0$ is the unit vector normal to the surface area $\partial\pazocal{D}$, and $\divv(\cdot)$ denotes  the divergence operator in the Lagrangian configuration. 
On the other hand, by using ${\delta{c}}=0$, and the fact that $\bm{\delta}\bm{v}$ is arbitrary, we get
\begin{equation*}
	\refstepcounter{equation} \latexlabel{Eq_strongmotionPrep1}
	\refstepcounter{equation} \latexlabel{Eq_strongmotionPrep2}
	\rho\dot{\bm{v} } = \divv(\boldsymbol{P} ) + \rho \bm{f}  \; \; \textmd{in} \; \pazocal{D} 
	\quad \textmd{and}  			\quad  
	\bm{P} \cdot \bm{n}_0 = \bm{t} \; \;  \textmd{in} \; \partial\pazocal{D}.
	\tag{\ref{Eq_strongmotionPrep1}-\ref{Eq_strongmotionPrep2}}
\end{equation*}

	\item The {first principle of thermodynamics} leads to the following equation:
\begin{eqnarray}
\label{dyn_objec}
	\rho \dot{e}  
  	=- \divv(\bm{q})
	+ \rho r 
	+ \bm{S}:\bm{\dot{E}} 
	+ k \dot{\varphi} 
	+ \bm{h} \cdot \Nabla (\dot{\varphi})\quad\text{in } \pazocal{D},
\end{eqnarray}
\noindent
where $e$ is the specific internal energy density; $\bm{q}$ is the heat flux vector field, $r$ is the specific heat source density and $\bm{E}$ is the Green-Lagrange strain tensor.

	\item The entropy inequality is also considered. 
As in Fabrizio et al. \cite{fabrizio2006thermodynamic},  and Boldrini et al. \cite{boldrini2016non}, the second principle of thermodynamics is here expressed in a generalized form of the Clausius-Duhem inequality \cite{truesdell2004non}, whose differential form is given by 
\begin{eqnarray}
   \label{eq_entropy_general}
   \rho \dot{\eta} \geq 
   - \divv(\bm{\Phi}) 
   + \rho \omega 
   \quad \textmd{in } \pazocal{D}.
\end{eqnarray}
\noindent
In this expression, the specific entropy density is denoted by $\eta$;
the general form of the total entropy flux is split as
$\bm{\Phi}= \bm{\Phi}_\theta + \bm{\Phi}_m$;
the term $\bm{\Phi}_\theta={\bm{q}}/{\theta}$ is the classical thermal entropy flux,
and $\theta>0$ is the absolute temperature;
 $\bm{\Phi}_m$ is a possible additional entropy flux due to other microscopic features.
The general form of the total specific entropy production term is also split as $\omega = \omega_\theta + \omega_m$,
where $\omega_\theta = r/\theta$ is the classical specific thermal entropy production and $\omega_m$ is a possible additional specific entropy production term due to other microscopic features.
In the present model these extra terms  may appear due to the damage mechanisms that lead to softening and fracture,  as well as to mechanisms related to memory effects.
For proper modeling, it is required that $\omega_m \geq 0$.
Expressions for $\Phi_m$ and $\omega_m$ will be obtained later on, when we will deal with the expressions for the constitutive relations.

\vspace{0.1cm}
By replacing the Helmholtz specific free-energy
\begin{eqnarray}\label{helmholtz}
	\psi = e - \theta \eta,
\end{eqnarray} 
\noindent
in Eq. \eqref{dyn_objec} and comparing it with inequality \eqref{eq_entropy_general}, we obtain
\begin{eqnarray}
\label{eq_inequality_entropy_complete}
	&& - \rho \left( \dot{\psi} 
	+ \dot{\theta} \eta \right) 
	+ \bm{S}:\bm{\dot{E}} 
	+ k \dot{\varphi} 
	+\bm{h}\cdot \Nabla \dot{\varphi}
	- \dfrac{1}{\theta} \bm{q} \cdot \Nabla (\theta) \nonumber\\
	&&
	+ \theta \divv(\Phi_m)	
	-\theta \rho \omega_m \geq 0 . 
\end{eqnarray}
\noindent
Inequality \eqref{eq_inequality_entropy_complete} must be satisfied for all physical admissible processes to ensure thermodynamic consistency. 
Some details about this aspect are considered in the next section.
\end{enumerate}

\subsection{General Model}
\label{GeneralModel}

We assume that we are dealing with a class of materials  with constitutive relations for $\psi$,
$\bm{S}$, $k$, $\bm{h}$ and $\bm{q}$ that depend on the state variables as follows:
\begin{equation*}
	\refstepcounter{equation} \latexlabel{eq_psi}
\psi :=	\psi(\Theta, \bm{E}, \Ha_t(\bm{E})),
	\tag{\ref{eq_psi}}
\end{equation*}

\begin{equation*}
	\refstepcounter{equation} \latexlabel{eq_split_T}
	\refstepcounter{equation} \latexlabel{eq_split_b}
	\bm{S} :=\bm{S} (\Theta, \bm{E}, \Ha_t(\bm{E}), \dot{\varphi}, \dot{\bm{E}}) ,
	\quad
	k :=  k(\Theta, \bm{E}, \Ha_t(\bm{E}), \dot{\varphi}, \dot{\bm{E}}) ,
	\tag{\ref{eq_split_T}-\ref{eq_split_b}}
\end{equation*}
\begin{equation*}
	\refstepcounter{equation} \latexlabel{eqh0}
	\refstepcounter{equation} \latexlabel{eqq0}
	{\bm{h}} :=  \bm{h} (\Theta, \bm{E}, \Ha_t(\bm{E}), \dot{\varphi}, \dot{\bm{E}}),
	\quad 
	\bm{q} := \bm{q}(\Theta, \bm{E}, \Ha_t(\bm{E}), \dot{\varphi}, \dot{\bm{E}}),
	\tag{\ref{eqh0}-\ref{eqq0}}
\end{equation*}
\noindent
where $\Theta=\{\theta , \varphi , \Nabla \theta , \Nabla \varphi \}$
and
$\Ha_t(\bm{E}):=\Ha(\bm{E})(\bm{p}, t)=\{\bm{E}_s:=\bm{E}(\bm{p}, s) \ \forall\ 0\leq s\leq t\}$
denotes the history
\footnote{
In the present work, we consider only the situation of bodies that are strain free for times $t$ preceding the initial time $t_0$;
that is, we always assume that $\mathbf{E}({t})=\mathbf{0}$ $\forall\ t<t_0$.
Thus, we take the strain history as
\begin{eqnarray*}
\mathcal{H}_t(\mathbf{E})
:=\mathcal{H}(\mathbf{E})(\mathbf{p}, t)
&=&\{\mathbf{E}_s=\mathbf{E}(\mathbf{p},t-s), 0<s <\infty)\}\nonumber\\
&=&\{\mathbf{E}_s=\mathbf{E}(\mathbf{p},s), 0<s <t)\}.
\end{eqnarray*}
\noindent
This definition is a particular case of \cite{fabrizio2014fractional},
and simplifies a bit the technical details.
We could, without too many difficulties, include in our model the complete past history of strain.
}
of the Green-Lagrange strain tensor $\bm{E}$ up to time $t$. 

The specific forms of the constitutive relations for the variables of Eqs. \eqref{eq_split_T}-\eqref{eqq0} will be expressed in terms of the specific free-energy density $\psi$ and the pseudo-potential of dissipation $\psi_d$, 
which are discussed in the next sections.

\subsubsection{General Form of the Free-Energy}
\label{PotentialMemoryEffects}

\begin{figure}[!h]
	\centering
	\includegraphics[scale=0.75]{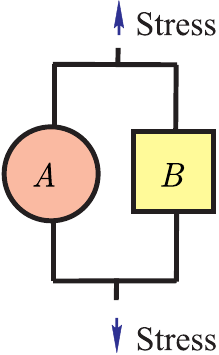} 
	\caption{\small {General rheological arrangement proposed in this work.}}
	\label{Fig_parallel_conection}
\end{figure}

The model here proposed can be compared with a rheological combination of two parts in parallel.
Part A is associated with the local strain effects, and part B is related with the memory strain effects; see Fig. \ref{Fig_parallel_conection}.

\noindent
Particular cases of this situation can be seen in Fig. \ref{Fig_particular_rheological_models}.

Based on the arrangement shown in Fig. \ref{Fig_parallel_conection}, we assumed that is proper to split the total specific free-energy density $\psi$ in two terms\footnote{
Christensen \cite{christensen2012theory}(pg. 265) presents a general free energy depending on the strain and its memory parts.
The total free energy $\varphi$ of Eq. (13) can be considered a particular case  of \cite{christensen2012theory}.
}:
\begin{eqnarray}
	\label{GeneralFormOfFreeEnergy}
	\psi(\Theta, \bm{E}, \Ha_t(\bm{E}))   
	= 
	\psi_c (\Theta , \bm{E}) 
	+ \psi_{m} \left(\Theta , \Ha_t(\bm{E})\right),
\end{eqnarray}

The term $\psi_c$ is the classical space-time pointwise potential related with part A of the rheological model, and $\psi_{m}$ is the potential that accounts for the memory effects in the strain field related with  part B. 
Both $\psi_c$ and $\psi_{m}$ are presented in this section in a general way. 
Specific equations for these functions are defined by the choice of the material to be modeled (see Sec. \ref{Particular_freeEnergy_pseudo} for the case of viscoelasticity).

The general form of the potential $\psi_{m}$ is chosen to be
\begin{eqnarray}
	\label{PotentialMemoryEffectsIncludingDamage}
	\psi_{m}(\varphi , \Ha_t(\bm{E}))
	:=  \displaystyle \frac{G_{m}(\varphi )}{\rho}
	\tilde{\psi}_{m}( \Ha_t(\bm{E})),
\end{eqnarray}
\noindent
where $G_{m}(\varphi) \geq 0$ is a suitable damage degradation function that will be particularized later on (see Sec. \ref{degradation_function}) and the potential $\tilde{\psi}_{m}$ is defined as
\begin{eqnarray}
\label{FirstFormMemoryEffects}
	\tilde{\psi}_{m}
	&:=&\tilde{\psi}_{m}(\Ha_t(\bm{E}))\nonumber\\
	&=&  \frac{1}{\Gamma(1-\alpha)}
	\left[
	\frac{  \mathcal{N}\left( \bm{E}_t, \bm{E}_0 \right) }{t^{\alpha}}
    + \alpha \int_0^{t}\frac{ \mathcal{N}\left( \bm{E}_t,\bm{E}_\tau \right)}{(t-\tau)^{1+\alpha}}\ \diff\tau
	\right],\nonumber\\
\end{eqnarray}
\noindent
where $\bm{E}_{(\cdot)}:=\bm{E}(\bm{p},\scriptstyle(\cdot)\displaystyle)$, $0 < \alpha < 1$, $\Gamma(\cdot)$ is the standard
Gamma function \cite{artin2015gamma} and $\mathcal{N}$ is a suitable continuous function of second-order symmetric tensors satisfying the conditions given in Appendix \ref{AppendixA}.
This appendix also presents the computations to obtain the derivative $\dot{{\psi}}_m$ and a property that is important for ensuring the validity of inequality \eqref{eq_inequality_entropy_complete}.
Appendix \ref{Examples_for_N} gives several examples of possible choices for $ \mathcal{N}$ satisfying the conditions stated in Appendix \ref{AppendixA}.
In particular, Examples 1-3 of Appendix \ref{Examples_for_N} also show particular choices for ${\psi}_m$, which lead to a constitutive stress/strain relation in terms of fractional derivatives.

It is important to emphasize that Eq. \eqref{PotentialMemoryEffectsIncludingDamage} depends only on $\varphi$ and the memory effects on the strain field. 
At the expense of simple, but longer computations, we could easily include in the mathematical model the dependence of $\tilde{\psi}_{m}$ on $\Theta$, and also on its respective memory effects, as suggested in the general form of Eq. \eqref{GeneralFormOfFreeEnergy}.
However, for simplicity of exposition, in this work we consider $\tilde{\psi}_{m}$ as presented in Eq. \eqref{FirstFormMemoryEffects}, depending just on memory effects of the strain field. 

By considering Eq. \eqref{GeneralFormOfFreeEnergy} and using the standard chain rule, we obtain the derivative of $\psi$ as
\begin{eqnarray}
\label{chain_rule}
	\dot{\psi} 
	&=&   \partial_\theta\psi_c \dot{\theta} 
	+ \partial_\varphi \psi_c \dot{\varphi} 
	+ \partial_{\Nabla \theta}\psi_c \dot{\overline{\Nabla \theta }} 
	+ \partial_{\Nabla \varphi}\psi_c \dot{\overline{\Nabla \varphi }}\nonumber\\  
	&&+ \partial_{\bm{E}}\psi_c : \dot{\bm{E}}
	+ \bm{S}_m  : \dot{\bm{E}}    + \frac{G'_{m}}{\rho}\tilde{\psi}_{m} \dot{\varphi}
	- R ,
\end{eqnarray}
\noindent
where $\partial_{(\cdot)}\psi_c$ represents the partial derivative of $\psi_c$ with respect to the subscribed variable, $G'_m:=\frac{\diff G_{m}}{\diff \varphi}$,
\begin{eqnarray}
	\bm{S}_m 
	&=&  \frac{G_{m}}{\rho\Gamma(1-\alpha)}
	\Big[
	\frac{  \partial_{\bm{E}_t} \mathcal{N}( \bm{E}_t, \bm{E}_0 ) }{t^{\alpha}}\nonumber\\
	&&+ \alpha \int_0^{t}\frac{ \partial_{\bm{E}_t} \mathcal{N}( \bm{E}_t, \bm{E}_\tau )}{(t -\tau)^{1+			\alpha}}\ \diff \tau
	\Big],
\end{eqnarray}
\noindent
and
\begin{eqnarray}
	\label{Eq_R}
	R 
	&=& \frac{\alpha G_{m}}{\rho\Gamma(1-\alpha)}
	\Big[   \frac{ \mathcal{N}( \bm{E}_t,\bm{E}_0 )}{t ^{1+\alpha}}\nonumber\\       
	&&+ (1+\alpha) \int_0^{t}\frac{ \mathcal{N}( \bm{E}_t,\bm{E}_\tau )}{(t -\tau)^{2+\alpha}}\ \diff			\tau \Big].
\end{eqnarray}

We observe that $R \geq 0$, due to the property (b) of function $\mathcal{N}$, as shown in Appendix \ref{AppendixA}.

\subsubsection{Expressions for the Constitutive Relations}
\label{ReductionEntropyInequality}

To obtain the general expressions of the constitutive relations, other than the free-energy,
we use the approach found in  Fr\'emond \cite{fremond2013non} and Fabrizio, Giorgi, Morro \cite{FGM2006}.
This approach is related to the Coleman-Noll procedure
and uses the free-energy and the pseudo-potential of dissipation, 
which  is a general way to satisfy the reduced form of the dissipation inequality,
to obtain expressions for the constitutive relations.
It consists of the following five steps.

In the first step,  the general class of materials is restricted by making assumptions on how the constitutive relations depend on the state variables and respective rates. Specific forms, and thus particular cases of materials,  are considered after the basic arguments are stated.
We have already accomplished this step in the introduction of the present Section \ref{GeneralModel}.

In the second step, we split the constitutive relations, other than that of the free-energy, in a 
non-dissipative (reversible) part, which depends only on the state variables, and another 
dissipative (irreversible) part, which may depend on the state variables and some of their rates.

In the third step, we use the entropy inequality (the generalized Clausius-Duhem inequality in our case) and the balance of internal energy to obtain a first form of the dissipation inequality.

In the fourth step, we separate the possibly non-dissipative (reversible) parts and use similar Coleman-Noll arguments to obtain the general expressions of the possibly non-dissipative (reversible) parts of the constitutive relations in terms of the free-energy being considered.
With these results, we are left with the reduced form of the dissipation inequality .

In the fifth and last step, the reduced dissipation inequality leads to the general expressions of the dissipative parts of the constitutive relations and of the extra-thermal entropy flux and production terms as functions of a pseudo-potential of dissipation (the free-energy also appears indirectly).

Once the general theory for the considered class of materials is obtained, specific cases are chosen 
by selecting specific forms for the free-energy and the pseudo-potential of dissipation.

Since we have already done the first step, we go directly to the second step.

We split each one of the constitutive relations in  Eqs. \eqref{eq_split_T}-\eqref{eqq0} in two parts as follows:
\begin{eqnarray}
	&&	\bm{S} = \bm{S}^{(r)} (\Theta, \bm{E}, \Ha_t(\bm{E})) + \bm{S}^{(ir)}  ( \Theta, \bm{E}, \dot{\varphi}, \dot{\bm{E}}) , \label{eq_split_Ts}
	\\
	&&	k = k^{(r)} (\Theta, \bm{E}, \Ha_t(\bm{E})) + k^{(ir)} (\Theta, \bm{E},\dot{\varphi}, \dot{\bm{E}}) , \label{eq_split_bs}
	\\
	&&	{\bm{h_0}} = \bm{h}^{(r)} (\Theta, \bm{E}, \Ha_t(\bm{E})) + \bm{h}^{(ir)} (\Theta, \bm{E},\dot{\varphi}, \dot{\bm{E}}), \label{eqh0s}
	\\
	&&	\bm{q_0} = \bm{q}^{(r)} (\Theta, \bm{E}, \Ha_t(\bm{E})) + \bm{q}^{(ir)} (\Theta, \bm{E},\dot{\varphi}, \dot{\bm{E}}). \label{eqq0s}
\end{eqnarray}
\noindent
As it will be clear from our next computations, 
the terms in the first part of each of these expressions, indicated by the superscript  $(\cdot)^{(r)}$, will be obtained by using the free-energy.
The terms in the second part, indicated by the superscript $(\cdot)^{(ir)}$, will be obtained with the help of a pseudo-potential of dissipation.
The first part terms are expected to give no contribution to the increasing of entropy of the system,
while the second part terms necessarily contribute to increasing of the entropy and are necessarily dissipative terms.
See Remark \ref{RemarkDissipativeNonDissipative} at the end of this subsection for further explanation on these aspects.

Following Fr\'emond \cite{fremond2013non}, we assume $\bm{h}^{(ir)} = \bm{0}$. 
Moreover, the heat flux is purely irreversible, then we take $\bm{q}^{(r)} = \bm{0}$.

By recalling that for any sufficiently smooth field $\chi(\bm{p}, t)$, the time derivative of the Lagrangian gradient operator is given by $ \dot{\overline{\Nabla \chi}}= \Nabla \dot{\chi}$, 
replacing Eq. \eqref{chain_rule} in the entropy inequality \eqref{eq_inequality_entropy_complete} and using the split of the variables into dissipative and nondissipative components given in Eqs. \eqref{eq_split_T}-\eqref{eqq0}, we obtain
\begin{eqnarray}
\label{inequality_en}
	&&-\rho\left(\eta+\partial_\theta \psi_c \right)\dot{\theta}
	+\left( k^{(r)}+k^{(ir)} -\rho\partial_\varphi \psi_c 
	- {G'_{m}} \tilde{\psi}_{m}\right)\dot{\varphi}
	\nonumber
	\\
	&&- \rho\partial_{\Nabla \theta}\psi_c \Nabla \dot{\theta}
	-\left(\rho\partial_{\Nabla \varphi}\psi_c -\bm{h}^{(r)}\right)\Nabla \dot{\varphi}
	\nonumber
	\\
	&&+\left(\bm{S}^{(r)}+\bm{S}^{(ir)}-\rho \partial_{\bm{E}} \psi_c  - \rho\bm{S}_m \right): \dot{\bm{E}}
	\nonumber
	\\ 
	&&-\frac{1}{\theta}\bm{q}^{(ir)}\cdot \Nabla \theta 
	+ \rho R  
	+ \theta \, \divv(\Phi_m) - \theta \rho \omega_m 
	\geq 0.
\end{eqnarray}
\noindent
We now require that the terms in the first three lines of the last inequality do no contribute to the increase of the entropy;
that is, we impose that

\begin{eqnarray}
\label{zero_terms}
	&&-\rho\left(\eta+\partial_\theta \psi_c \right)\dot{\theta}
	+\left( k^{(r)} -\rho\partial_\varphi \psi_c  
	- {G'_{m}} \tilde{\psi}_{m}\right)\dot{\varphi}
	\nonumber 
	\\ 
	&&
	- \rho\partial_{\Nabla \theta}\psi_c \Nabla \dot{\theta}
	-\left(\rho\partial_{\Nabla \varphi}\psi_c -\bm{h}^{(r)}\right)\Nabla \dot{\varphi}
	\nonumber 
	\\ 
	&&
	+\left(\bm{S}^{(r)}-\rho \partial_{\bm{E}} \psi_c  - \rho\bm{S}_m \right):\dot{\bm{E}} 
	=0.
\end{eqnarray}
\noindent
Since $\dot{\theta}$, $\dot{\varphi}$, $\Nabla \dot{\varphi}$, $\dot{\bm{E}}$ and $\Nabla \dot{\theta}$ of Eq. \eqref{zero_terms} are arbitrary and independent, the classical Coleman-Noll approach leads to
\begin{equation*}
	\refstepcounter{equation} \latexlabel{eta0}
	\refstepcounter{equation} \latexlabel{b0}
	\eta = -\partial_\theta \psi_c ,\quad
	k^{(r)}=\rho\partial_\varphi\psi_c + G'_{m} \tilde{\psi}_{m},
	\tag{\ref{eta0}-\ref{b0}}
\end{equation*}
\begin{equation*}
	\refstepcounter{equation} \latexlabel{partialnabla}
	\refstepcounter{equation} \latexlabel{h0r}
	\partial_{\Nabla \theta}\psi_c = 0,\quad
	\bm{h}^{(r)}=\bm{h}=\rho\partial_{\Nabla \varphi}\psi_c,
	\tag{\ref{partialnabla}-\ref{h0r}}
\end{equation*}
\noindent
and
\begin{eqnarray}
	\label{Sclass}
	\bm{S}^{(r)}=\rho \partial_{\bm{E}}\psi_c  + \rho\bm{S}_m.
\end{eqnarray}

By replacing the above relations in inequality \eqref{inequality_en}, 
and remembering that $\theta>0$, we obtain
\begin{equation}
\label{DissipIneqMod1}
\begin{array}{r}
\displaystyle
	\frac{k^{(ir)}}{\theta}\dot{\varphi}
	+\frac{1}{{\theta}}\bm{S}^{(ir)}:\bm{\dot{E}} 
	- \frac{\bm{q}^{(ir)}}{\theta^2}\cdot \Nabla \theta
	\\
\displaystyle
	  + \divv(\Phi_m)  + \frac{\rho}{\theta} R - \rho \omega_m  
	  \geq 0.
\end{array}
\end{equation}

Since we want to develop the simplest possible theory, we reduce this last expression by taking 
the additional flux of entropy and the specific entropy production due to to microscopic features other than thermal ones respectively as 
$\Phi_m = 0$ {and} $\omega_m = R/{\theta}$.
It means that $R$ is related to the possible extra specific entropy production due to microscopic features other than the thermal ones; in the present model they are related to the damage mechanisms that lead to softening and fracture.

By using these results in \eqref{DissipIneqMod1}, we then are left with the following reduced form of the dissipation inequality:
\begin{eqnarray}
\label{inequality_3}
	\frac{k^{(ir)}}{\theta}\dot{\varphi}
	+\frac{1}{\theta}\bm{S}^{(ir)}:\bm{\dot{E}} 
	- \frac{\bm{q}^{(ir)}}{\theta^2}\cdot \Nabla \theta
	\geq 0.
\end{eqnarray}

\subsubsection{General Pseudo-Potential of Dissipation}

To ensure inequality (\ref{inequality_3}), it is enough to take the coefficients
$\displaystyle \frac{k^{(ir)}}{\theta}$,
$\displaystyle \frac{1}{\theta}\bm{S}^{(ir)}$
and 
 $\displaystyle - \frac{\bm{q}^{(ir)}}{\theta^2}$
respectively as the derivatives of the pseudo-potential $\psi_{d}$ with respect to 
$\dot{\varphi}$, $\dot{\bm{E}}$ and $\Nabla \theta$.


This pseudo-potential of dissipation is a nonnegative  functional that in the present situation has the general expression
\begin{eqnarray}
	\psi_{d}:=\psi_{d}(\dot{\varphi}, \dot{\bm{E}},\Nabla \theta, \tilde{\Theta})\geq 0 ,
\end{eqnarray}
\noindent
and satisfies $\psi_{d}(0, \bm{0}, \bm{0}, \tilde{\Theta})=0$ where $\tilde{\Theta} = \{ \theta, \varphi, \nabla \varphi, \bm{E} \}$.
Moreover, it must be continuous and convex with respect to the independent variables $\dot{\varphi}$, $\dot{\bm{E}}$ and $\Nabla \theta$. 

We then obtain
\begin{equation*}
	\refstepcounter{equation} \latexlabel{irrev1}
	\refstepcounter{equation} \latexlabel{irrev2}
	\begin{array}{c}
	\displaystyle
	k^{(ir)} = \theta\partial_{\dot{\varphi}} \psi_{d}, 
	\quad
	\bm{S}^{(ir)} =  \theta\partial_{\dot{\bm{E}}} \psi_{d}, 
	\end{array}
	\tag{\ref{irrev1}-\ref{irrev2}}
\end{equation*}
\noindent
and
\begin{eqnarray}
\label{irrev3}
	\bm{q}^{(ir)} = - \theta^2 \partial_{\Nabla \theta} \psi_{d}.
\end{eqnarray}
\noindent

If $\psi_{d}$ is non-differentiable, then we must work with subdifferentials.

\begin{remark}
\label{RemarkDissipativeNonDissipative}
Dissipation is related to the increase of the entropy.
Thus, terms appearing in constitutive relations for a certain material are said to be either dissipative or non-dissipative according  to they respectively do or do not contribute to the increase of the entropy.

The procedure we used in this subsection to obtain the general expressions for the constitutive relations in terms of the free-energy and the pseudo-potential of dissipation is related to what is called the Coleman-Noll procedure.
In this procedure, usually one expects to obtain the non-dissipative terms of the constitutive relations in the first step of the arguments by using the free-energy;
then one is left with what is called the reduced form of the dissipation inequality,
which, by using suitable pseudo-potentials of dissipation, give necessarily dissipative terms.

However, the terms obtained in the first step of this procedure are only guaranteed to be non-dissipative if there are no extra (non-thermal) non-negative sources of entropy associated to them.
In fact, when there are extra (non-thermal) non-negative sources of entropy,
the entropy can increase due to the presence of those sources.
In this case, those terms of the constitutive relations obtained in the first step of the procedure
that contribute to those extra non-negative sources of entropy,
although derived from the free-energy, 
are in fact dissipative.

This is exactly the situation of our model, where there is an extra source of entropy $\omega_m = R/{\theta}$,
where $R \geq 0$ is given in Eq. \eqref{Eq_R} and depends on the memory terms of the free-energy.
Thus, the terms in the constitutive relations related to the memory part of the free-energy are dissipative even though they are not derived from a pseudo-potential of dissipation.

In particular, it can also be seen in the particularized model proposed in the following section, where the viscoelasticity is modeled by using fractional fractional derivative. 
In this case, an intuitive argument to understand the dissipative contribution of terms with fractional derivatives is the following:
due to the nature of a fractional derivative element, with interpolates between the behavior of a spring (non-dissipative) and a dashpot (dissipative), it always includes some dissipation.
\end{remark}

\subsection{Viscoelastic Model}
\label{Particular_freeEnergy_pseudo}

The framework derived until this moment is general in the sense that the appropriate choices for the free-energy potential, $\psi$, and the pseudo-potential of dissipation, $\psi_d$, consider several classes of materials.
Now, we want to particularize this model for viscoelastic materials. 
For that purpose, we start describing the particular form of the pseudo-potential of dissipation that we will consider in this work.

\subsubsection{A Particular Pseudo-Potential of Dissipation}

A possible choice for $\psi_{d}$, satisfying the conditions described in the previous subsection, is the following:
\begin{eqnarray}
\label{pseudo_d1}
	\psi_{d}(\dot{\varphi}, \dot{\bm{E}},\Nabla \theta, \tilde{ \Theta})
	&=&\frac{\tilde{\lambda}(\tilde{\Theta})}{2}|\dot{\varphi}|^2
	+\frac{\tilde{b}(\tilde{\Theta})}{2}|\dot{\bm{E}}|^2
	\nonumber \\ &&
	+\frac{\tilde{c}(\tilde{\Theta})}{2}\Nabla \theta \cdot\bm{C}^{-1}\Nabla \theta.
\end{eqnarray} 
\noindent
The inverse of parameter $\tilde{\lambda}$ corresponds to the rate of change in damage $\varphi$ \cite{boldrini2016non}, and we take it as
\begin{eqnarray}
\label{um_sobre_lambda}
	\frac{1}{\tilde{\lambda}}=\frac{{c}_\lambda}{(1+\tilde{\delta}-\varphi)^\zeta}>0,
\end{eqnarray}
\noindent
where ${c}_\lambda$ and $\zeta$ are positive material parameters and $\tilde{\delta}$ is a small perturbation to avoid singularity when $\varphi=1$. 
Moreover, {$\tilde{b}\geq 0$} and $\tilde{c}\geq0$ correspond to the viscous damping and the heat conductivity of the material, respectively \cite{boldrini2016non}.

By considering this pseudo-potential of dissipation and Eqs. \eqref{irrev1}-\eqref{irrev3}, we obtain
\begin{equation*}
	\refstepcounter{equation} \latexlabel{b0ir}
	\refstepcounter{equation} \latexlabel{sclass}
	k^{(ir)} = \theta\tilde{\lambda}(\tilde{\Theta})\dot{\varphi},
	\quad
	\bm{S}^{(ir)} = \theta \tilde{b}(\tilde{\Theta})\bm{\dot{E}} ,
	\tag{\ref{b0ir}-\ref{sclass}}
\end{equation*}
\noindent
and
\begin{eqnarray}
	\label{q0ir}
	\bm{q}^{(ir)}= - \theta^2 \tilde{c}(\tilde{\Theta})\bm{C}^{-1}\Nabla \theta,
\end{eqnarray}
where $\bm{C}=2\bm{E}+\bm{I}$ is the right Cauchy-Green strain tensor. 

\begin{remark}
\label{RemarkPseudoPotentialExplanation}
The form of last term in expression \eqref{pseudo_d1} deserves an explanation.

The term related to the gradient of the temperature in the pseudo-potential of dissipation gives rise to a diffusion of temperature term in the energy equation.
To be physically correct,  this term must be related to some variant of Fick's law,
which requires that  the actual diffusive heat flow must be in a direction determined by the  Eulerian gradient of the temperature.
That is, the actual diffusive heat flow directions must be expressed in terms of $\nabla_{\mathbf{x}}(\theta)$,
where the subscript $(\cdot)_{\mathbf{x}}$ indicates the operation done in Eulerian coordinates.

For instance, to obtain the simple case of isotropic thermal diffusion, with thermal diffusion coefficient $\tilde{c}(\tilde{\Theta}) >0$,
the corresponding term in the pseudo-potential of potential, when written in Eulerian coordinates,
must be have a term of form
\[
\psi_{d}^{(\theta)}:= \frac{1}{2} \tilde{c}(\tilde{\Theta}) |\nabla_{\mathbf{x}}\theta|^2.
\]

We observe that many articles dealing with this issue use a Lagrangian version of this energy obtained by just replacing the Eulerian gradient by the Lagrangian gradient.
However, we think that this replacement is not quite correct because the resulting Lagrangian form does not generate the correct diffusion term in the energy equation.

The last term in \eqref{pseudo_d1} gives the physically correct diffusion of temperature because it is simply the previous $\psi_{d}^{(\theta)}$ written in Lagrangian coordinates; this change of variable is done because the theoretical framework in the present work is based on Lagrangian description.

In fact, we have $\nabla_{\mathbf{x}} \theta = \mathbf{F}^{-t} \nabla \theta$,
where the last gradient is with respect to Lagrangian coordinates.
By substituting this in the previous expression of $\psi_{d}^{(\theta)}$ and performing some simple computations, we obtain
\[
\psi_{d}^{(\theta)} = \frac{\tilde{c}(\tilde{\Theta})}{2} |\nabla_{\mathbf{x}} \theta|^2
 = \frac{\tilde{c}(\tilde{\Theta})}{2}\nabla \theta \cdot\mathbf{C}^{-1} \nabla \theta ,
\]

\noindent
which is exactly the last term appearing in \eqref{pseudo_d1}.
\end{remark}


\subsubsection{A Particular Free-Energy Potential}
\label{sec_free_energy}

The local free-energy density $\psi_c$ is decomposed in three parts related to the hyperelastic deformation, represented by $\psi_h$, purely thermal effects, given by $\psi_\theta$, and damage contributions, considered in $\pazocal{I}$. 
Therefore, the volumetric density of the part of the free-energy independent of memory effects is given by
\begin{eqnarray}\label{free_energy}
	\rho\psi_c (\theta, \varphi, \Nabla \varphi, \bm{E})
	&=&
	G_h (\varphi) \psi_h (\bm{E})+\psi_\theta(\theta)
	+\pazocal{I}(\varphi,\Nabla \varphi, \bm{E}),\nonumber
\end{eqnarray}
\noindent
where $G_{h}(\varphi) \geq 0$ is a suitable damage degradation function of the hyperelastic part of the free-energy which will be particularized later on (see Sec. \ref{degradation_function}). 
Note that, as in the case of  Eq. \eqref{PotentialMemoryEffectsIncludingDamage},
 the volumetric density of elastic energy $\psi_h$ is multiplied by this degradation function.

The hyperelastic energy { density} for a compressible Neo-Hookean material is given by \cite{bonet2008nonlinear}
\begin{eqnarray}\label{hyperelastic}
	\psi_h
	&=&\frac{\mu}{2}[\mbox{tr}(\bm{C})-3]-\mu\ln[\det(\bm{C})]^{\frac{1}{2}}
	+\frac{\lambda}{2}[\ln(\det(\bm{C}))^{\frac{1}{2}}]^2,\nonumber\\
\end{eqnarray}
\noindent
where $\mu$ and $\lambda$ are the Lam\'e material parameters.
The nonlinear elastic behavior of Eq. \eqref{hyperelastic} requires the consideration of finite strain in this model.

The thermal part of the free-energy is assumed to be \cite{fremond2013non}:
\begin{eqnarray}
	\psi_\theta=c_v\theta\ln\theta,
\end{eqnarray}
\noindent
where $c_v$ is the heat capacity. 

The damage contribution is given by
\begin{eqnarray}
\label{FreeEnergyParcelDamageVariable}
	\pazocal{I}
	&=& g_c\left(\frac{\gamma}{2}\Nabla \varphi \cdot \bm{C}^{-1}\Nabla \varphi +\frac{1}{\gamma}H(\varphi)			\right).
\end{eqnarray}
\noindent
The Griffith fracture energy $g_c$ is assumed positive and constant; $H(\varphi) = \frac{\varphi^2}{2}$ is the potential for $\varphi \in [0,1]$ \cite{boldrini2016non}.
The parameter $\gamma>0$ is related with the width of the fractured layers. 
According to \citep{haveroth2020non}, smaller values for $\gamma$ must lead to a less diffuse crack path,
and sufficiently small values of $\gamma$
lead to sharp cracks \cite{bourdin2008variational}.
In addition, $\delta $ can also be related to crack propagation speed, being faster for larger values of $\delta$ (see \citep{haveroth2020non}).

Finally, concerning the free-energy with memory effects $\psi_m$, we assume the definition of Eq. \eqref{PotentialMemoryEffectsIncludingDamage} with $\tilde{\psi}_m$ given by
\begin{eqnarray}
\label{pseudod2_text}
	\tilde{\psi}_{m}  
	&=&\frac{\kappa}{\rho}\bigg[\frac{\left(\bm{E}_t-\bm{E}_0\right):\mathbfcal{{A}} : \left(\bm{E}_t-			\bm{E}_0\right)}{t^{\alpha}}\nonumber\\
	&&+ \alpha \int_0^{t}\frac{\left(\bm{E}_t-\bm{E}_\tau\right): \mathbfcal{{A}} : \left(\bm{E}_t-				\bm{E}_\tau\right)}{(t - \tau)^{1+\alpha}} \diff \tau	\bigg],
\end{eqnarray}
\noindent
as suggested in Eq. \eqref{pseudod2}. 
Here, $\kappa={1}/(2\Gamma(1-\alpha))$ and the specific form of the fourth order symmetric tensor $\mathbfcal{{A}}$ will be described in Sec. \ref{ResultsDiscussion}.
Equation \eqref{pseudod2_text} leads to
\begin{eqnarray}
\label{SmUSED}
	\bm{S}_m 
	&=& \displaystyle \frac{G_{m}}{\rho} \left[  \mathbfcal{{A}} :{_0\mathrm{D}_t}^\alpha (\bm{E}_t) 			\right.+\kappa\frac{\left(\bm{E}_t-\bm{E}_0\right): \partial_{\bm{E}_t} \mathbfcal{{A}}: 					\left(\bm{E}_t-\bm{E}_\tau\right)}{t^\alpha}\nonumber\\
	&& \displaystyle \left. + \alpha\kappa \int_0^{t}\frac{\left(\bm{E}_t-\bm{E}_\tau\right): 					\partial_{\bm{E_t}} \mathbfcal{{A}} : \left(\bm{E}_t-\bm{E}_\tau\right)}{(t - \tau)^{1+\alpha}}\ 			\diff \tau \right]
\end{eqnarray}
\noindent
where ${_0\mathrm{D}_t}^\alpha \bm{E}_t$ is the Caputo fractional derivative of $\bm{E}$ (see Appendix \ref{AppendixA} for details and comments on how to obtain Eq. \eqref{SmUSED}).

\begin{remark}
As for the case of the term commented in Remark \ref{RemarkPseudoPotentialExplanation},
the first term of Eq. \eqref{FreeEnergyParcelDamageVariable} also deserves some explanation.

We recall that a standard physical assumption 
in phase field models is that part of the energy associated to the damage process accumulates in 
transition layers of the damage variable.
This assumption brings a contribution to the free-energy that depends on the gradient of $\varphi$.

It is important, however, to understand that, to correspond to the physical situation, 
such parcel of the free-energy must depend directly on the gradient of $\varphi$ in the actual (deformed) configuration of the body, not directly on the gradient with respect to the reference configuration,
which, in principle may be arbitrary.

The simplest case to be considered is that in which that parcel of the free-energy has an expression depending on the Eulerian gradient of $\varphi$ as follows:
\[
	\psi^{(\varphi)}
	= g_c \frac{\gamma}{2} |\nabla_{\mathbf{x}} \varphi |^2.
\]

Many articles dealing with this issue use a Lagrangian version of this parcel by just replacing the Eulerian gradient by the Lagrangian gradient, and this leads to a simple form of the equation for the damage variable $\varphi$ where a simpler damage diffusion term $\Delta \varphi$ appears.

Sometimes this term $\Delta \varphi$ appears directly in the equation for the damage variable, without mention of the corresponding pseudo-potential of dissipation, because it is taught simple as an artificial, but convenient, smoothing approximation.

But, as in our previous remark, we think that these approaches is not quite physically correct 
because, as we previously said,  energy can in fact accumulate in transitions layers and must be considered with energy parcels similar to the previous $\psi^{(\varphi)}$.

Another way to see the difficulty of considering a free-energy with a term depending on the square of the norm of the damage gradient with respect to the variables in the reference configuration is following.
Let us  consider two different reference configurations, related by a change of variables that is not a simple rotation.
Then, the use of a free-energy with a term depending on the square of the norm of the damage gradient with respect to the variables in the reference configuration leads in both cases to a diffusion term $\Delta \varphi$, obviously with derivatives in terms of variables associated to each reference configuration.
However, these diffusion terms are not correctly related by the changing of variables relating these two reference configurations,
and they will correspond to different patterns of damage spreading in the actual deformed configuration, and thus to different physical predictions.
This aspect may be not usual.
We also observe that this difficulty does not appear in the bulk part of the free-energy, that  is, the part depending on the pointwise values of $\varphi$;
the parts of the  driving-forces associated to the bulk free-energy in  two different reference configurations would be correctly related.

Thus, we think that the correct way to proceed in a Lagrangian framework is to rewrite $\psi^{(\varphi)}$  
in terms of the Lagrangian gradient by using the relation between Eulerian and Lagrangian gradients.

By using 
$\nabla_{\mathbf{x}} \varphi = \mathbf{F}^{-t} \nabla \varphi$
in $\psi^{(\varphi)}$, after some simple computations, we obtain the Lagrangian expression:
\[
	\psi^{(\varphi)}
	= g_c \frac{\gamma}{2} |\nabla_{\mathbf{x}} \varphi |^2
	= g_c \frac{\gamma}{2} \nabla \varphi \cdot \mathbf{C}^{-1} \nabla \varphi  ,
\]

\noindent
which is exactly the expression of the first term appearing in Eq. \eqref{FreeEnergyParcelDamageVariable}.
\end{remark}

Finally, our derivation can easily modified to obtain the usual diffusion model, replacing 
$\bm{C}$ by $\bm{I}$ in  \ref{FreeEnergyParcelDamageVariable}). Simpler computation and 
usual linear dissipation term for damage spreading are achieved.


{\subsubsection{Viscoelastic Stress}}
\label{viscoelastic_stress}

\noindent
Under the conditions of the previous subsection, Eqs. \eqref{Sclass}, \eqref{irrev2} and \eqref{SmUSED} lead to the following second Piola-Kirchhoff stress tensor:
\begin{eqnarray}
	\label{stress}
	\bm{S}
	=\bm{S}^{(r)}+\bm{S}^{(ir)} 
	={\rho\partial_{\bm{E}}\psi_c}^{} 
	+ \rho\bm{S}_m
	+{\theta}\partial_{\bm{\dot{E}}}\psi_{d}.
\end{eqnarray}

Taking into account Eq \eqref{stress}, the expression for $\bm{S}_m$ given by Eq. \eqref{SmUSED}, the local free-energy $\psi_c$ of Eq. \eqref{free_energy} and the pseudo-potential of dissipation $\psi_{d}$ of Eq. \eqref{pseudo_d1}, the complete expression for the second Piola-Kirchhoff stress tensor is given by
\begin{eqnarray}
	\label{final_second_piola}
	\bm{S}
	&=&G_h\left(\mu(\bm{I}-\bm{C}^{-1})+\lambda \ln(\det \bm{C})^{\frac{1}{2}}\bm{C}^{-1}\right)
	+\theta\tilde{b}\bm{\dot{E}}\nonumber\\
	&&-g_c\gamma(\bm{C}^{-1}\Nabla (\varphi))\otimes (\bm{C}^{-1}\Nabla (\varphi))
	+{G_{m}} \left[  \mathbfcal{{A}} :{_0\mathrm{D}_t}^\alpha (\bm{E}_t) \right.\nonumber\\
	&&+\kappa\frac{\left(\bm{E}_t-\bm{E}_\tau\right): \partial_{\bm{E}_t} \mathbfcal{{A}} : \left(\bm{E}		_t-\bm{E}_\tau\right)}{t^\alpha}\nonumber\\
	&&+ \alpha\kappa\int_0^{t}\frac{\left(\bm{E}_t-\bm{E}_\tau\right): \partial_{\bm{E}_t} 						\mathbfcal{{A}}: \left(\bm{E}_t-\bm{E}_\tau\right)}{(t - \tau)^{1+\alpha}}\ \diff\tau \bigg]  .
\end{eqnarray}
For the sake of simplicity, we assumed that the hyperelastic and the memory parts of the stress tensor degenerate in the same way, that is, $G(\varphi):= G_h(\varphi) = G_m(\varphi)$. Section \ref{degradation_function} describes the degradation function $G(\varphi)$.

\begin{figure*}[!ht]
	\centering
	\begin{subfigure}{.32\textwidth}
	\centering
		\includegraphics[scale=0.75]{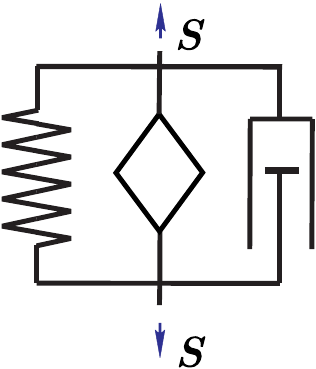} 
		\caption{\small {General model.}}
		\label{general_model}
	\end{subfigure}
	\begin{subfigure}{.32\textwidth}
		\centering
		\includegraphics[scale=0.75]{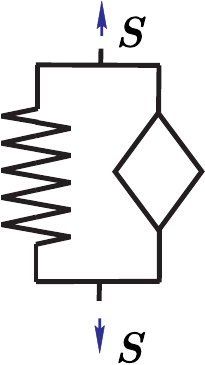} 
		\caption{\small {Modified fractional Kelvin-Voigt model.}}
		\label{kelvin_frac}
	\end{subfigure}
	\begin{subfigure}{.32\textwidth}
	\centering
		\includegraphics[scale=0.75]{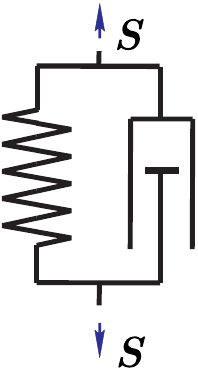} 
		\caption{\small {Modified Kelvin-Voigt model.}}
		\label{trad_kelvin}
	\end{subfigure}
	\caption{\small {Possible interpretations of our model in the one-dimensional case for $\mathbfcal{{A}}$ constant. The spring represents the hyperelastic contribution of the Neo-Hookean material and the dashpot gives the viscous damping. The fractional rheological element is represented by the rhombus.}}
		\label{Fig_particular_rheological_models}
\end{figure*}

If {the tensor} $\mathbfcal{{A}}$ is constant, then the last two terms of Eq. \eqref{final_second_piola} are null.
In this case, the one-dimensional version of the model can be represented by the rheological mechanism of Fig. \ref{general_model}.
The spring represents the hyperelastic contribution of the Neo-Hookean material and the dashpot represents the viscous dissipative damping given by the term $\theta\tilde{b}\bm{\dot{E}}$.
The fractional rheological element, represented by the rhombus, is called spring-pot \cite{koeller1984applications} and results in a nondissipative viscoelastic counterpart, whose behavior is governed by $\mathbfcal{{A}}$ and $\alpha$.
Here, the degradation function $G(\varphi)$ indicates that the spring and the spring-pot include damage effects.
Additionally, if $\mathbfcal{{A}}$ is constant and $\tilde{b}=0$, then we obtain the modified fractional Kelvin-Voigt\footnote{In this work, we refer to modified fractional Kelvin-Voigt when the spring represents a Neo-Hookean spring to account for hyperelasticity. 
If the spring represents the traditional linear elastic material, then we refer to the traditional Kelvin-Voigt model.} model of Fig. \ref{kelvin_frac}.
On the other hand, if $\mathbfcal{{A}}$ is constant and viscoelastic effects due to the fractional component are not considered, we recover the modified Kelvin-Voigt model of Fig. \ref{trad_kelvin} {which includes thermal effects}.
In other words, by using the appropriate simplification, Eq. \eqref{stress} can describe several material behaviors. 
Additionally, the last two terms in Eq. \eqref{final_second_piola} come from the consideration of memory effects in $\bm{S}_m$. 
In fact, these terms do not contribute very much to the evaluation of stress and can be disregarded in several cases. 
A complete study on this subject is presented in Sec. \ref{Evaluation_of_stress_terms}.


Even for the one-dimensional case, it is important to emphasize that for small strain, the Neo-Hookean spring becomes the traditional linear elastic spring. 
In this case, if $\mathbfcal{{A}}$ and $\theta$ are constants and no damage is considered, the model described in this work recovers the usual fractional Kelvin-Voigt model, largely discussed in the literature \cite{lewandowski2010identification,xu2015equivalent,farno2018comparison}. 
Section \ref{Displacement_Rod} presents an example where this simplification is considered.
In fact, for that case, Eq. \eqref{final_second_piola} is simplified for a widely known equation, for which the thermodynamics were addressed by Lion \cite{lion1997thermodynamics}.
A free-energy potential was even derived with physical justification and the corresponding mechanical dissipation potential {was} obtained.
{In the present work}, we generalize the hypothesis including the possibility of nonlinear dependence of $\mathbfcal{{A}}$ on $\bm{E}$, large strain and thermal effects. 

\vspace{0.2cm}
\subsubsection{Degradation Function}
\label{degradation_function}
\begin{figure*}[!h]
\centering
	\begin{subfigure}{.32\textwidth}
	\centering
	\vspace{0.2cm}
		\includegraphics[scale=0.45]{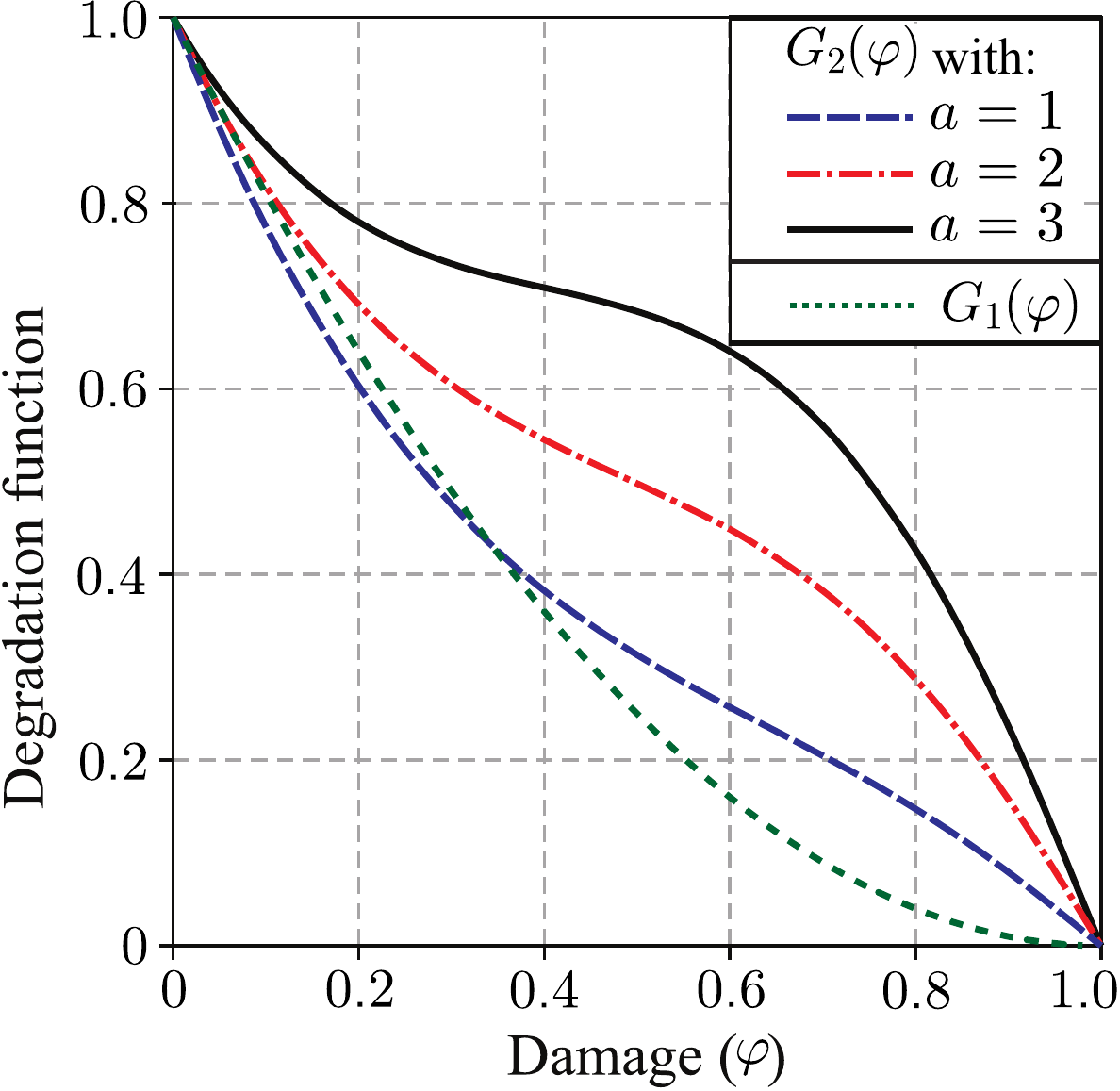} 
		\vspace{-0.2cm} 
		\caption{\small{Degradation function $G_2$ for some values of $a$ and $b=c=1$.}}
		\label{newdegrada}
	\end{subfigure}
	\begin{subfigure}{.32\textwidth}
	\centering
	\vspace{0.2cm}
		\includegraphics[scale=0.45]{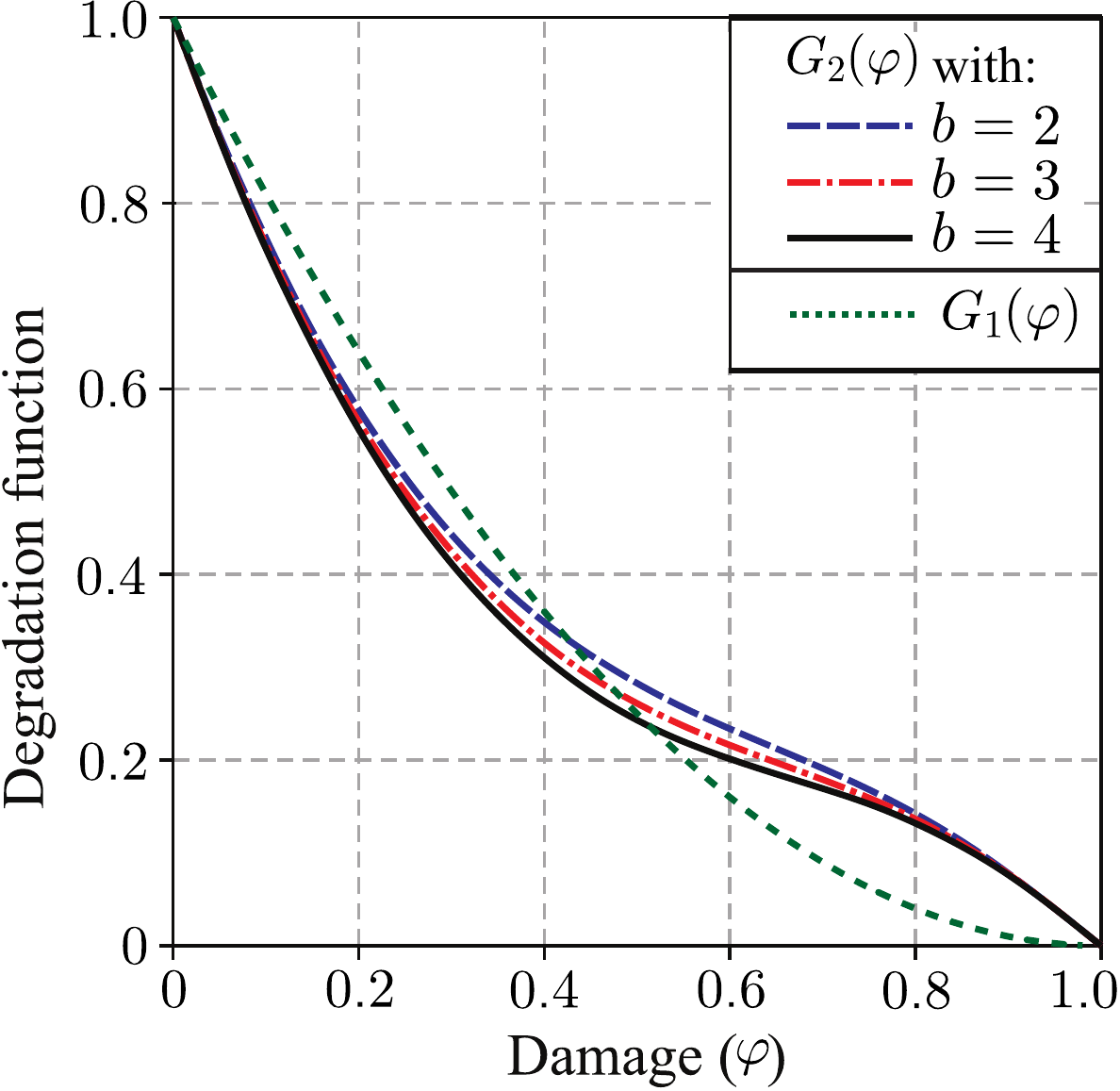}
		\vspace{-0.1cm} 
		\caption{\small{Degradation function $G_2$ for some values of $b$ and $a=c=1$.}}
		\label{newdegradb}
	\end{subfigure}
	\begin{subfigure}{.32\textwidth}
	\centering
	\vspace{0.2cm}
		\includegraphics[scale=0.45]{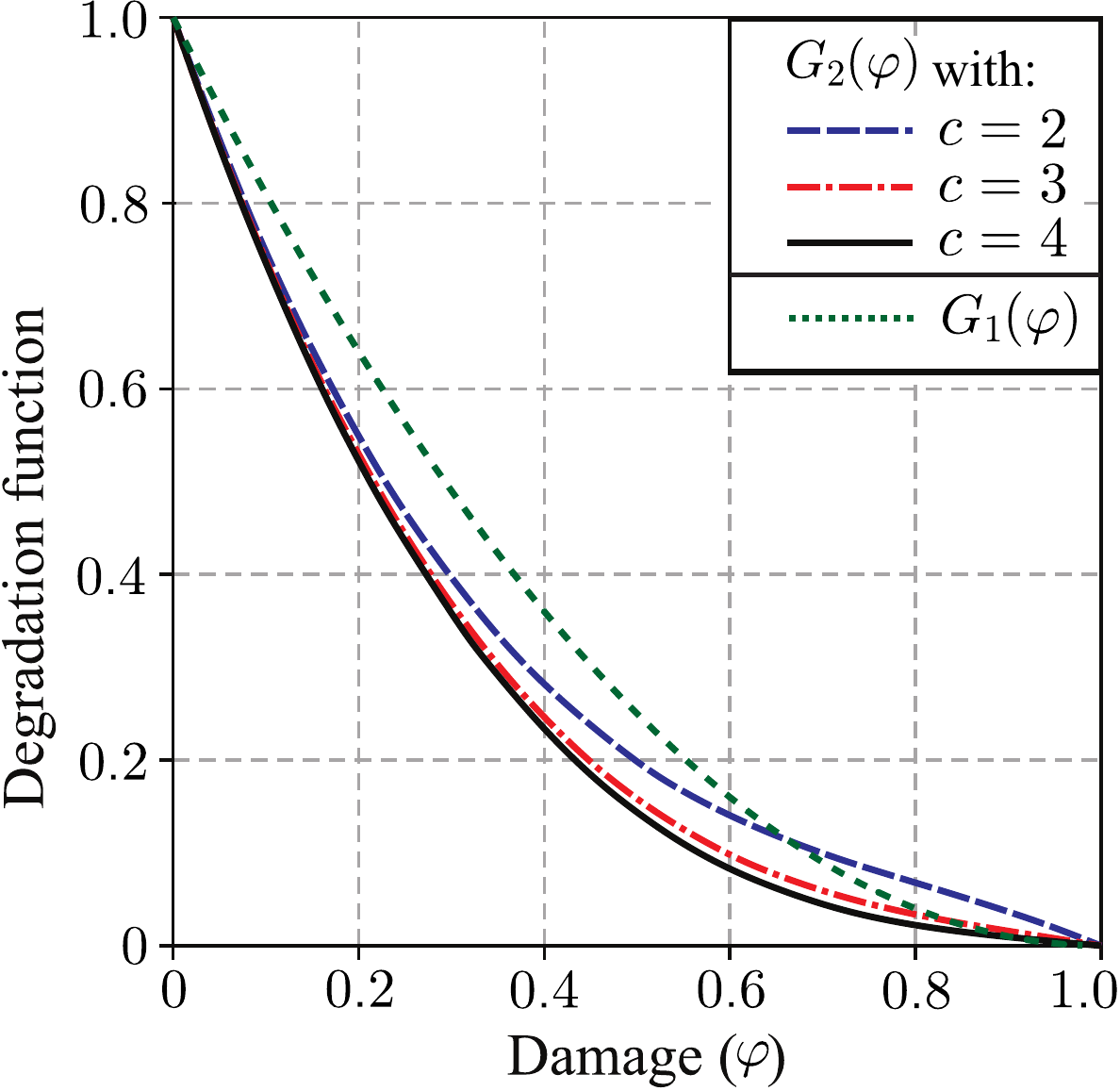} 
		\vspace{-0.1cm} 
		\caption{\small{Degradation function $G_2$ for some values of $c$ and $a=b=1$.}}
		\label{newdegradc}
	\end{subfigure}
	\caption{\small{Degradation functions.}}
\label{quadraticdegrad}
\end{figure*}

The degradation function $G(\varphi)$ couples the damage to the material properties {and} models the change in stiffness between the undamaged and fractured states. 
The material response of the damage is mainly controlled by the degradation function which must satisfy the conditions \cite{MieheWH2010}:
\begin{equation*}
	\refstepcounter{equation} \latexlabel{degrad1}
	\refstepcounter{equation} \latexlabel{degrad2}
	G(\varphi)>0,\quad
	G(\varphi)= \begin{cases} 
	1 & \mbox{if } \varphi=0,\\ 
	0 & \mbox{if } \varphi=1, 
	\end{cases}
	\tag{\ref{degrad1}-\ref{degrad2}}
\end{equation*}
\noindent
and
\begin{eqnarray}
	\label{degrad3}
	G'(1)=0.
\end{eqnarray}
The condition expressed by Eq. \eqref{degrad3} ensures that the part of the driving force associated to the hyperelastic interaction appears in the equation for the evolution of damage $\varphi$. 
There are many proposals for this function \cite{kuhn2015degradation,borden2016phase,haveroth2020non}, which in turn depend on the material. 
Firsty, we follow  Miehe et al. \cite{MieheHW2010} and use the quadratic function:
\begin{eqnarray}
\label{functionG1}
	G(\varphi):=G_1=(1-\varphi)^2.
\end{eqnarray}
\noindent
Indeed, this expression is one of the most frequently found in the literature to describe the degradation function for crack modeling, but it yields a significant degradation of stiffness, as can be seen in Fig. \ref{quadraticdegrad}. 
This is not desirable when modeling viscoelastic materials which present a different behavior in fracture.

\begin{remark}
\label{Remark_viscoelastic_damage}
The damage process due to loading for viscoelastic materials generally occurs in two steps: slippage of the chains and chain separation \cite{daniels1989polymers}.
Differently from the case of metals, most of the viscoelastic materials are made up of long molecular chains \cite{anderson1994fracture}, and the speed of the slippage process depend on the considered material. 

Let us consider the damage evolution for the case of polymers (the standard example of viscoelastic material).
Under tensile stress, a rather fast chain separation process occurs; these chain separation may lead to nucleation and coalescence of voids  and to a certain amount of stiffness degradation;
then slippage along the chains occurs in a process leading to small decrease in stiffness;
as  the slippage of chains increases, the localized stress level also increases;
when the stress on a small chain segment  is larger than the bound strength can sustain,
chain breaking occurs, leading again to void nucleation and coalescence of voids,
resulting again in a rather fast stiffness degradation.

This process may lead to coalescence of voids and evolve until  fracture
\cite{daniels1989polymers, anderson1994fracture, kuksenko2013fracture}.
Additionally, according to Christensen \cite{christensen2012theory}, viscoelastic mateials can undergo 
local failure, even at moderate strain levels. Local instabilities, as slippage, often dominate the subsequent 
behavior. 

In other words, the damage process for a viscoelastic material can not be restricted to the nucleation and coalescence of voids, once the slippage is an important part of the process.

\end{remark}

Taking into account the damage process explained in Remark \ref{Remark_viscoelastic_damage},

we propose a new degradation function:
\begin{eqnarray}
\label{functionG2}
	G(\varphi):=G_2=(1-\varphi)^3+\frac{a{\varphi}^{d}(1-\varphi)^{d}}{1+b(\varphi-c)^2},
\end{eqnarray}
where parameters $a,b$ and $c$ are assumed to be positive and $d=1.05$. 
Exponent $d$ could be a variable considered as an additional parameter to be determined. 
However, in this work we consider this value fixed to obtain the desirable behavior in $G_2$. 
Figures \ref{newdegrada}-\ref{newdegradc} show the effect of variation of $a,b$ and $c$ in this function.
In fact, $G_2$ imposes different transitions in the damage process, when compared with the quadratic function $G_1$.

The change in the concavity of $G_2$ creates a region where the damage grows slowly. 
We consider that this region is related with the slippage of the chains, as described in Remark \ref{Remark_viscoelastic_damage}.
The time between the slippage and the fracture can be adjusted by the parameters $a$, $b$ and $c$, used to define $G_2$ (see Eq. 51).
It allows more flexibility for the modeling of different materials, once the variation of $a$, $b$ and $c$ can change the region of $G_2$ related with the slippage accordingly.
As we see in Remark \ref{Remark_viscoelastic_damage}, this behavior is best suited for viscoelastic materials because it agrees with the micro-structural evolution in strain processes. 
A test comparing functions $G_1$ and $G_2$ is presented in Sec. \ref{largestrains}.

\subsubsection{Final Governing Equations for the Viscoelastic Model}\label{final_governing}

Considering the aspects discussed previously, the final governing equations for the evolution of motion, damage and temperature in a body with viscoelastic behavior can be resumed as follows.

\begin{enumerate}
	\item  The equation of motion is given by the balance of linear momentum of Eq. \eqref{Eq_strongmotionPrep1}; that is,
\begin{equation} 
	\label{strongmotion}
	\rho\dot{\bm{v} } = \divv(\boldsymbol{P} ) + \rho \bm{f} .
\end{equation}

\noindent
We recall that in this equation $\bm{v}$ is the velocity field and that $\bm{P} = \bm{F} \bm{S}$, with
the constitutive relation stress/strain for $\bm{S}$ given by Eq. \eqref{final_second_piola}.

	\item In order to obtain the equation for the damage evolution, we replace $k$ and $\bm{h}$ in Eq. \eqref{eq_inertia_power_02a} using Eqs.  \eqref{eq_split_b}, \eqref{b0}, \eqref{h0r} and \eqref{b0ir}. Then
\begin{eqnarray}
	\label{damage_intermediate}
	\theta \partial_{\dot{\varphi}}\psi_{d}
	&=&\divv(\rho\partial_{\Nabla \varphi}\psi_c) 
	- \rho\partial_\varphi\psi_c 
	- G'_{m}(\varphi ) \tilde{\psi}_{m}.
\end{eqnarray}
\noindent
The above equation is written in terms of the pseudo-potential of dissipation $\psi_{d}$, the free-energy $\psi_c$ and the free-energy with memory effects $\psi_m$, given by Eqs. \eqref{pseudo_d1}, \eqref{free_energy} and \eqref{pseudod2_text}, respectively.
Replacing the derivatives of Eq. \eqref{damage_intermediate} by the corresponding expressions and recalling that $ G' = G'_{h} = G'_{m}$, we obtain
\begin{eqnarray}
	\label{final_damage}
	\dot{\varphi}&=&\frac{1}{\tilde{\lambda}\theta}\divv \left( g_c\tilde{\gamma}\bm{C}^{-1}\Nabla \varphi \right) 
	-\frac{g_c H'(\varphi)}{\gamma  \tilde{\lambda}\theta} \nonumber\\  
	&&- \frac{G'}{\tilde{\lambda}\theta} \left(\psi_h +\tilde{\psi}_{m}\right).
\end{eqnarray}

	\item The expression for the temperature evolution is obtained by considering Eq. \eqref{dyn_objec}. 
We replace the expressions for $\bm{q}$, ${k}_0$ and $\bm{h}$ given in Eqs. \eqref{q0ir}, \eqref{b0}, \eqref{h0r} and \eqref{b0ir} to obtain 
\begin{eqnarray}
	\rho\dot{e}_0
	&=& -\divv(\theta\partial_{\Nabla \theta} \psi_{d})+\rho r
	+\bm{S}:	\bm{\dot{E}}\nonumber \\ &&
	+\left(\rho\partial_\varphi\psi_c + \partial_{\dot{\varphi}} \psi_{d} + {G'} \tilde{\psi}_{m}\right)		\dot{\varphi}\nonumber \\ &&
	+\rho\partial_{\Nabla \varphi}\psi_c \cdot \Nabla \dot{\varphi}.
	\label{temperature_eq_1}
\end{eqnarray}
	By taking this expression and using the Helmholtz specific free-energy of Eq. \eqref{helmholtz},  we have
\begin{eqnarray}\label{final_temperature}
	-\rho \theta \partial^2_\theta\psi \dot{\theta}
	&=&\divv (\theta \partial_{\Nabla (\theta)}\psi_{d})
	+\rho_0r_0\nonumber \\ 
	&&+ \left(\rho\partial_\theta\partial_\varphi \psi_c
	+\partial_{\dot{\varphi}}\psi_d
	\right)\dot{\varphi}\nonumber \\ 
	&&+ \rho\theta \partial_{\Nabla \varphi} \partial_\theta\psi_c : \Nabla \dot{\varphi}
	+\rho R\nonumber \\ &&
	+\left(\rho\theta\partial_\theta \partial_{\bm{E}}\psi+ \partial_{\dot{\bm{E}}}\psi_{d}
	\right):\dot{\bm{E}}.
\end{eqnarray}
	The free-energy {functional} $\psi$ and the pseudo-potential of dissipation $\psi_d$, defined by Eqs. \eqref{free_energy} and \eqref{pseudo_d1}, respectively, lead to the final governing equation for the temperature evolution
\begin{eqnarray}\label{final_temperature2}
	\dot{\theta}&=&
	\frac{1}{c_v}\divv\left(\theta\tilde{c}\bm{C}^{-1}\Nabla \theta\right)
	+\frac{\tilde{\lambda}}{c_v}|\dot{\varphi}|^2
	+\frac{\rho r}{c_v}\nonumber \\ &&		
	+\theta\frac{\tilde{b}}{c_v}|\bm{\dot{E}}|^2 + \frac{\rho R}{c_v}.
\end{eqnarray}

\end{enumerate}

Equations \eqref{final_second_piola}, \eqref{strongmotion}, \eqref{final_damage} and \eqref{final_temperature2} constitute a nonlinear system of differential equations with fractional derivatives and memory terms. 
The numerical approximation used to solve this system is presented in the next section.

\begin{remark}
\label{Remark_irrevers_damage}
Note that the governing equations, as presented in this section, does not necessarily ensure the irreversibility of the damage;
that is it does not guarantee that $\dot{\varphi} \geq 0$.
This means that the model, as proposed up to now, allows the possibility of healing,
a behavior that can be, in fact, found in some real materials \cite{hayes2007self,li2010thermomechanical, li2014self}.
 
 However damage irreversibility can easily be incorporated in the model.
 From the theoretical point of view,
 as mentioned by Miehe \cite{MieheWH2010} and Boldrini \cite{boldrini2016non}, a possibility for this 
is to adapt the model by adding a multi-valued convex functional $U(\dot{\varphi})$ to the pseudo-potential given in Eq. (38), $U(z)= [0,+\infty)$ if $z<0$ and $U(z)=0$ if $z\geq 0$.
By working with subdifferentials (see \cite{boldrini2016non} for details), we obtain $\partial_\varphi U(\dot{\varphi})= (-\infty,0]$ if $\dot{\varphi}<0$ and $\partial_\varphi U(\dot{\varphi})= 0$  if $\dot{\varphi}\geq 0$. 
It adds a Lagrange multiplier to the stress tensor and, consequently, for the damage equation ensuring that $\dot{\varphi}\geq 0$ \cite{haveroth2020non}.

This seems an intricate approach  to guarantee damage irreversibility, 
but, due to simple form of the required restriction, $\dot{\varphi} \geq 0$, 
which simple means that $\varphi$ cannot decrease in time,
from the practical point of view
irreversibility can be easily implemented in numerical simulations.

It is enough to use a kind of predict-corrector procedure as follows.
Assuming known the state values at time step $n$,
we use the evolution equation without the additional term in the pseudo-potential to predict the damage values
at time step $n+1$,
obtaining a predict value  $\varphi^*_{n+1}$ for the damage variable. 
Next, for each node of the mesh, we compare $\varphi^*_{n+1}$ to $\varphi_{n}$:
if at that node $\varphi^*_{n+1} \geq \varphi_{n}$ then we take $\varphi_{n+1} = \varphi^*_{n+1}$;
otherwisel
$\varphi_{n+1} = \varphi_{n}$.
This in fact is a practical implementation of the above theoretical approach,
which does not require to compute the Lagrange multipliers due to the simple form of the required constraint.

In the simulations to be presented later on, we guarantee damage irreversibility by imposing the numerical constraint as just described. See also Sec. \ref{Sec_Eq_damage}.


There are other possibilities in the literature to impose damage irreversibility:  one could use either the history of elastic energy as in \cite{miehe2010phase} or the penalty criteria as in \cite{haveroth2020non}.
\end{remark}

{\section{Numerical Approximation}}

\noindent
This section presents the numerical approximation used to solve the nonlinear system of equations summarized in Sec. \ref{final_governing}. 
The global method concerns on the application of a semi-implicit/{explicit} time integration scheme coupled with the Newton-Raphson method \cite{haveroth2018comparison}. 

The semi-implicit/{explicit} scheme consists in solving
each equation of the system individually by using a suitable implicit time integration method, resulting in a significant computational economy when compared with usual coupled methods to solve nonlinear systems of equations.

Since the temperature was fixed for the numerical simulations presented in this work,
it is enough to explain how the damage variable and displacement are evolved from a time-step to the next.
This is done as follows: we solve the damage equation to obtain the updated damage variable by using the backward Euler method for time discretization and the Newton-Raphson procedure to handle the nonlinearities; 
at this stage we use as input the known displacement of the previous time-step. 
Next, the just updated damage variable is kept fixed and used as input in the equation of motion, which is solved by the standard Newmark method also combined with the Newton-Raphson procedure. 
It results in the updated displacement, velocity and acceleration. 


The time interval $[0,T]$ is divided into $k$ intervals considering the time step  $\Delta t =t_n-t_{n-1}$, with $n=1,\cdots,k+1$. 
The time discretization is indicated in the updated variables for the time $t_{n+1}$ by using the subscript $(\cdot)_{n+1}$.

The spatial discretization is done for two-dimensional finite element meshes. 
We consider a division of the domain $\pazocal{D}$ into $m$ elements 
$\pazocal{D}_q$ with $q=1,\cdots,m$, {where} $\pazocal{D}=\cup_{q=1}^m \pazocal{D}_q$ and $\pazocal{D}_i\cap \pazocal{D}_j=\emptyset$, for $i \neq  j$.
The approximation of the vector $\bm{z}$ and scalar $z$ fields in each $q$-th element are written as a superposition of the local nodal basis function $N_i$ (Lagrange polynomials) as
\begin{equation*}
	\refstepcounter{equation} \latexlabel{bw}
	\refstepcounter{equation} \latexlabel{w}
	{z}^q\simeq\bm{N}\tilde{{z}}^q
	\quad \textmd{and} \quad
	\bm{z}^q\simeq\hat{\bm{N}}\tilde{\bm{z}}^q,
	\tag{\ref{bw}-\ref{w}}
\end{equation*}
\noindent
with $i=1,\cdots,\vartheta$ and $\vartheta$ the number of element nodes.
Here, the tilde symbol represents the nodal values of the field of interest.
The matrices ${\bm{N}}$ and $\hat{\bm{N}}$ {are} given by
\begin{eqnarray}
	\label{noriginal}
	\bm{N}=\begin{bmatrix}
	N_1 & N_2 & \cdots & N_\vartheta
	\end{bmatrix},
\end{eqnarray}
\noindent
and
\begin{eqnarray}
	\label{form_functions}
	\hat{\bm{N}}=\begin{bmatrix}
	N_1 & 0   & N_2 & 0   & \cdots & N_\vartheta & 0 \\ 
	0   & N_1 & 0   & N_2 & \cdots & 0     & N_\vartheta
	\end{bmatrix} .
\end{eqnarray}
\noindent
The approximation for the gradient operator is given by the global derivatives of the interpolation functions to the $x$ and $y$ directions and organized as
\begin{eqnarray}
	\label{B1}
	{\bm{B}}=\begin{bmatrix}
	N_{1,x} & N_{2,x} & \cdots & N_{\vartheta,x} \\ 
	N_{1,y} & N_{1,y} & \cdots & N_{\vartheta,y}
	\end{bmatrix},
\end{eqnarray}
\noindent
and
\begin{eqnarray}
	\label{derivative_form_funtions}
	\hat{\bm{B}}=\begin{bmatrix}
	N_{1,x} & 0       & N_{2,x}   & 0       &\cdots     & N_{\vartheta,x} & 0 \\
	N_{1,y} & 0       & N_{2,y}   & 0       &\cdots     & N_{\vartheta,y} & 0\\  
	0       & N_{1,x} &  0        & N_{2,y} &\cdots     & 0         &  N_{\vartheta,x}\\
	0       & N_{1,y} &  0        & N_{2,x} &\cdots     &0          & N_{\vartheta,y}
	\end{bmatrix}.
\end{eqnarray}

Details concerning the linearization and numerical discretization  for each equation of the governing system are presented below.

\vspace{0.5cm}
\subsection{Equation of Motion}
\label{discretization_motion}

The evolution of motion is given by  the balance of linear momentum in Eq. \eqref{strongmotion}. 
Using  finite elements, we must work with its corresponding weak form, 
which can be obtained in the standard way by taking the inner product of Eq. \eqref{strongmotion} with any virtual velocity  $\bm{\delta}\bm{v}$ and doing integration by parts:
\begin{eqnarray} 
\label{motion0}
	\int_{\pazocal{D}} \rho\dot{\bm{v}}\cdot {\delta}\bm{v} \ \diff \pazocal{D} 
	&=& -{ \int_{\pazocal{D}} ( \boldsymbol{F}  \bm{S} ) : \Nabla \left({\delta}\bm{v} \right)\ \diff \pazocal{D}}
	+{\int_{\pazocal{D}} \rho \bm{f} \cdot {\delta}\bm{v}\ \diff \pazocal{D}}\nonumber \\ 
	&&	{+ \int_{\partial \pazocal{D}} \bm{t} \cdot {\delta}\bm{v}\ \diff \left(\partial\pazocal{D}\right)}.
\end{eqnarray}		
\noindent
Considering the symmetry of $\bm{S}$ we have
\begin{eqnarray}
 \boldsymbol{F}  \bm{S} : \Nabla \left({\delta}\bm{v} \right)
&=&   \bm{S}  : \bm{F}^t\Nabla \left({\delta}\bm{v} \right)
= \bm{S}  : \bm{F}^t\delta \dot{ \bm{F}}\nonumber\\
&=&  \bm{S}  : \frac{1}{2}\left(\bm{F}^t\delta \dot{ \bm{F}} + \delta \dot{ \bm{F}}^t \bm{F} \right ),
\end{eqnarray}
\noindent
where $\delta \dot{ \bm{F}}=\nabla \left(\bm{\delta}\bm{v}\right)$.
Using the last equation, it is possible to rewrite Eq. \eqref{motion0} in terms of $\bm{S}$ and ${\delta}\bm{\dot{E}}$ as
\begin{eqnarray} 
	\label{motion1}
	\int_{\pazocal{D}} \rho\dot{\bm{v}}\cdot {\delta}\bm{v} \ \diff \pazocal{D} 
	&=& -{ \int_{\pazocal{D}} \bm{S} : {\delta}\bm{\dot{E}} \ \diff \pazocal{D}}
	+{\int_{\pazocal{D}} \rho \bm{f} \cdot {\delta}\bm{v}\ \diff \pazocal{D}}\nonumber \\ 
	&&	{+ \int_{\partial \pazocal{D}} \bm{t} \cdot {\delta}\bm{v}\ \diff \left(\partial\pazocal{D}\right)},
\end{eqnarray}
\noindent		
with
\begin{eqnarray}
\label{dotEvirtual}
	{\bm{\delta}\dot{\bm{E}}(\bm{u}	)}
	= \frac{1}{2} \left[\boldsymbol{F}^t (\bm{u})  \delta \dot{ \bm{F}}  
	+ \delta \dot{ \bm{F}}^t \bm{F} (\bm{u}) \right],
\end{eqnarray}
\noindent
where ${\bm{\delta}\dot{\bm{E}}(\bm{u})}$ is the time rate of the Green Lagrange virtual strain tensor.

The numerical solution of Eq. \eqref{motion1} involves three steps: application of the Newmark method for time discretization; linearization of the nonlinear terms; and application of the finite element method for space discretization.

\subsubsection{Newmark Method for the Equation of Motion}

In the Newmark method (see \cite{lindfield2012numerical}, p. 266), the acceleration at time step $n+1$ is updated using the following relation:
\begin{eqnarray}
	\label{newm1}
	\dot{\bm{v}}_{n+1}=\ddot{\bm{u}}_{n+1} = a_1(\bm{u}_{n+1}-\bm{u}_n)-a_2 \dot{\bm{u}}_n-a_3\ddot{\bm{u}}_n.
\end{eqnarray}
\noindent
Constants $a_i$, with $i=1,\cdots, 3$, are given by
\begin{equation*}
	\refstepcounter{equation} \latexlabel{constnewm1}
	\refstepcounter{equation} \latexlabel{constnewm2}
	\refstepcounter{equation} \latexlabel{constnewm3}
	a_1=\frac{1}{\tilde{\beta}\Delta t^2}, \quad
	a_2=\frac{1}{\tilde{\beta}\Delta t}, \quad 
	a_3=\frac{1-2\tilde{\beta}}{2\tilde{\beta}}, 
	\tag{\ref{constnewm1}-\ref{constnewm3}}
\end{equation*}
\noindent
where $\tilde{\beta}$ is the Newmark constant. 

By replacing Eq. \eqref{newm1} in Eq. \eqref{motion1}, we obtain the following expression for the residue of the time discretization of the equation of motion:
\begin{eqnarray}
	\label{motion_time_disc}
	\bm{R}_{n+1}&=&\int_{\pazocal{D}} \left(a_1(\bm{u}_{n+1}-\bm{u}_n)-a_2 \dot{\bm{u}}_n
	-a_3\ddot{\bm{u}}_n\right)\cdot \delta\bm{{v}}\ \diff \pazocal{D}
	\nonumber
	\\
	&&  +{\frac{1}{\rho} \int_{\pazocal{D}} {\bm{S}(\bm{E}(\bm{u}_{n+1}))} : {\bm{\delta}\dot{\bm{E}}(\bm{u}	_{n+1})}}\ \diff \pazocal{D} 
	\nonumber
	\\
	&&-{\int_{\pazocal{D}} \bm{f}_{n+1} . \delta{\bm{v}}\ \diff \pazocal{D}}
	{-\frac{1}{\rho}\int_{\partial \pazocal{D}} \bm{t}_{n+1} . \delta{\bm{v}}\ \diff \left(\partial				\pazocal{D}\right)},
\end{eqnarray} 
\noindent
where for shortness of notation we did not make explicit the dependence on the other variables except the mechanical ones.

\subsubsection{Linearization of the Weak Form}

Let us evaluate the directional derivative of $\bm{R}_{n+1}$ with respect to the displacement at $\bm{u}_{n+1}$ in the direction of the displacement increment 
$\bm{w}_{n+1} = \Delta \bm{u}_{n+1}$, denoted by  $\bm{D}_{\bm{w}_{n+1}} \bm{R}_{n+1} := \nabla \bm{R}_{n+1} \cdot \bm{w}_{n+1}$ (see \cite{bhatti2006advanced}).
 
The directional derivative of the first term in the right hand side of Eq. \eqref{motion_time_disc} is given by
\begin{eqnarray}
	\label{linearized_dynamical}
	&&\bm{D}_{\bm{w}_{n+1}} \left(\int_{\pazocal{D}} \left(a_1(\bm{u}_{n+1}-\bm{u}_n)-a_2 	\dot{\bm{u}}_n-a_3\ddot{\bm{u}}_n\right)\cdot {\delta}\bm{v}\ \diff \pazocal{D}\right)\nonumber\\
	&&= a_1\int_{\pazocal{D}}\bm{w}_{n+1} \cdot \bm{\delta}\bm{v}\ \diff \pazocal{D}.
\end{eqnarray}
\noindent

Next, the body force term is considered, where $\bm{f}_{n+1}$ represents the updated body forces in the initial configuration that are not affected by the displacement. 
Thus, the directional derivative  with respect to displacement variation $\bm{w}_{n+1}$ is zero. 
Similarly, we assume that the surface loads $\bm{t}$ do not depend on the deformation, then its directional derivative is also zero.
 
The second term on the right hand side of Eq. \eqref{motion_time_disc} has two sources of nonlinearities from the displacement field. 
In fact, the stress tensor $\bm{S}$ depends on strain, which in turn depends on the displacements;
furthermore, from expression (\ref{dotEvirtual}), we see that the time rate of the  Green-Lagrange virtual strain tensor is also a function of displacement.

For the computations that follow, we observe that the directional derivative of $\boldsymbol{F}_{n+1}$ with respect to the displacement along an increment of displacement $\bm{w}_{n+1}$ is given by
\begin{eqnarray}
\label{DerivativeF}
	\bm{D}_{\bm{w}_{n+1}} (\boldsymbol{F}_{n+1}) 
	= \Nabla (\bm{w}_{n+1}).
\end{eqnarray}
\noindent
Moreover,  for the sake of simplicity, we denote
\[\bm{S}(\bm{E}(\bm{u}_{n+1})) := \bm{S}_{n+1},\qquad
\bm{F}(\bm{u}_{n+1}):=\bm{F}_{n+1},\] and
\[\bm{\delta}\dot{\bm{E}}_{n+1}  := \bm{\delta}\dot{\bm{E}}(\bm{u}	_{n+1}). \]
Then, the directional derivative of the second term in the right hand side of Eq. \eqref{motion_time_disc} with respect of displacement along a displacement increment $\bm{w}_{n+1}$ is obtained by using the product rule differentiation as follows:
\begin{eqnarray}\label{internal_power_linearization}
	&&\bm{D}_{\bm{w}_{n+1}}\left(\frac{1}{\rho} \int_{\pazocal{D}} 
	\bm{S}_{n+1}  : {\delta}\dot{\bm{E}}_{n+1} \ \diff  \pazocal{D}\right)\nonumber\\
	&&   
	= \frac{1}{\rho} \int_{\pazocal{D}}  \bm{D}_{\bm{w}_{n+1}} (\bm{S}_{n+1})  : \bm{\delta}\dot{\bm{E}}_{n+1} \ \diff  \pazocal{D} \nonumber\\
    && +  \frac{1}{\rho} \int_{\pazocal{D}}  
    \bm{S}_{n+1}  : \bm{D}_{\bm{w}_{n+1}} ( {\bm{\delta}\dot{\bm{E}}}_{n+1}  )  \ \diff  \pazocal{D}	.   
\end{eqnarray}

Now, from (\ref{dotEvirtual}) and (\ref{DerivativeF}) we have
\begin{eqnarray}
	&&\bm{D}_{\bm{w}_{n+1}} ( {\delta}\dot{\bm{E}}_{n+1} ) \nonumber\\
	&&=	\frac{1}{2} \left[ \Nabla (\bm{w}_{n+1})^t   \delta\dot{\bm{F}}_{n+1}  
	+ \delta\dot{\bm{F}}_{n+1}^t  \Nabla (\bm{w}_{n+1})  \right].
\end{eqnarray}

For the linearization of the second Piola-Kirchhoff stress tensor, the chain rule is used for differentiation to obtain
\begin{eqnarray}
	\label{DerivativeS}
	\bm{D}_{\bm{w}_{n+1}}\left(\bm{S}_{n+1}\right) 
	&=& \left(\frac{\partial \bm{S}}{\partial \bm{E}}\right)_{n+1} : \bm{D}_{\bm{w}_{n+1}} (\bm{E}_{n			+1})\nonumber\\
	&=&\mathbfcal{{D}}_{n+1}:\bm{D}_{\bm{w}_{n+1}} (\bm{E}_{n+1}),
\end{eqnarray}
\noindent
where $\frac{\partial \bm{S}}{\partial \bm{E}}\bigg|_{n+1}=\mathbfcal{{D}}_{n+1}$ is the fourth order, symmetric and positive-definite tangent stiffness tensor. 
The components of $\mathbfcal{{D}}_{n+1}$ are obtained by differentiating the constitutive relation given in Eq. \eqref{final_second_piola} to $\bm{E}_{n+1}$. 
In this work, $\mathbfcal{{D}}_{n+1}$ is calculated by using complex derivatives as explained in Sec. \ref{tensorA}.

The linearization of ${\bm{E}}_{n+1}$ can be obtained from 
\begin{eqnarray}
	\label{E_forH}
	\bm{E}
	=\frac{1}{2}\left[\nabla(\bm{u})^t\nabla(\bm{u})
	+\nabla(\bm{u})^t
	+\nabla(\bm{u})
	\right]
\end{eqnarray}
\noindent
and written as
\begin{eqnarray}
	\bm{D}_{\bm{w}_{n+1}} ( {\bm{E}}_{n+1} )  
	=\frac{1}{2} ( \Nabla (\bm{w}_{n+1})^t   {\bm{F}}_{n+1}  + {\bm{F}}_{n+1}^t  \Nabla (\bm{w}_{n+1})  ).
\end{eqnarray}

Now, we consider the symmetry of $\mathbfcal{{D}}$ and $\bm{S}$;
after some algebraic manipulations, we obtain
\begin{eqnarray}
	\label{eq1}
	\boldsymbol{S} : \frac{1}{2}\left[
	\delta \dot{\bm{F}}^t\nabla(\bm{w}) +  \nabla(\bm{w})\delta\dot{\bm{F}}	\right]
	=\delta \dot{\bm{F}} : \nabla(\bm{w})\bm{S},
\end{eqnarray}
\noindent
and
\begin{eqnarray}
	\label{eq2}
	&&\frac{1}{2}\left[\delta \dot{\bm{F}}^t \bm{F} + \bm{F}^t\delta \dot{\bm{F}}\right]
	:\mathbfcal{{D}}
	:\frac{1}{2} ( \Nabla (\bm{w})^t   {\bm{F}}  + {\bm{F}}^t  \Nabla (\bm{w})  )\nonumber\\
	&&= \bm{F}^t \delta \dot{\bm{F}}
	:\mathbfcal{D}
	:\bm{F}^t \nabla(\bm{w}).
\end{eqnarray}
\noindent
Therefore,
\begin{eqnarray}
	\label{internal_power_linearization2}
	&&\bm{D}_{\bm{w}_{n+1}}\left(\frac{1}{\rho} \int_{\pazocal{D}} {\bm{S}_{n+1}}:						
	\delta	\dot{\bm{E}}_{n+1}\ \diff \pazocal{D}\right)\nonumber\\
	&&=\frac{1}{\rho}\int_{\pazocal{D}}\delta \dot{\bm{F}}^t_{n+1} : \nabla(\bm{w}_{n+1})\bm{S}_{n+1}\ 	\diff {\pazocal{D}}\nonumber\\
	&& +\frac{1}{\rho}\int_{\pazocal{D}}\bm{F}^t_{n+1} \delta \dot{\bm{F}}_{n+1}
	:\mathbfcal{D}
	:\bm{F}^t_{n+1} \nabla(\bm{w}_{n+1})\ \diff {\pazocal{D}}.
\end{eqnarray}

The final linearized form of Eq. \eqref{motion_time_disc} is obtained from Eqs. \eqref{linearized_dynamical} and \eqref{internal_power_linearization2} and given by
\begin{eqnarray}
\label{linearized_internal_power}
	&&\bm{D}_{\bm{w}_{n+1}} \left(\bm{R}_{n+1}\right) 
	= a_1\int_{\pazocal{D}}\bm{w}_{n+1} \cdot\hat{\bm{v}} \diff \pazocal{D}\nonumber\\
	&& +\frac{1}{\rho}\int_{\pazocal{D}}\delta \dot{\bm{F}}^t _{n+1}: \nabla(\bm{w}_{n+1})\bm{S}_{n+1}\ 	\diff {\pazocal{D}}\nonumber\\
	&&+\frac{1}{\rho}\int_{\pazocal{D}}\bm{F}^t_{n+1} \delta \dot{\bm{F}}_{n+1}
	:\mathbfcal{D}_{n+1}
	:\bm{F}^t_{n+1} \nabla(\bm{w}_{n+1})\ \diff {\pazocal{D}}.
\end{eqnarray}

\subsubsection{Approximation by Finite Elements}

The finite element method (FEM) is applied to the previous equations to obtain the final discretized expressions for the residue vector and the Jacobian matrix. 
We consider the spatial approximations
\begin{equation*}
	\refstepcounter{equation} \latexlabel{appprox1}
	\refstepcounter{equation} \latexlabel{appprox2}
	\refstepcounter{equation} \latexlabel{appprox3}
	\bm{u}_{(\cdot)}\simeq \hat{\bm{N}}\tilde{\bm{u}}_{(\cdot)},\quad
	\delta{\bm{v}}_{(\cdot)}\simeq \hat{\bm{N}}\delta\tilde{\bm{v}}_{(\cdot)},\quad
	\tag{\ref{appprox1}-\ref{appprox3}}
\end{equation*}
\begin{equation*}
	\refstepcounter{equation} \latexlabel{appprox4}
	\refstepcounter{equation} \latexlabel{appprox5}
	{\bm{t}}_{(\cdot)}\simeq \hat{\bm{N}}\tilde{\bm{t}}_{(\cdot)},
	\quad
	{\bm{f}}_{(\cdot)}\simeq \hat{\bm{N}}\tilde{\bm{f}}_{(\cdot)},
	\tag{\ref{appprox4}-\ref{appprox5}}
\end{equation*}
\noindent
where the matrices $\hat{\bm{N}}$ and $\hat{\bm{B}}$ are given in Eqs. \eqref{form_functions} and \eqref{derivative_form_funtions}, respectively.
By using an equivalent product of matrices (see details in Bhatti \cite{bhatti2006advanced}, p.496)) the residue of Eq. \eqref{motion_time_disc} is approximated for each $q$-th element by 
\begin{eqnarray}
 \label{Eq_Residum_motion}
	\bm{R}^q_{n+1} &&\simeq
	\bm{M}^q
	\left(
	a_1(\hat{\bm{u}}_{n+1}^q + \hat{\bm{u}}_n^q
	-a_2\dot{\hat{\bm{u}}}_n^q - a_3\ddot{\hat{\bm{u}}}_n^q) +\hat{\bm{f}}_{n+1}^q
	\right)\nonumber\\
	&&+\frac{1}{\rho}\int_{\pazocal{D}^q} \hat{\bm{B}}^t\bar{\bm{F}}^t_{n+1}\bm{s}_{n+1}^q    
	\ \diff \pazocal{D}^q + BT^q,
\end{eqnarray}
\noindent
where the element mass matrix is
\begin{eqnarray}
	\bm{M}^q = \int_{\pazocal{D}^q} \hat{\bm{N}}^t\hat{\bm{N}} \ \diff \pazocal{D}^q,
\end{eqnarray}
\noindent
$\bm{s}$ is a vector form of the tensor $\bm{S}$, $\bar{\bm{F}}$ is obtained from $\bm{F}$ (see Appendix D) and $BT^q$ are the boundary terms which may depend, for instance, on stresses and displacements.

We obtain the Jacobian matrix $\bm{J}^{q}_{n+1}$ by deriving the residue $\bm{R}^q_{n+1}$ with 
respect to $\bm{w}_{n+1}$:
\begin{eqnarray}
	\label{Eq_jac_motion}
	\bm{J}_{n+1}^q
	&&= \bm{M}^q a_1
	+\frac{1}{\rho}\int_{\pazocal{D}}\hat{\bm{B}}^t\bar{\bm{S}}_{n+1}^q\hat{\bm{B}}\ \diff				\pazocal{D}\nonumber\\
	&&+\frac{1}{\rho}\int_{\pazocal{D}}\hat{\bm{B}}^t\bar{\bm{F}}_{n+1}^t\mathbf{D}_{n+1}^q\bar{\bm{F}}_{n+1}			\hat{\bm{B}}\ \diff\pazocal{D},
\end{eqnarray}
\noindent
where $\bar{\bm{S}}$ is a block-diagonal and symmetric matrix constructed from $\bm{S}$,  $\bar{\bm{F}}$ is obtained from $\bm{F}$ (for details, see Appendix D); and $\mathbf{D}$ is a symmetric matrix representing the double contraction of the fourth order tensor $\mathbfcal{D}$.

Finally, we must solve the final linearized system
\begin{eqnarray}
	\bm{J}_{n+1,i} \Delta \bm{u}_{n+1,i}= \bm{R}_{n+1,i},
\end{eqnarray}
where $i$ is the Newton-Rhapson iteration.
A new approximative solution for $\bm{u}_{i+1}$ is given by
\begin{eqnarray}
	\bm{u}_{n+1,i+1} = \bm{u}_{n+1,i} +  \Delta \bm{u}_{n+1,i}.
\end{eqnarray}
\noindent
The procedure is repeated until $||\bm{u}_{n+1,i+1}-\bm{u}_{n+1,i}||\leq \epsilon$, where $\epsilon$ is a prescribed tolerance and outputs the updated value $\bm{u}_{n+1}$.

\subsubsection{Evaluation of the Tangent Stifness Tensor $\mathbfcal{{D}}$}
\label{tensorA}

The constitutive tensor ${\mathbfcal{D}}$, {that appears in Eq. \eqref{DerivativeS}}, is defined in \cite{bonet2008nonlinear} as
\begin{eqnarray}
	\mathbfcal{{D}}
	:=\frac{\partial \bm{S}}{\partial \bm{E}}
	=\frac{1}{2}\frac{\partial \bm{S}}{\partial\bm{C}},
\end{eqnarray}
\noindent
resulting in a fourth order symmetric positive definite tensor. 
In order to obtain a suitable matrix multiplication in Eq. \eqref{Eq_jac_motion}, tensor $\mathbfcal{{D}}$ is rewritten as a symmetric matrix $\mathbf{D}$.

The derivative of $\bm{S}$ in relation to $\bm{E}$ (or $\bm{C}$) must be calculated by using Eq. \eqref{final_second_piola}. 
The difficulty in deriving the final expression for $\mathbf{D}$ is evident, since Eq. \eqref{final_second_piola} has many nonlinear dependencies on $\bm{E}$.
In order to overcome this issue, we perform a numerical complex derivative for each component $D_{pq}$ using the relation \cite{haveroth2015application}
\begin{eqnarray}
	\label{complexx}
	({{D}}_{n+1})_{pq}
	=\frac{\partial {S}_p}{\partial {C}_q}=\frac{\text{Im}{\left(\bm{S}_{n+1}\left((\bm{C}_{n+1})+i\hat{\delta}			\right)_p\right)_q}}{\hat{\delta}},
\end{eqnarray}
\noindent
where $p,q=1,2,3$, $\hat{\delta}$ is a small perturbation ($\hat{\delta} \in[10^{-100},
10^{-300}]$) and $i$ the imaginary unit.
This method presents advantages due to the single term in the numerator of Eq. \eqref{complexx}. 
It avoids the instability related to cancellation error inherent to all real valued finite difference approximations. 
Furthermore, the complex finite difference method is more accurate when compared with the real valued finite difference method.

\subsection{Equation of Damage}
\label{Sec_Eq_damage}

Consider the damage evolution given in Eq. \eqref{final_damage}. 
Firstly, we apply the backward Euler method for time discretization obtaining
\begin{eqnarray}
	\varphi_{n+1}
	&=&	\varphi_n+\frac{\Delta t}{\tilde{\lambda}_{n+1}\theta_{n+1}}
	\bigg[\divv \left(g_c\gamma\bm{C}_{n+1}^{-1}\Nabla \varphi_{n+1} \right)\nonumber\\
	&&-\frac{g_c}{\gamma}H'_{n+1}
	-{G'_{n+1}}\left({\psi}_h 
	+\tilde{\psi}_m\right)_{n+1}\bigg],
\end{eqnarray}
\noindent
where 
\[\left({\psi}_h +\tilde{\psi}_m\right)_{n+1}:={\psi}_h (\bm{E}_{n+1}) +\tilde{\psi}_m(\bm{E}_{n+1}),\]
\[G'_{n+1}:=G'(\varphi_{n+1}),\] and
\begin{eqnarray}
\label{Eq_deriv_H}
	H'_{n+1}:=H'(\varphi_{n+1})=\varphi_{n+1},
\end{eqnarray}
\noindent
according to the definition of $H$ in Sec. \ref{sec_free_energy}.

Before applying the spatial discretization by finite elements, consider the following modification for the divergent term in the previous equation:
\begin{eqnarray}
	\varphi_{n+1}
	&=&\varphi_n
	+\Delta t g_c\gamma
	\left[\divv \left(\frac{1}{\tilde{\lambda}_{n+1}\theta_{n+1}}\bm{C}_{n+1}^{-1}\Nabla (\varphi_{n+1})\right)\right. \nonumber\\
	&&\left.\left.-\Nabla\left(\frac{1}{\tilde{\lambda}_{n+1}\theta_{n+1}}   \right)\cdot\bm{C}_{n+1}^{-1}\Nabla (\varphi_{n+1}) \right.\right]\\
	&&-\frac{\Delta t}{\tilde{\lambda}_{n+1}\theta_{n+1}}\left[\frac{g_c}{\gamma}H'_{n+1}
	+{G'_{n+1}}\left({\psi}_h + \tilde{\psi}_m\right)\right].\nonumber
\end{eqnarray}
Then, the weak form for the damage phase-field evolution is obtained by multiplying the previous equation by a suitable scalar test function ${\omega}$ and integrating over the domain $\pazocal{D}$. 
Therefore,
\begin{eqnarray}
	\label{weakdamage}
	&&\int_{\pazocal{D}}\varphi_{n+1} {\omega}\ \diff \pazocal{D}
	=\int_{\pazocal{D}}\varphi_n {\omega} { \diff \pazocal{D}}\nonumber\\
	&&+\Delta t g_c\gamma \int_{\pazocal{D}}\divv \left(\frac{1}{\tilde{\lambda}_{n+1}\theta_{n+1}}				\bm{C}_{n+1}^{-1} \Nabla \varphi_{n+1} \right) {\omega}\ \diff \pazocal{D}\nonumber\\
	&&-\Delta t g_c\gamma \int_{\pazocal{D}}\Nabla\left(\frac{1}{\tilde{\lambda}_{n+1}\theta_{n+1}} \right)\cdot\bm{C}_{n+1}^{-1}					\Nabla \varphi_{n+1}   {\omega}\ \diff \pazocal{D}\nonumber\\
	&&-\frac{\Delta t g_c}{\gamma }\int_{\pazocal{D}}\frac{1}{\tilde{\lambda}_{n+1}\theta_{n+1}}H'_{n+1} {\omega}\ \diff 				\pazocal{D}\nonumber\\
	&& -\Delta t \int_{\pazocal{D}} \frac{G'_{n+1}}{\tilde{\lambda}_{n+1}\theta_{n+1}}\left(\psi_h + \tilde{\psi}_m\right) {\omega} \ \diff \pazocal{D}.
\end{eqnarray}

Now, consider the gradient properties and Eq. \eqref{um_sobre_lambda} to write
\begin{eqnarray}
\label{Eq_property_grad}
	\Nabla\left(\frac{1}{\tilde{\lambda}\theta}\right)
	&=&\frac{1}{\theta}\Nabla\left(\frac{1}{\tilde{\lambda}}\right)
	+\frac{1}{\tilde{\lambda}}\Nabla\left(\frac{1}{\theta}\right)\nonumber\\
	&=&\frac{\zeta{c}_\lambda}{\theta(1+\tilde{\delta}-\varphi)^{\zeta+1}}\Nabla \varphi
	-\frac{1}{\tilde{\lambda}\theta^2}\Nabla \theta.
\end{eqnarray}
\noindent
By using Eqs. \eqref{Eq_deriv_H} and \eqref{Eq_property_grad} into Eq. \eqref{weakdamage}, applying the Green's theorem and assuming $\theta$ and $\lambda$ delayed (in order to avoid numerical instability), we obtain
\begin{eqnarray}
	\label{Eq_residue_before}
	&&\int_{\pazocal{D}}\varphi_{n+1}{\omega}\ \diff \pazocal{D}
	=\int_{\pazocal{D}}\varphi_{n+1} {\omega} \diff\pazocal{D}\nonumber\\
	&&-\Delta t g_c\gamma\int_{\pazocal{D}}\frac{1}{\tilde{\lambda}_{n}\theta_{n}} \bm{C}_{n+1}^{-1} \Nabla \varphi_{n+1} \cdot \Nabla {\omega} \ \diff \pazocal{D}\nonumber\\
	&&-\Delta t g_c{\gamma}\zeta{c}_\lambda\int_{\pazocal{D}}
	\frac{\Nabla \varphi_n \cdot\left( \bm{C}_{n+1}^{-1}\Nabla \varphi_{n+1} \right)}{\theta_{n}(1+\tilde{\delta}-\varphi_N)^{\zeta+1}}  {\omega}\ \diff 	\pazocal{D}\nonumber\\
	&&+\Delta t g_c\gamma\int_{\pazocal{D}}\frac{1}{\tilde{\lambda}_n\theta_n^2}\Nabla \theta_n \cdot \left(\bm{C}^{-1}_{n+1} \Nabla \varphi_{n+1}\right) {\omega}\ \diff \pazocal{D}\nonumber\\
	&&-\frac{\Delta t g_c}{\gamma}\int_{\pazocal{D}}\frac{1}{\tilde{\lambda}_n\theta_n}\varphi_{n+1}{\omega}\ 				\diff \pazocal{D}\nonumber\\
	&&-\Delta t\int_{\pazocal{D}}\frac{1}{ \tilde{\lambda}_n\theta_n}G'_{n+1}			\left(\psi_h + \tilde{\psi}_m\right)_{n+1}{\omega}\ \diff \pazocal{D}.
\end{eqnarray}

{Adopting the spatial approximations}
\begin{equation*}
	\refstepcounter{equation} \latexlabel{approx11}
	\refstepcounter{equation} \latexlabel{approx22}
	\varphi_{(\cdot)}^q\simeq\bm{N}\tilde{\varphi}_{(\cdot)}^q,
	\quad \Nabla(\varphi_{(\cdot)}^q)\simeq \bm{B} \tilde{\varphi_{(\cdot)}}^q,
	\tag{\ref{approx11}-\ref{approx22}}
\end{equation*}
\begin{equation*}
	\refstepcounter{equation} \latexlabel{approx55}
	\refstepcounter{equation} \latexlabel{approx66}
	\theta_{(\cdot)}^q\simeq\bm{N}\tilde{\theta}_{(\cdot)}^q,
	\quad\Nabla (\theta^q) \simeq\bm{B}\tilde{\theta}_{(\cdot)}^q,
	\tag{\ref{approx55}-\ref{approx66}}
\end{equation*}
\begin{equation*}
	\refstepcounter{equation} \latexlabel{approx33}
	\refstepcounter{equation} \latexlabel{approx44}
	w^q\simeq\bm{N}\tilde{w}^q,
	\quad\Nabla (w^q )\simeq\bm{B}\tilde{w}^q,
	\tag{\ref{approx11}-\ref{approx44}}
\end{equation*}
\noindent
where the matrices $\bm{N}$ and $\bm{B}$ are given in Eqs. \eqref{noriginal} and \eqref{B1}, respectivelly, and making $\nabla(\varphi)$ delayed in the third and fourth terms of the right hand side of Eq. \eqref{Eq_residue_before} (to avoid non-symmetric Jacobian matrix), we obtain the residue for each  element $q$ at time step $n+1$ for the damage equation as
\begin{eqnarray}
\label{Eq_residue_damage}
	&&\bm{R}_{n+1}^{q, \textmd{damage}}
	=\int_{\pazocal{D}_q} \bm{N}^t \bm{N}\left[
	\left(1+\frac{\Delta t g_c}{\gamma {\tilde{\lambda}_n^q}\bm{N}\tilde{\theta}_{n}^q}\right)\tilde{\varphi}_{n+1}^q
	-\tilde{\varphi}_n^q
	\right] \ \diff \pazocal{D}_q\nonumber\\
	&&+{\Delta t g_c\tilde{\gamma}}\int_{\pazocal{D}_q} \frac{\bm{B}^t\left(\bm{C}_{n+1}^{q}\right)^{-1}			\bm{B}\tilde{\varphi}_{n+1}^q}{\tilde{\lambda}_n^q\bm{N}\tilde{\theta}_n^q}\ \diff 	\pazocal{D}_q\nonumber\\
	&&+\Delta t g_c \gamma \zeta{c_\lambda}\int_{\pazocal{D}_q}\frac{\bm{N}^t\left(\tilde{\varphi}^{q}_{n}\right)^t\bm{B}^t\left(\bm{C}_{n+1}^{q}\right)^{-t}\bm{B}\tilde{\varphi}_{n}^q}{\bm{N}\tilde{\theta}_n^q(1+\tilde{\delta}-\bm{N}			\tilde{\varphi}_n^q)^{\zeta+1}}\ \diff \pazocal{D}_q\nonumber\\
	&&-\Delta t g_c\gamma\int_{\pazocal{D}_q}\frac{\bm{N}^t\left(\tilde{\varphi}^{q}_{n}\right)^t\bm{B}^t\left(\bm{C}_{n+1}^{q}\right)^{-t} \bm{B}\tilde{\theta}_n^q}{\tilde{\lambda}_n^q(\bm{N}\tilde{\theta}_{n}^{q})^2} \diff \pazocal{D}_q\nonumber\\
	&&+{\Delta t} \int_{\pazocal{D}_q}\frac{\bm{N}^t \left(G_{n+1}^q\right)' \left(\psi_h+\tilde{\psi}_m\right)_{n+1}^q}{\tilde{\lambda}_n^q\bm{N}\tilde{\theta}_{n}^q} \ \diff \pazocal{D}_q,
\end{eqnarray}
\noindent
where
\begin{eqnarray}
	\frac{1}{\tilde{\lambda}_n^q}
	=\frac{c_\lambda}{(1+\tilde{\delta}-\varphi_n^q)^{\zeta}} 
	\simeq \frac{c_\lambda}{(1+\tilde{\delta}-\bm{N}\tilde{\varphi}_n^q)^{\zeta}},
\end{eqnarray}
\noindent
due to Eq. \eqref{um_sobre_lambda}, and $G_{n+1}^q:=G(\varphi_{n+1}^q)$.
The respective Jacobian matrix $\bm{J}^{q, \textmd{damage}}_{n+1}$ is obtained by deriving Eq. \eqref{Eq_residue_damage} to  $\varphi_{n+1}^k$:
\begin{eqnarray}
	&&\bm{J}^{q,\textmd{damage}}_{n+1}
	=\int_{\pazocal{D}_q} \bm{N}^t \bm{N}\left(
	1+\frac{\Delta t g_c}{\gamma {\tilde{\lambda}_n}\bm{N}\tilde{\theta}_n^q}
	\right) \ \diff \pazocal{D}_q\nonumber\\
	&&+{\Delta t g_c\gamma}\int_{\pazocal{D}_q}\frac{\bm{B}^t\bm{C}_{n+1}^{-1}\bm{B}}{\tilde{\lambda}_n\bm{N}\tilde{\theta}_{n}^q}\ 			\diff \pazocal{D}\nonumber\\
	&&+{\Delta t} \int_{\pazocal{D}_q}\frac{\bm{N}^t \left(G_{n+1}^q\right)'' \left(\psi_h+\tilde{\psi}_m\right)_{n+1}^q}{\tilde{\lambda}_n\bm{N}\tilde{\theta}_{n}^q} \ \diff \pazocal{D}_q.				\nonumber\\
\end{eqnarray}

Finally, for each time step $n+1$, we solve iteractively the global linearized system
\begin{eqnarray}
\bm{J}^{\textmd{damage}}_{n+1,i} \Delta {\varphi}_{n+1,i}= -\bm{R}^{\textmd{damage}}_{n+1,i},
\end{eqnarray}
\noindent
where $i$ is the Newton-Raphson iteration.
Matrix $\bm{J}^{\textmd{damage}}_{n+1,i}$ is the global Jacobian matrix obtained by assembling each $q$-th local Jacobian $\bm{J}^{q, \textmd{damage}}_{n+1,i}$.
Similarly, the global residue $\bm{R}^{\textmd{damage}}_{n+1,i}$ is obtained by assembling the local residue vector $\bm{R}^{q,\textmd{damage}}_{n+1,i}$.

A new approximate solution for ${\varphi}_{n+1}$ is given by
\begin{eqnarray}
{\varphi}_{n+1, i+1} = {\varphi}_{n+1,i} +  \Delta {\varphi}_{n+1,i}.
\end{eqnarray}
\noindent
The procedure is repeated until $||\varphi_{n+1,i+1}-\varphi_{n+1,i}||\leq \epsilon$, where $\epsilon$ is a prescribed tolerance and outputs the updated values $\varphi_{n+1}$.
In the first iteration of each time step, we adopt
$\varphi_{n+1,0}:=\varphi_{n}$, where $\varphi_n$ is the damage value of the previous step. 

\vspace{0.2cm}
As discussed in Remark \ref{Remark_irrevers_damage}, damage irreversibility
will be imposed by using a kind of predict-corrector procedure as follows.
Assuming known the state values at time step $n$,
we use the evolution equation without the additional term in the pseudo-potential to predict the damage values
at time step $n+1$,
obtaining a predict value  $\varphi^*_{n+1}$ for the damage variable. 
Next, for each node of the mesh, we compare $\varphi^*_{n+1}$ to $\varphi_{n}$:
if at that node $\varphi^*_{n+1} \geq \varphi_{n}$ then we take $\varphi_{n+1} := \varphi^*_{n+1}$;
otherwise
$\varphi_{n+1} := \varphi_{n}$.

Furthermore, for simplicity we prescribe $\varphi_0=0$ (undamaged material) to start the analysis;
we could take any given damage state to initiate the evolution.

\subsection{Numerical fractional derivative}
\label{fractional_numerical}

\noindent
Oldhan and Spanier \cite{oldham1974} {used} the numerical algorithm G1 to calculate fractional derivatives. 
The expression for this approximation is given by
\begin{eqnarray}
	\label{riet}
	{\left. _0\mathrm{D}^\alpha_{t} f(t) \right|}_{G1} 
	=
	{(\Delta t)}^{-\alpha}\sum_{m=0}^{N-1}A_{m+1}f_m ,
\end{eqnarray}
where the coefficients $A_{m+1}$ are given by
\begin{eqnarray}
	\label{ret}
	A_{m+1} 
	= \frac{\Gamma(m-\alpha)}{\Gamma(-\alpha)\Gamma(m+1)}
	= \dfrac{m - 1 - \alpha}{m} A_m.
\end{eqnarray}
Herein, $\Delta t=t/N$ is the time increment, $N \in [1, \infty)$ is the number of time steps and $f_m = f(t-m\Delta t)$. 
If $f(0)=0$, then the algorithm G1 can be used as an approximation for the Caputo fractional derivative. 
For strain free materials in the initial time ($\bm{E}(0)=\bm{0}$), as considered in this work, we calculate the fractional derivative of Eq. \eqref{final_second_piola} by using the algorithm G1. 

\section{Results and Discussion}
\label{ResultsDiscussion}

This section presents some results and comments for the model proposed in this work. 
Initially, the one-dimensional version of the model is used to simulate tensile tests in a viscoelastic rod. Next, two-dimensional examples are considered, including comparison with experimental data.

\subsection{Viscoelastic Rod}
\label{viscoelastic_rod}

\begin{figure}[!h]
	\centering
	\begin{subfigure}{.49\textwidth}
		\centering
		\includegraphics[scale=0.4]{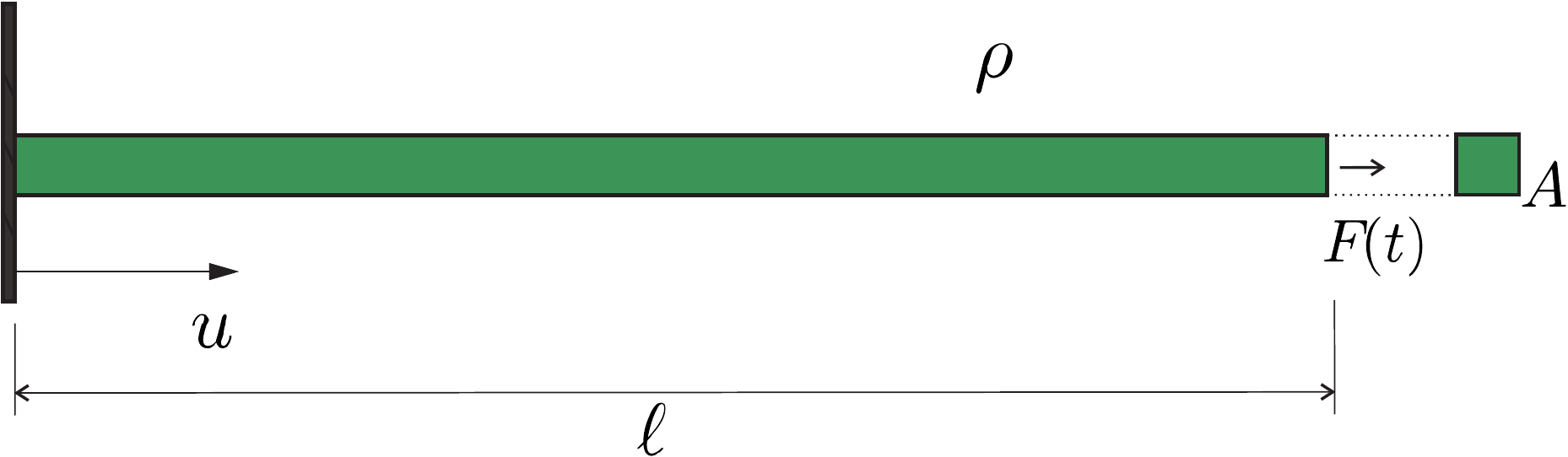}
		\caption{\small{Conditions for the polyoxymethylene rod.}}
		\label{rod}
	\end{subfigure}
	\begin{subfigure}{0.49\textwidth}
		\centering
		\hspace{-0.7cm}
		\includegraphics[scale=0.4]{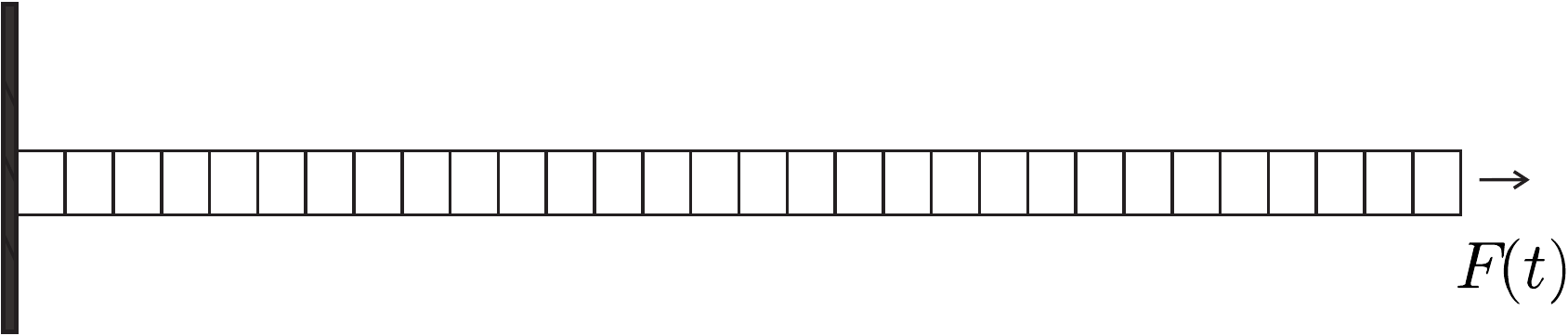}
		\caption{\small{Mesh used in the two-dimensional simulation.}}
		\label{rod2}
	\end{subfigure}
	\caption{\small{Polyoxymethylene rod.}}
\end{figure}

Consider a polyoxymethylene viscoelastic rod \cite{schimidt2006} with density $\rho=1420\ kg/m^3$, length $\ell=2\ m$ and a squared cross section of area $A=176.71459\ mm^2$ fixed at $x=0$ and subject to an external force given by $F(t)$ at $x=\ell$ (see Fig. \ref{rod}). 
We promoted one-dimensional dynamic tensile tests to evaluate the contribution of the terms in the stress, given by Eq. \eqref{final_second_piola}, and to study the behavior of the displacement concerning the fractional viscoelastic parameters.
The rod is discretized into 30 equally spaced elements with 2 integration points for the application of the finite element method. 
The time discretization is considered by using the Newmark method with $\beta=0.25$ and time increment $\Delta t=1\times 10^{-4}s$. 
The tolerance for the Newton-Raphson method is $10^{-8}$. 
We also consider no viscous dissipative damping, i.e. $\tilde {b}=0$, neither damage effects ($\tilde{c}=0$). 

\subsubsection{Evaluation of the Stress Terms}
\label{Evaluation_of_stress_terms}

{The constitutive stress/strain equation for our model is given by Eq. \eqref{final_second_piola}.} 
The inclusion of the last two terms {in this equation} represents a great computation cost in simulations, especially due to the singular integral that appears in the last one. 
Thus, we perform tests to evaluate the relevance of these terms and establish the situations in which they do not influence significantly in the results for stress.   
To this end, we refer to Eq. \eqref{final_second_piola} as the complete stress $\bm{S} =\bm{S}_c$. 
On the other hand, the expression obtained disregarding the last two terms of Eq. \eqref{final_second_piola} is called partial stress $\bm{S}_p$.

Tensor $\mathbfcal{{A}}$ is reduced to a scalar for the one-dimensional case, whose corresponding equation is chosen to be
 \begin{eqnarray}\label{A1}
	\mathbfcal{{A}}
	:=\mathbfcal{{A}}_1 
	= \overline{\lambda}\bm{C}^{-1}\otimes \bm{C}^{-1}+2(\overline{\mu}-					\overline{\lambda}\ln(J))\bm{I},
\end{eqnarray}
where
\begin{equation*}
	\refstepcounter{equation} \latexlabel{mlabda}
	\refstepcounter{equation} \latexlabel{mmu}
	\overline{\lambda}=\frac{p\nu}{(1+\nu)(1-2\nu)},
	\quad
	\overline{\mu}=\frac{p}{2(1+\nu)},
	\tag{\ref{mlabda}-\ref{mmu}}
\end{equation*}
and $p$ is a viscoelastic material parameter similar to the Young's modulus for the elastic case. 
Parameters $\overline{\lambda}$ and $\overline{\mu}$ are modifications of the usual Lam\'e constants and $\nu$ is the Poisson ratio. 
In fact, tensor $\mathbfcal{{A}}_1$ is a generalization of the traditional elastic tensor written in terms of the Young's modulus. 
We assume $p=21.46\times 10^6\ N/m^2s^\alpha$ and Young's modulus $E_Y=1430.1\times 10^6\ Pa$.
The numerical fractional derivative is calculated by using the algorithm G1 given in Eq.\eqref{riet}.

The resulting strain at the right end of the specimen for the second integration point at the end of the simulation is considered. 
The difference between $\bm{S}_c$ and $\bm{S}_p$ is calculated by using the mean square difference (MSD):
\begin{eqnarray}
\text{MSD} = \sqrt{\frac{1}{N}\sum_{i=1}^N\frac{(S_{c,i}-S_{p,i})^2}{S_{c,i}^2}},
\end{eqnarray} 
where $S_{c,i}$ and $S_{p,i}$ are the components of the complete and the partial stress, respectively, for each time step $i=1, \cdots, N$.

\begin{figure}[!h]
\centering
	\begin{subfigure}{.45\textwidth}
		\includegraphics[scale=0.45]{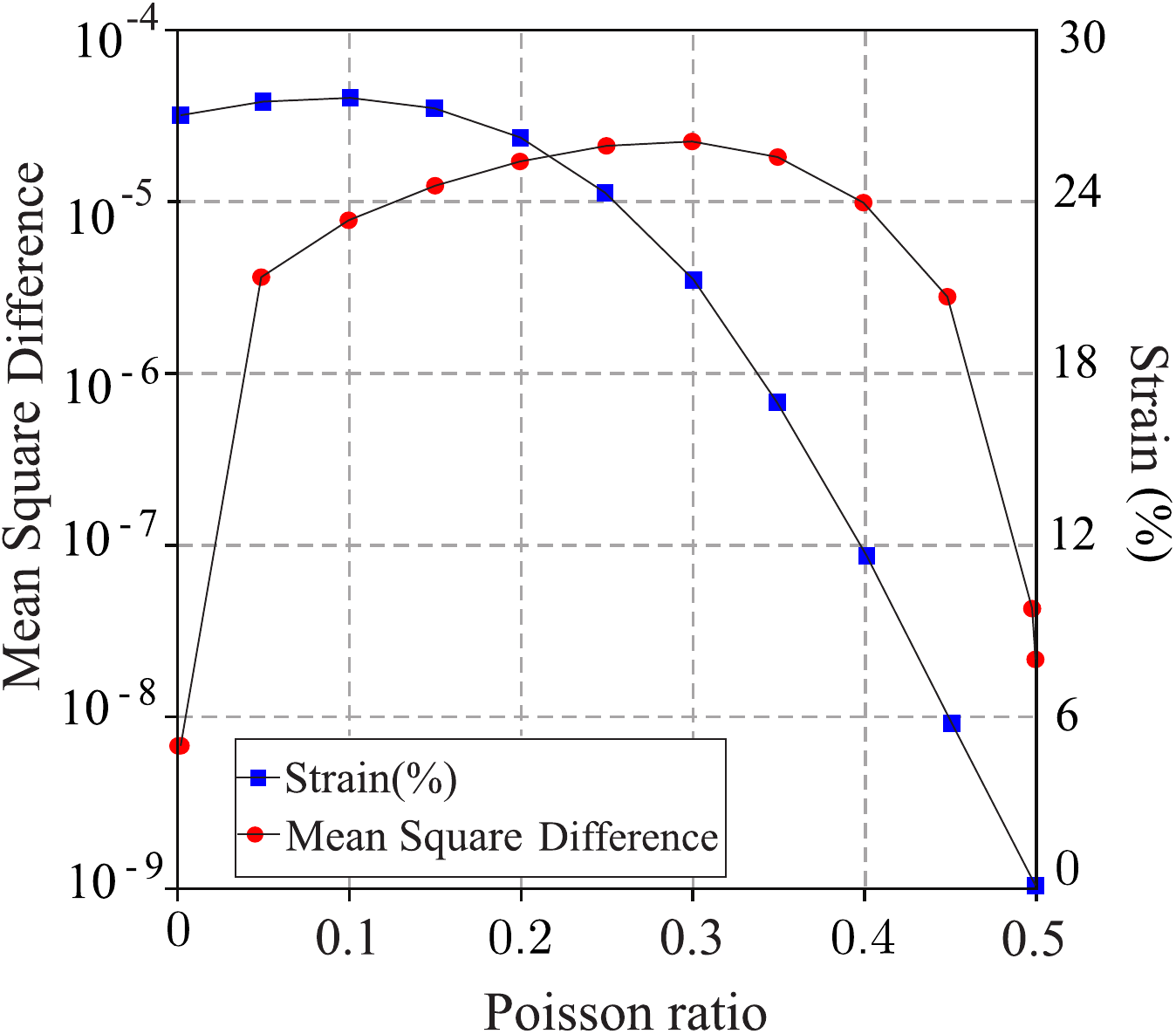} 
		\caption{\small{Left axis: MSD between total stress and partial stress for the variation of $\nu$. Right axis: percentage of strain for the variation of $\nu$}}\vspace{0.25cm}
		\label{errorpoisson}
	\end{subfigure}
	\begin{subfigure}{.45\textwidth}
		\includegraphics[scale=0.45]{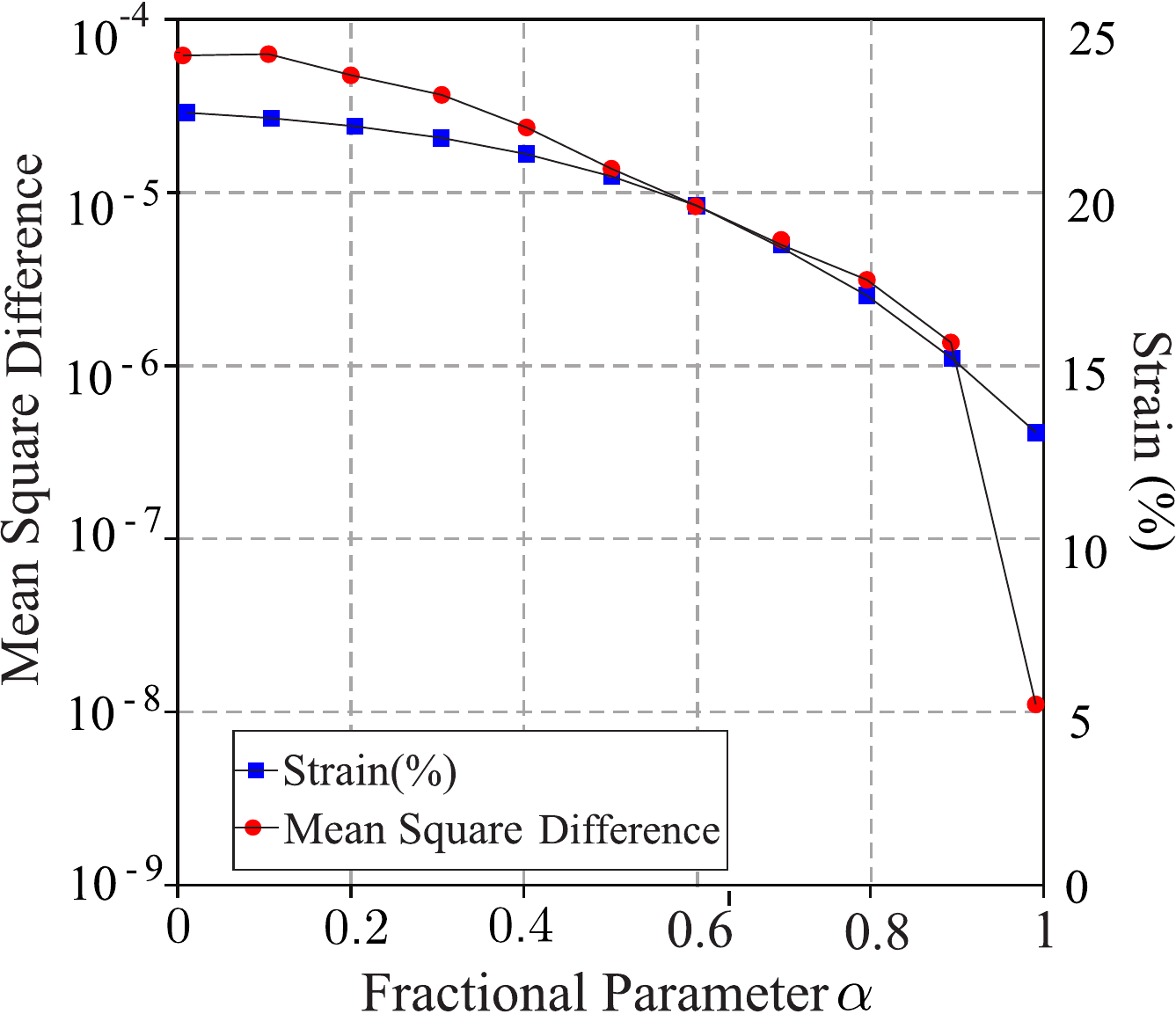}
		\caption{\small{Left axis: MSD between total stress and partial stress for $\alpha$ variation. Right axis: percentage of strain for $\alpha$ variation.}}
		\label{erroralpha}
	\end{subfigure}
	\caption{\small{Evaluation of the mean square error and the strain level for variation of Poisson's ratio $\mu$ and the fractional parameter $\alpha$.}}
\end{figure}

Figure \ref{errorpoisson} presents the MSD between $\bm{S}_c$ and $\bm{S}_p$ and the strain percentage for simulations with $\nu$ varying between  0.0001 and 0.4999. 
The total time of analysis is $0.05s$, the applied force is $10\times 10^5\ N$ and the fractional derivative order is $\alpha = 0.5$. 
It can be observed that when $\nu$ approximates to 0.3 the difference increases, i.e., for values near to 0.3 the influence of the last two terms in $S_c$ is more significant than for the remaining values. 
On the other hand, the specimen presents a strain bigger than $21\%$ for values of $\nu$ smaller than 0.3, and decreases to $0.011\%$ when $\nu$ approximates to $0.5$. 
This behavior was expected once $\nu=0.5$ corresponds to an incompressible material. 
Furthermore, the strain $\bm{E}$ in the last two terms of Eq. \eqref{final_second_piola} has less influence for small strain. 
Although higher values of $\nu$ indicate an increase in $\partial_{\bm{E}} \mathbfcal{A}$, it also leads to the reduction of the strain.

Figure \ref{erroralpha} presents the MSD and the resulting strain for $\alpha$ varying between 0.001 and 0.999. 
The analysis time is $0.05\ s$, the applied force is $10\times 10^5\ N$ and the Poisson's ratio is $0.3$.  
As $\alpha$ approaches to 1, both the MSD  and the percentage of strain decreases. 

Table \ref{Tab_FOrce_varying} presents the MSD and the percentage of strain for different values of loads and final times. 
In this case, the fractional derivative order is $\alpha=0.5$ and the results are presented for $\nu=0.3$ and $\nu=0.45$.

\begin{table}[!h]
	\small
	\centering
	\caption{\small{MSD and percentage of strain for the the applied forces.}}
	\begin{tabular}{|c|c|c|}
	\hline 
	\hline
	\multicolumn{3}{|c|}{Final Time (s) = 0.05 and $\nu= 0.3$} \\
	\hline 
	Force ($kN$) & Strain (\%) & MSD ($\times 10^{-6}$)\\
	\hline 
	$200$ & $2.92$ & $2.6615$\\
	$400$ & $6.30$ & $5.3166$ \\
	$800$ & $15.07$ & $10.664$ \\
	$1000$ & $21.02$  &  $13.498$\\  
	\hline
	\hline
	\multicolumn{3}{|c|}{Final Time (s) = 0.05 and $\nu=  0.45$}\\
	\hline
	Force ($kN$) & Strain (\%) & MSD ($\times 10^{-6}$) \\
	\hline
	$400$ & $2.05$ & $1.5109$ \\
	$1000$ & $5.70$  &  $3.8064$ \\
	$2000$ & $14.41$ & $7.9290$ \\
	$2500$ & $21.70$  &  $23.942$ \\
	\hline
	\hline
	\multicolumn{3}{|c|}{$\nu=  0.45$}  \\
	\hline
	\multicolumn{3}{|c|}{Force $1000 kN$}  \\ 
	\hline
	Final Time(s) & Strain (\%) & MSD ($\times 10^{-6}$)  \\
	\hline
	0.075 & $9.62$ & $ 3.1115$   \\ 
	0.05 & $5.70$  &  $2.7057$ \\ 
	0.1 & $14.65$ & $3.8064$ \\
	\hline
	\hline
	\end{tabular}
	\label{Tab_FOrce_varying}
\end{table}

In all the cases analyzed, the MSD is not bigger than $1.0 \times 10^{-4}$, even when the percentage of strain is large (bigger than $5\%$). 
It implies that the influence of the last two terms in Eq. \eqref{final_second_piola} can be neglect for these cases without significant influence on the stress evaluation.
Then, the remaining simulations presented in this paper are performed by disregarding these terms. 
In other words, we consider $\bm{S}=\bm{S}_p$ for the analyses that follows. 

\subsubsection{Displacement of the Rod}
\label{Displacement_Rod}

Now, we use the model proposed in Sec. \ref{development_of_the_model} to describe the dynamic response of the viscoelastic rod of Fig. \ref{rod} when it is subject to a force
\begin{eqnarray}
	F(t) 
	= \begin{cases} 
	0, & \mbox{if }  t=0 \\ 
	100\ N, & \mbox{if }  t>0 
	\end{cases},
\end{eqnarray}
during $0.1\ s$.
It is done in order to check the displacement behavior and the viscoelastic effect induced by the fractional derivative. 
Considering the magnitude of the applied force, the analysis time and the bar dimensions, we consider small strain regime, and Eq. \eqref{final_second_piola} is simplified, replacing the Neo-Hookean by a linear spring. 
We also remember that damage is not considered. 
Furthermore, we assume that the material does not have nonlinearities due to the fourth order tensor $\mathbfcal{A}$, that appears in the stress/strain relation \eqref{final_second_piola}.
Then, it can be simplified to a scalar parameter $p$, with the same purpose of $\mathbfcal{A}$ in weighting the fractional derivative.
We consider the Poisson's ratio $\nu = 0.39$ and use the numerical fractional derivative algorithm G1 of Eq. \eqref{riet}.

\begin{figure*}[!h]
	\centering
	\includegraphics[scale=0.45]{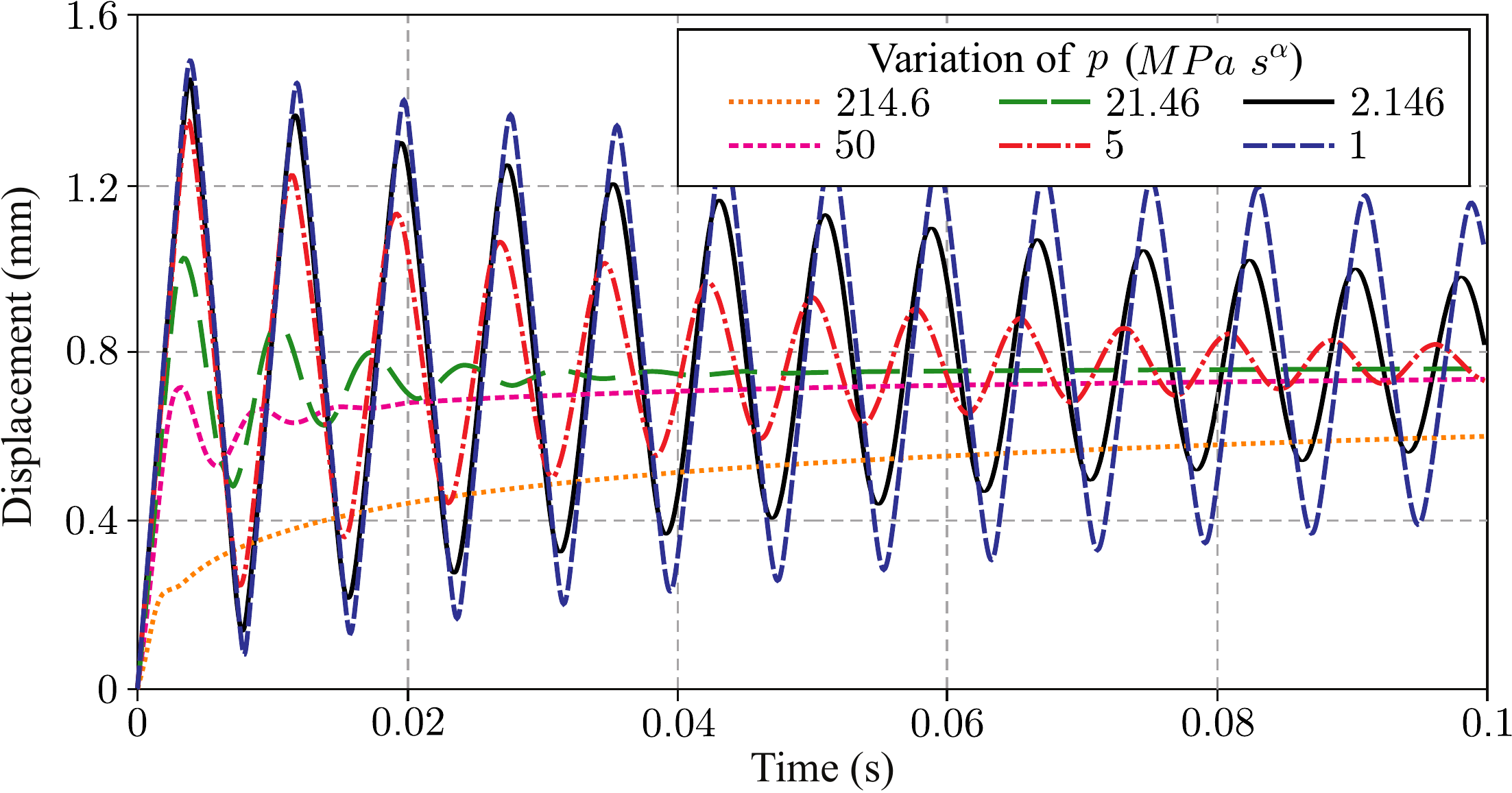} 
	\caption{\small{Displacement at the free end of the viscoelastic rod for $\alpha=0.5$ and different values of $p$ $\left[(N/m^2)s^\alpha\right]$.}}
	\label{comp1}
\end{figure*}

The application of the force $F(t)$ results in an oscillatory displacement at the free end of the rod. 
This behavior can be seen in Fig. \ref{comp1} for some values of $p$ and $\alpha=0.5$. 
When $p$ increases, the damping effect also grows. 
This was expected because $p$ weights the viscoelastic behavior.
These results are important because they give qualitative information on how to control the viscoelasticity effects by changing parameter $p$. 

\subsubsection{Two-Dimensional Case}
\label{two_dimentsional_calse}

The results shown in Fig. \ref{comp1} were used as reference to extend our model for the two-dimensional case including the possibility of large strains (i.e., adopting the hyperelastic Neo-Hookean spring as shown in Eq. \eqref{final_second_piola}).

Two alternatives for the tensor $\mathbfcal{{A}}$ were tested here. 
The first one is $\mathbfcal{{A}}=\mathbfcal{{A}}_1$ as in Eq. \eqref{A1}, and the second tensor proposed is $\mathbfcal{{A}}:=\mathbfcal{{A}}_2$ where $\mathbfcal{{A}}_2$ is a fourth order tensor with  ${A}(1,1,1,1)=p$ and ${A}(i,j,k,l)=0$ for {other components}.

The viscoelastic rod is discretized using a mesh of 30 quadratic elements as shown in Fig. \ref{rod2}, and two integration points for each element.
This corresponds to an equivalent two-dimensional version of the problem considered in the previous section.

Plane stress state is used and remaining information is the same as that used in the previous section. 

Although this problem {is in the small strain regime}, hyper elasticity was included aiming to test the model for more general problems. 

Table \ref{table1} shows the mean square difference (MSD) for the displacement at the free end of the rod between the one and the two-dimensional models, for some values of $\alpha$. The MSD is calculated by
\begin{eqnarray}
	\text{MSD}=\sqrt{\frac{1}{N}\sum_{i=1}^N\frac{\left(d_i-\bar{d}_i \right)^2}{\bar{d}_i^2}},
\end{eqnarray}
{where $d_i$ and $\bar{d}_i$ are the displacements for the one and the two-dimensional cases, respectively, $i=1,\cdots,N$ and $N$ is the number of time steps.}

\begin{table}[!h]
	\small
	\centering
	\caption{\small{Mean square difference (MSD) between the one and the two-dimensional models for the displacement of the free end of the rod.}}
	\begin{tabular}{|l|l|l|l|}
	\hline
	\multicolumn{4}{|c|}{MSD}\\ 
	\hline 
	\hline 
	$\alpha$  & $p$ & \multicolumn{1}{c|}{$\mathbfcal{A}_1$} & \multicolumn{1}{l|}{$\mathbfcal{A}_2$} \\ 	\hline 			
	\hline
	0.00794  & $214.6\times 10^6$& $0.0097$  & $0.0097$        \\ 
	\hline
	0.2 & $214.6\times 10^6$ &$  9.7113\times 10^{-5}$ &  $9.7232\times 10^{-5}$     \\ 
	\hline
	0.5 &  $214.6\times 10^4$& 	$2.2396\times 10^{-4}$	&  	$2.2416\times 10^{-4}$ \\ 
	\hline
	0.5 & $214.6\times 10^6$& $1.0576\times 10^{-5}$ &   $1.0559\times 10^{-5}$ 	\\ 
	\hline
	0.7 & $214.6\times 10^4$& $8.0830\times 10^{-5}$& 	$8.0877\times 10^{-5}$ 	  \\ 
	\hline
	0.9 & $214.6\times 10^4$& $2.3415\times 10^{-5}$ &   $7.2279\times 10^{-14}$		\\ 
	\hline
	\end{tabular}
	\label{table1}
\end{table}

The magnitude of the error presented in Tab. \ref{table1} shows that both proposals for $\mathbfcal{{A}}$ lead to a reliable two-dimensional extension. 
Since tensor $\mathbfcal{{A}}_1$ can be considered a natural extension for the usual elastic tensor, it will be used to simulate the viscoelastic materials in the examples that follow.

\subsection{I-shaped Viscoelastic Specimen}

\begin{figure*}[!h]
\centering
	\begin{subfigure}{.7\textwidth}
		\includegraphics[scale=0.39]{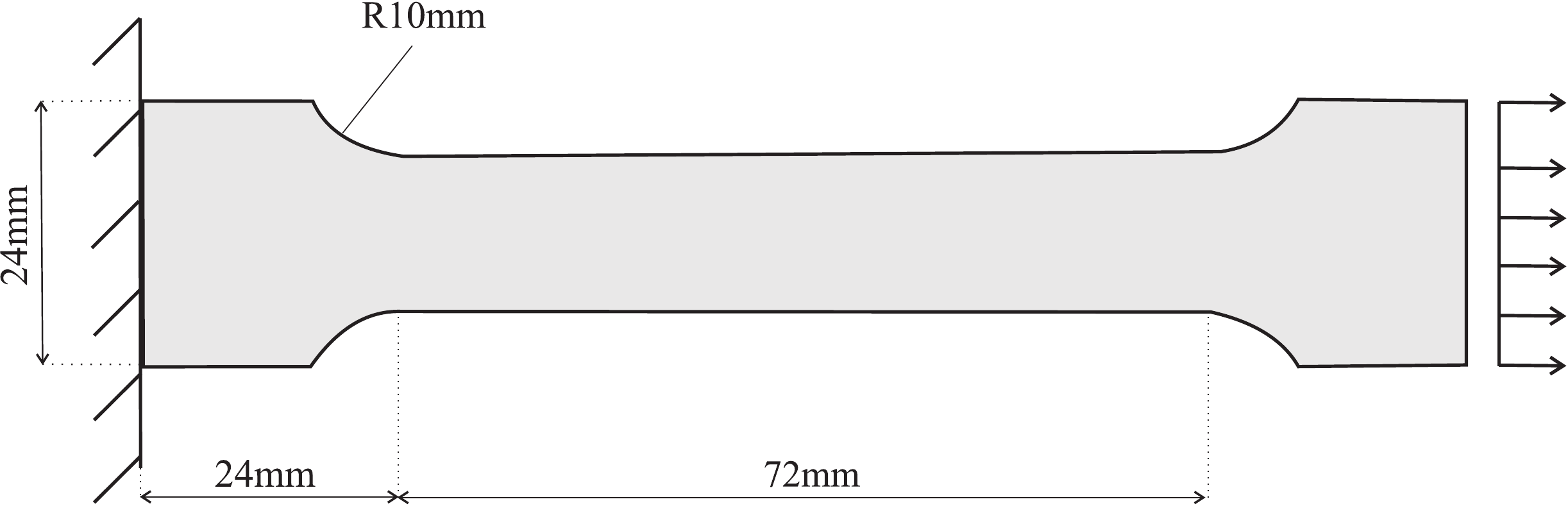} 
		\caption{\small{Dimensions of the specimen.}}
		\label{ishaped}
	\end{subfigure}\\
	\begin{subfigure}{.7\textwidth}
	\hspace{0.7cm}
		\includegraphics[scale=0.385]{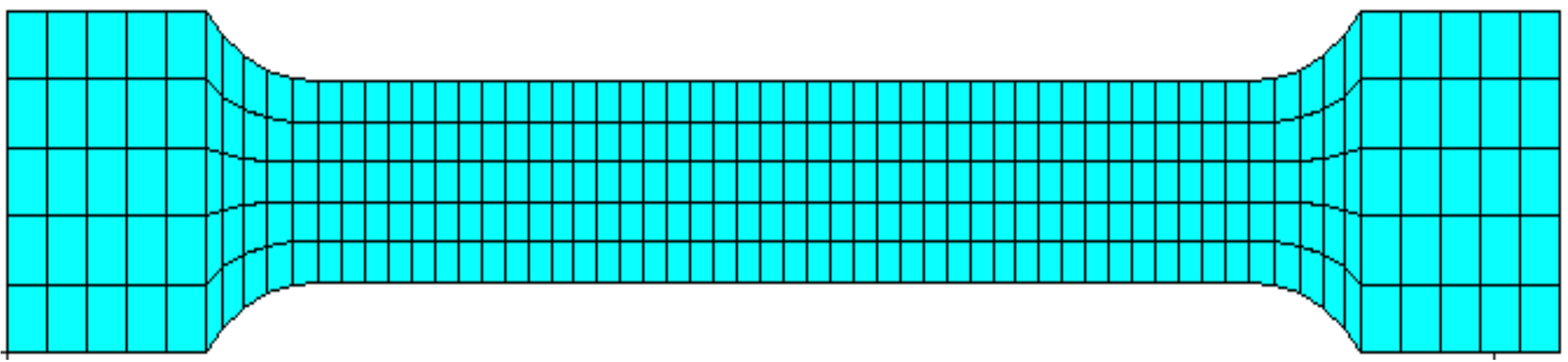} 
		\caption{\small{Mesh used in the simulations.}}
		\label{ishaped2}
	\end{subfigure}
	\caption{\small{I-shaped viscoelastic specimen.}}
\end{figure*}

This section presents numerical results for an I-shaped viscoelastic specimen, without voids, modeled as a plane stress state, whose dimensions are given in Fig. \ref{ishaped}.
The adopted finite element mesh has 300 linear squared elements, as shown in Fig. \ref{ishaped2}, and the time step is $\Delta t=10^{-3}\ s$. 
Other geometric and material parameters are thickness $t=0.132934\ mm$, Griffith coefficient $g_c=4000\ N/m$, fracture layer width
$\gamma= 0.025\ mm$, Young's modulus $E=69\times 10^9\ Pa$, Poisson's ratio ${\nu}=0.33$ and density $\rho=2700\ Kg/m^3$. 
Some of these material parameters are chosen in order to simulate a general viscoelastic hard-strong polymeric material. 
For the results presented in this section, the inspection point corresponds to the center of the specimen. 
We use the numerical fractional derivative algorithm G1 of Eq. \eqref{riet}.

\subsubsection{Loading-Unloading Test}

Firstly, we perform a loading-unloading test in order to
check the dynamic response of the motion equation. 
The specimen of Fig. \ref{ishaped} is fixed on the left end and subject to an incremental distributed load on the opposite end with rate $5.0 \times 10^6\ N /s$ until time $t=0.8\ s$. 
After that, unloading is performed with the same rate in the opposite direction. 
Damage effect is not considered. 
The tolerance of the Newton-Raphson procedure is $1\times 10^{-8}$.

\begin{figure}[!h]
\centering
	\begin{subfigure}{.33\textwidth}
		\includegraphics[scale=0.42]{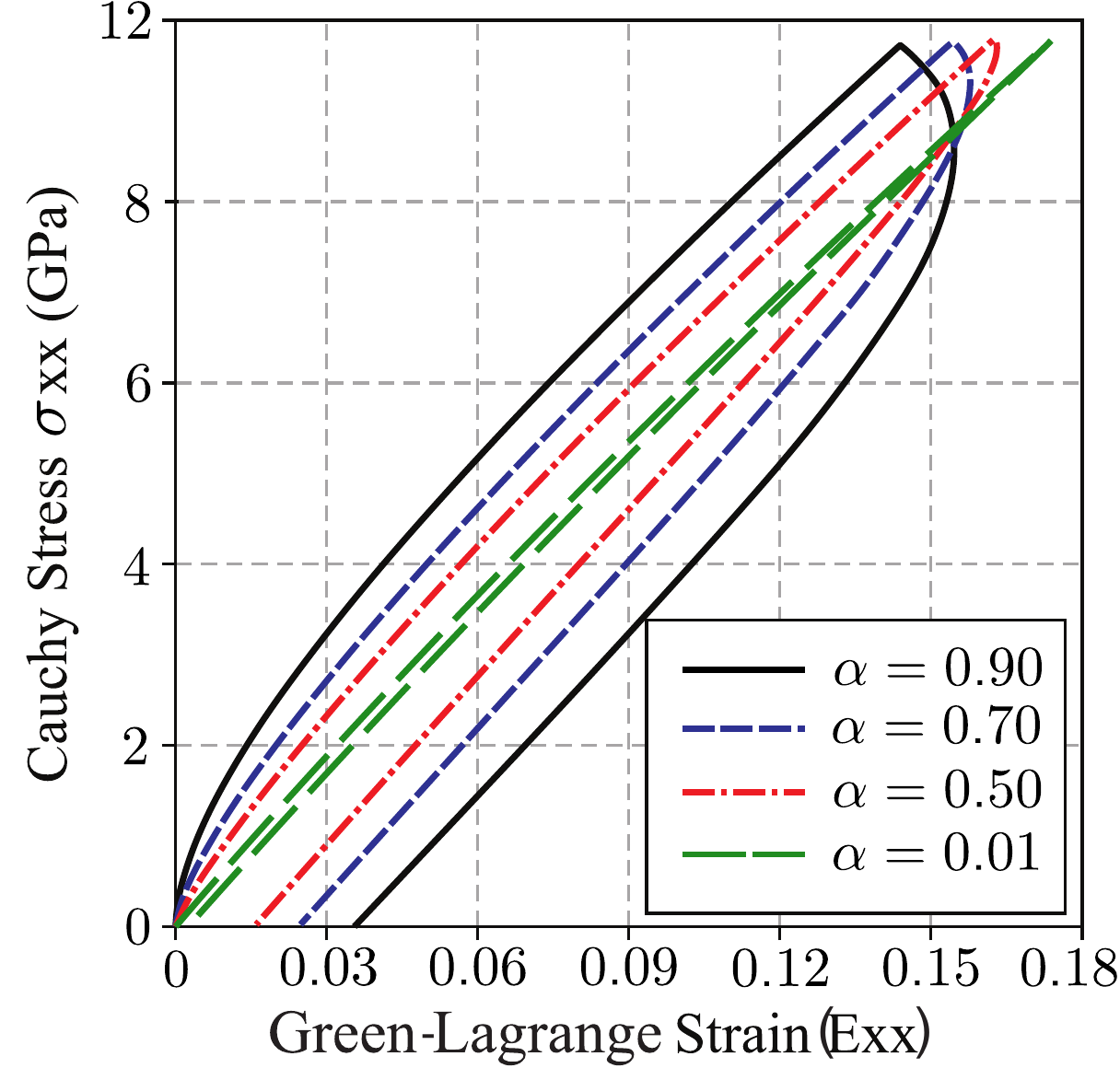} 
		\caption{\small{Load-unload test.}}\vspace{0.25cm}
		\label{load_unload}
	\end{subfigure}
	\begin{subfigure}{.33\textwidth}
		\vspace{1pt}
		\includegraphics[scale=0.42]{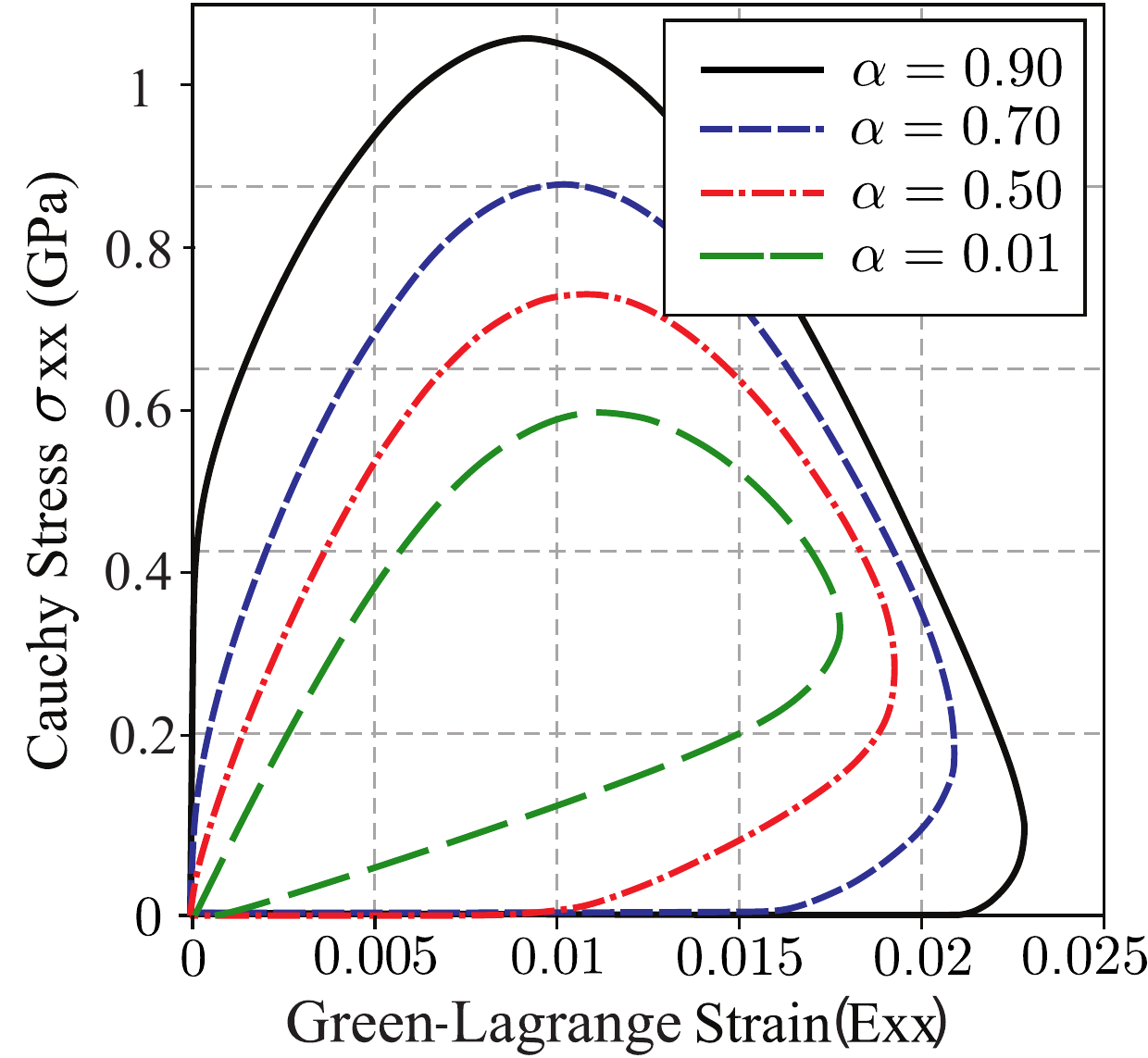} 
		\caption{\small{Tensile test until the specimen breaks.}}
		\label{with_damage}
	\end{subfigure}
	\caption{\small{Stress/strain diagram for some values of $\alpha$ in the horizontal direction.}}
\end{figure}

Figure \ref{load_unload} shows the stress/strain diagram in the horizontal direction for some values of $\alpha$ and $p=214.6\times 10^4\ N/m^2s^\alpha$. 
We remark that the residual strains shown in Fig. \ref{load_unload} were not prescribed in our model; they depend on the variation of $\alpha$.
The residual strain is larger when $\alpha$ is closer to 1. 
This was expected, because when $\alpha$ increases, the viscous effect grows and the elastic recovery decreases.
The behavior of the curves agrees with the literature for viscoelastic material under a loading-unloading process \cite{zhang1997nonlinear}. 

\subsubsection{Tensile Test with Damage Evolution}

A simple tensile test is performed for the specimen shown in Fig. \ref{ishaped}. 
In this case, the effect of the damage evolution is included and an incremental displacement $u=1.0\times 10^{-5}\ mm/t.s$ is applied until the specimen breaks. 
Most of the parameters adopted in this case are the same as the previous section, except for $p=69\times 10^8\ N/m^2s^\alpha$ and $\tilde{c}=10^{-7}\ m^2/Ns$. 
The latter corresponds to the rate of the damage increase and appears in Eq. \eqref{pseudo_d1}.
The tolerance of the Newton-Raphson procedure is $1.0\times 10^{-8}$ for the motion equation and $1.0\times 10^{-3}$ for the damage equation. 
Figure \ref{with_damage} shows the stress/strain curves for some values of $\alpha$. 
Here, the nonlinearity is different from the previous case due to the several effects considered. 
As a result of the viscoelastic behavior, when the specimen breaks the stress/strain curves return to the origin.

\subsection{Fitting with Experimental Data - Loading-Unloading Test}\label{eperimentalfit}

\noindent
In this section, we describe the fitting of experimental results by using the model proposed in Sec. \ref{development_of_the_model}. 
The experimental data are obtained from the work carried out by Dusunceli and Colak \cite{dusunceli2008effects}, who performed loading-unloading tensile tests to describe some properties of high density polyethylene (HDPE). 
The specimens used are collected from extruded PE100 pipes whose dimensions of the samples are shown in Fig. \ref{experimental_specimen} following the ISO 6259-1 and ISO 6259-3 standards.
It is important to say that we do not have access to the exact experimental points.
Then, in order to promote qualitative comparisons, we collected this information directly from \cite{dusunceli2008effects}.

The performance of our model in fitting the experimental data is tested for the case of small ($\leq 5\%$) and large strains ($> 5\%$), as presented next.

\subsubsection{Small Strain}
\label{smallStrains}

\begin{figure*}[!h]
\centering
	\begin{subfigure}{.7\textwidth}
		\includegraphics[scale=0.43]{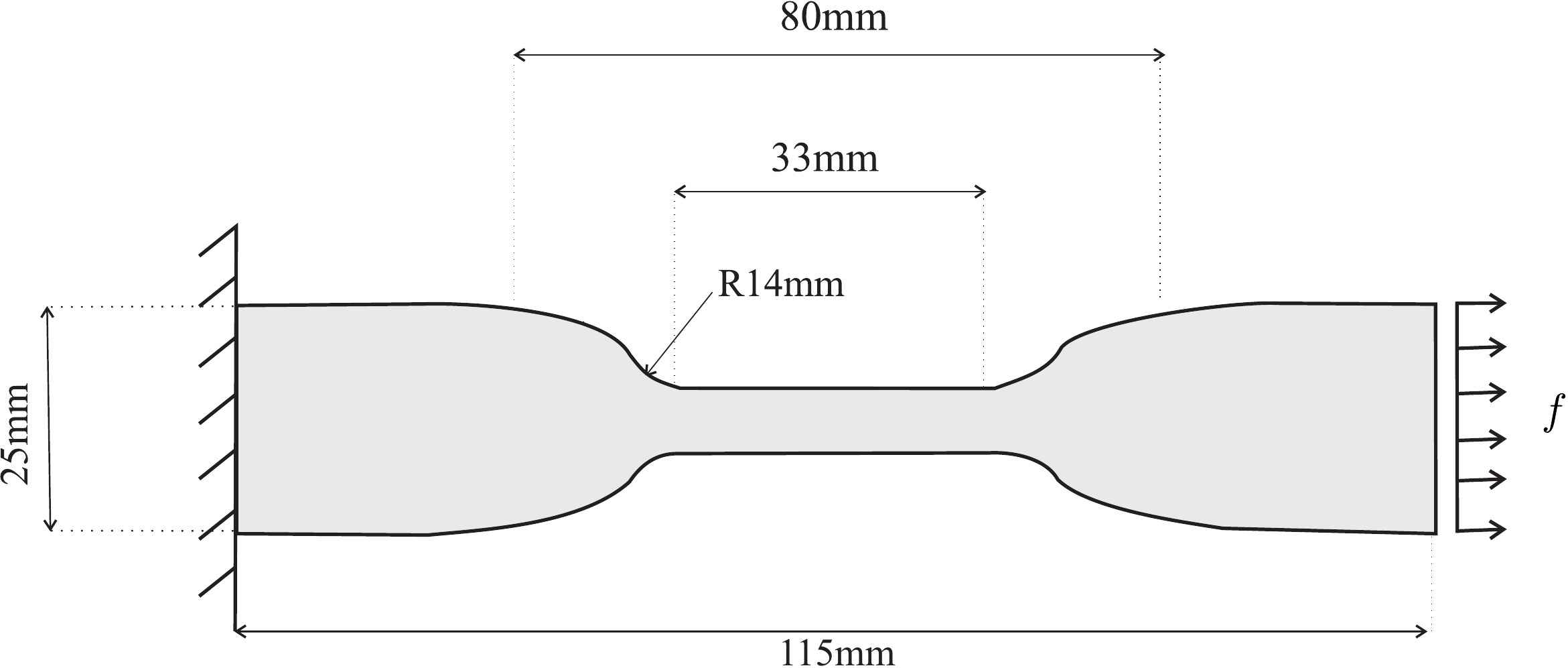} 
		\caption{\small{Dimensions of the specimen.}}
		\label{experimental_specimen}
	\end{subfigure}\\
	\begin{subfigure}{.7\textwidth}
	\hspace{1.4cm}
		\includegraphics[scale=0.42]{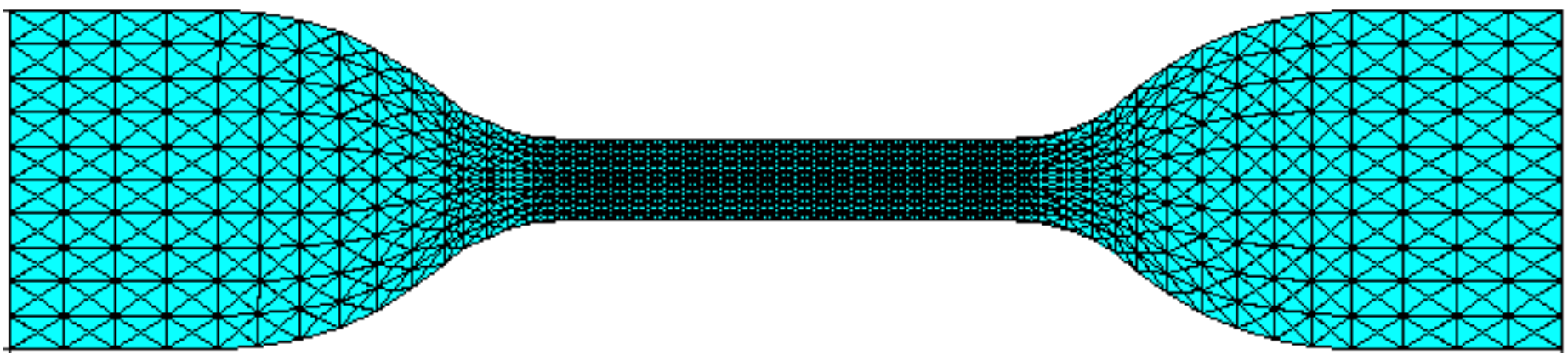} 
		\caption{\small{Mesh used.}}
		\label{experimental_specimen2}
	\end{subfigure}
	\caption{\small{Sample of HDPE used for load-unload tests.}}
\end{figure*}

\noindent
In the loading-unloading tensile test, the specimen of Fig. \ref{experimental_specimen} is fixed at one end and a uniaxial load in direction $x$ with a constant strain rate of $1\times 10^{-4}$ is applied on the other end. 
When the specimen achieves 5\% strain, an unloading is performed with the same strain rate in the opposite direction. The temperature is constant at 24$^\circ$ C.

We ran quasi-static simulations in order to reproduce the experimental procedure described above for  plane strain state.
The effects of damage are included by using Eq. \eqref{final_damage} with the degradation function $G_1$ of Eq. \eqref{functionG1}. 
The fractional derivative is calculated using the Algorithm G1 presented in Sec. \ref{fractional_numerical}. 
For the results presented in this section, the considered point corresponds to the center of the specimen. 

We perform a curve fitting based on identifying the parameters.
The effects of the variation of a particular parameter were investigated by a series of tests.
Once the influence of this parameter on the stress/strain curve is established, we checked values which lead to the intended behavior for the fitting.

It is important to emphasize that the fitting process is performed just for the loading case.
By imposing the opposite strain rate for the simulation, 
we predict the unloading results which can then be compared with the experimental unloading results.

The material parameters identified in this procedure were Young's modulus $E= 0.8\times\ 10^8Pa$; rate of the damage propagation $\tilde{c}=0.18\times 10^{-2}\ m^2/N.s$; $p=0.56\times 10^9\ N/m^2s^\alpha$; and the order of the fractional derivative $\alpha=0.3$. 
Other required constants are 
fracture toughness $f_t=0.89\times 10^6\ Pa.m^{\frac{1}{2}}$; fracture layer width $\gamma= 0.006\ mm$; Poisson's ratio ${\nu}=0.45$; density $\rho=0.954\ g/m^3$; and $\zeta=1$. 
The fracture toughness is used to calculate the Griffith constant $g_c$ using the relation \cite{franccois1998mechanical}
\begin{eqnarray}
	g_c = f_t^2 \frac{(1.0 - \tilde{\mu}^2)}{E}.
\end{eqnarray}
The tolerance of the Newton-Raphson procedure is $10^{-12}$ for both motion and damage equations.
A finite element mesh of 2240 linear triangular elements is considered as shown in Fig. \ref{experimental_specimen2} and the time step is $\Delta t=0.1\ s$. 

\begin{figure}[t]
\centering
	\includegraphics[scale=0.44]{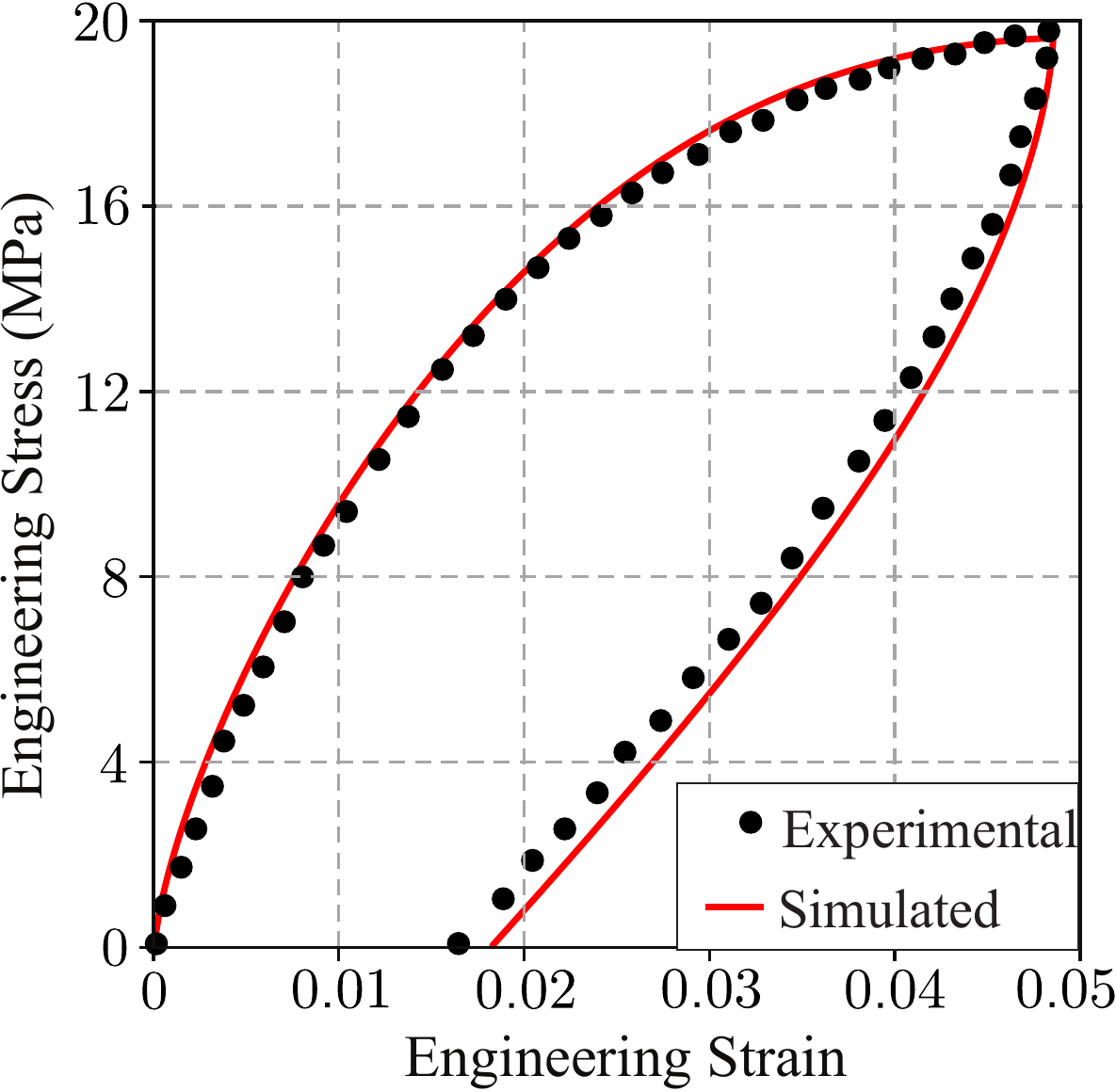}
	\caption{\small{Stress/strain relation in a load/unload test for the HDPE.  Degradation function $G_1$ was used in the case of small strain.}}
	\label{paper_fitting2} 
\end{figure}

\begin{figure}[t]
	\centering
	\includegraphics[scale=0.5]{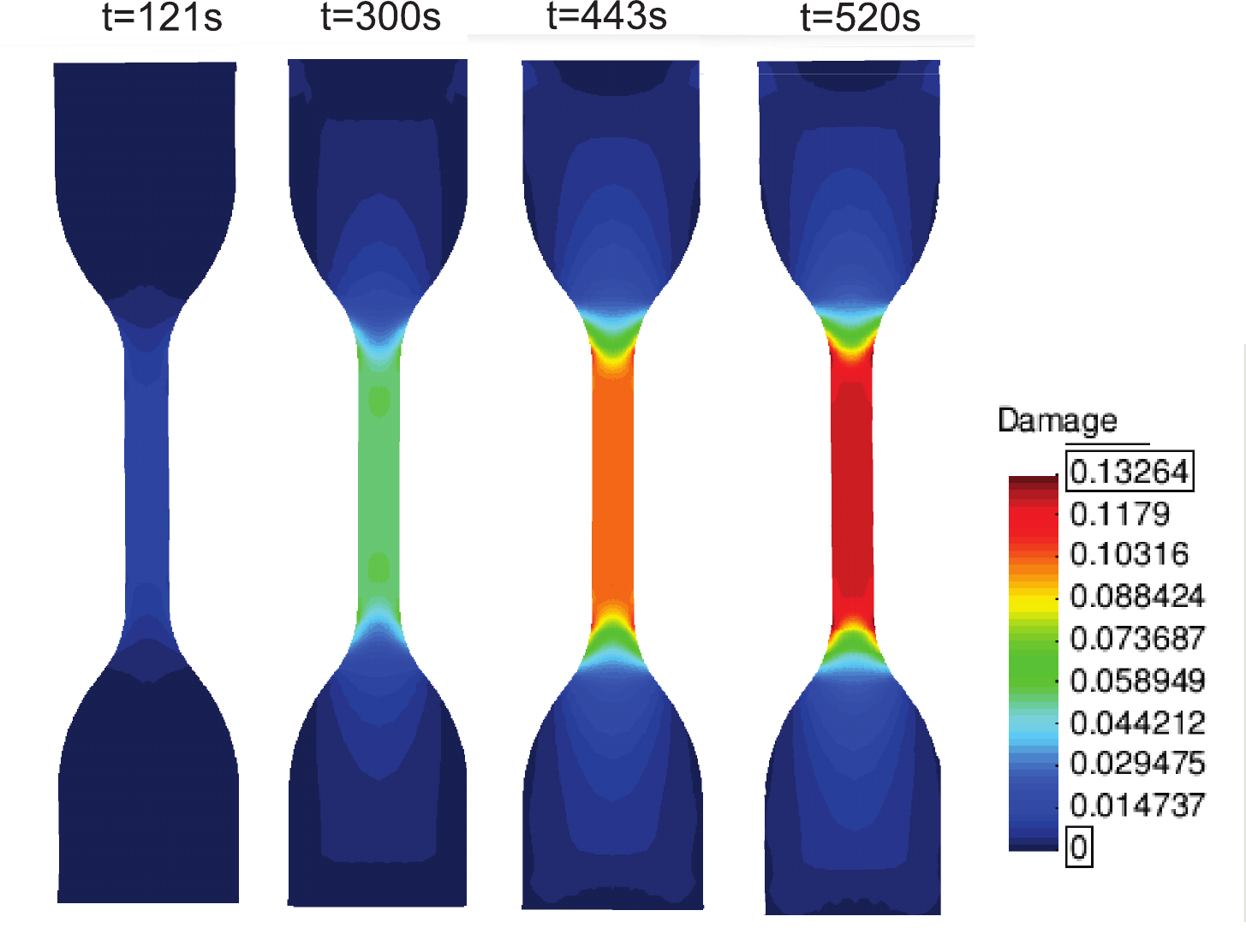} 
	\caption{\small{Damage distribution in the specimen until 5\% strain by using degradation function $G_1$.}}
	\label{damage_exp}
\end{figure}

Figure \ref{paper_fitting2} presents the comparison between the stress/strain curves obtained in the simulation and the experimental data. 
The model promoted good curve fitting for the loading process and recovered the curve pattern in the unloading. 
It is advantageous in relation to many models presented in the literature for viscoelastic materials, which are not able to recover the unloading process properly. 
{Furthermore,} the model presented in this work allows the strain process and damage evolution to be coupled.
At the end of the loading process, the level of degradation induced by the damage affects the strain in the unload process, yielding a residual strain similar with the experimental data.
Figure \ref{damage_exp} shows the damage evolution.

{We also ran the same test with $\Delta t=10^{-3}\ s$ in order to check the reliability of the analysis obtaining the same qualitative results.
We observed some numerical issues when using the fractional derivative algorithm G1 with very small values of $\Delta t$
that brings difficulties in using automatic optimization procedures for parameter identification.
Algorithm G1 is attractive due to its simple implementation, but it has high computational cost.

Since the purpose of the present work is to verify the effectiveness of the proposed viscoelastic model, 
we do not investigate these computational aspects of the fractional derivatives in the present paper.
In future works, we will consider more economic fractional derivative algorithms.}

\subsubsection{Large Strain}
\label{largestrains}

\noindent
We now consider the extension of the previous test for the case of large strain. 
The specimen shown in Fig. \ref{experimental_specimen} is subject to an uniaxial load in the $x$ direction until the sample achieves 15\% strain, then an unloading is performed with the same strain rate. 

\begin{figure}[!h]
\centering
	\begin{subfigure}{.33\textwidth}
	\includegraphics[scale=0.42]{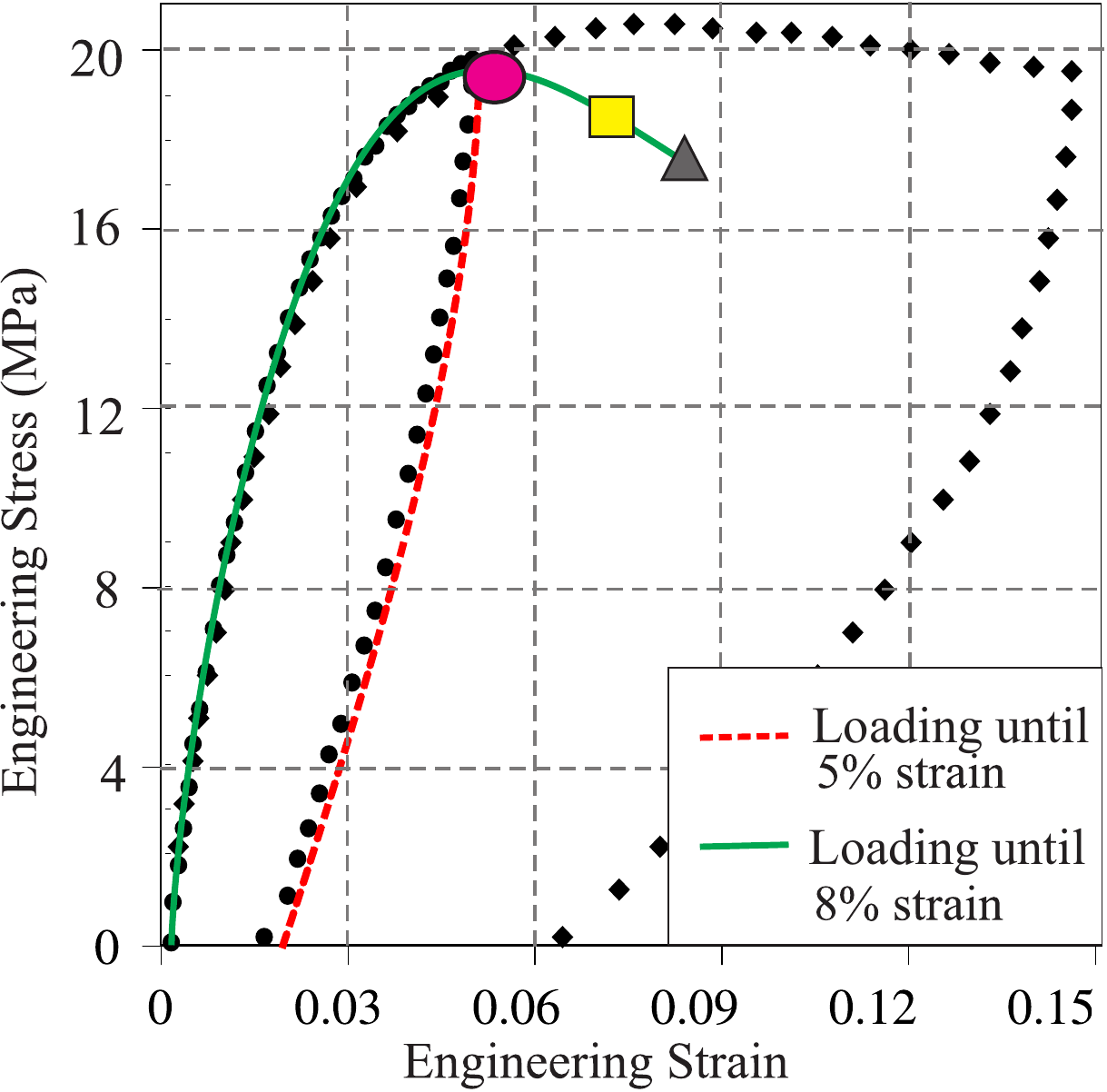} 
	\caption{\small{Stress/strain curves using the degradation function $G_1$. Dashed line represents the loading/unloading process for 5\% strain. Solid line represents the loading until 8\% strain.}}\vspace{0.25cm}
	\label{large1}
	\end{subfigure}
	\begin{subfigure}{.33\textwidth}
	\includegraphics[scale=0.42]{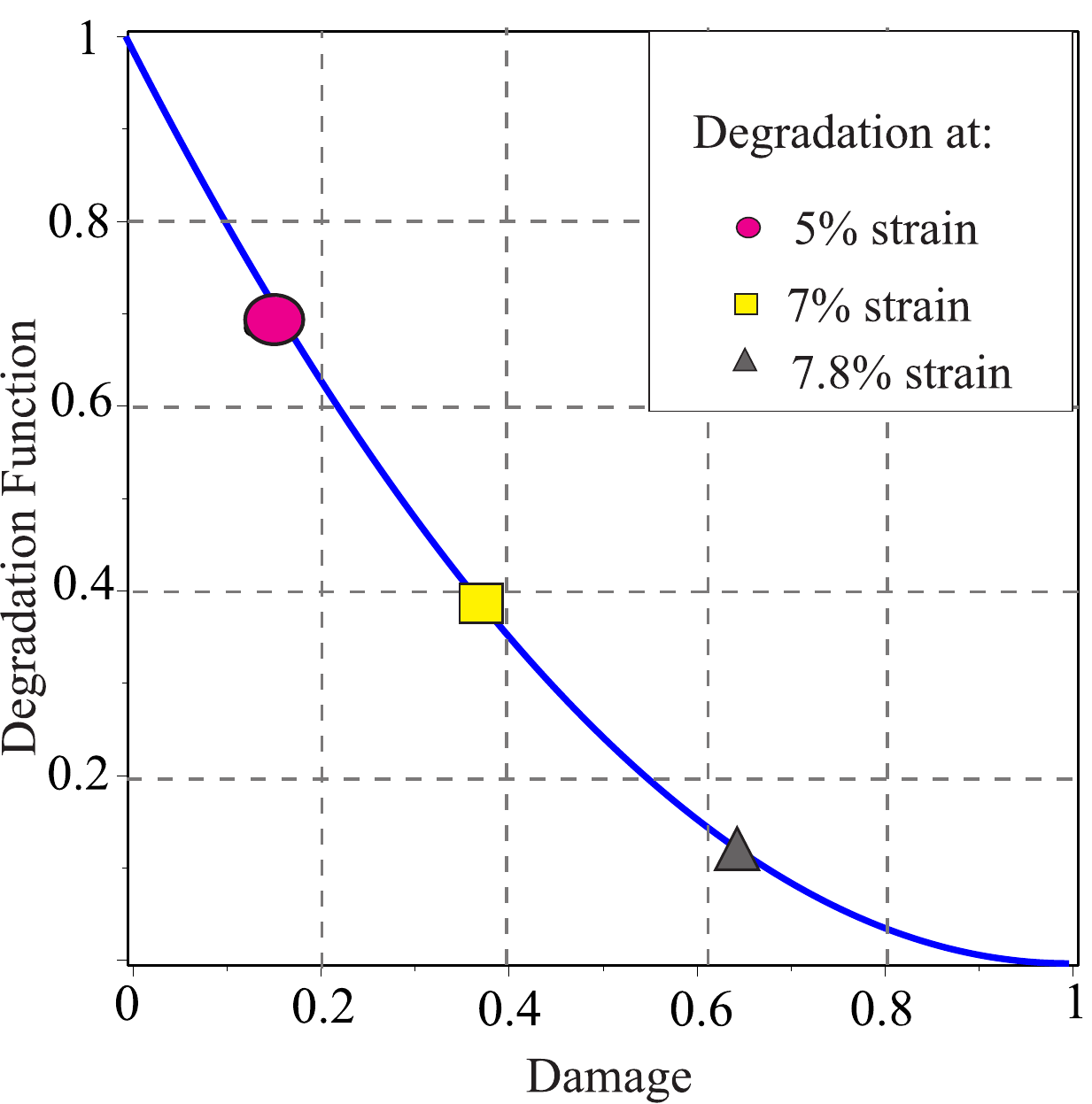} 
	\caption{\small{Degradation function $G_1$.}}
	\label{quadraticdegradfunc}
	\end{subfigure}
	\caption{\small{ Stress/strain curve in the tensile test using $G_1$ for the HDPE \ref{large1}  and the associated evolution of the degradation function $G_1$ \ref{quadraticdegradfunc}. 
	The geometric symbols correlate the percentage of strains with the corresponding damage values for the degradation function.}}
\end{figure}

We evaluated the fitting of the experimental data by using the same conditions and parameters identified in Sec. \ref{smallStrains}. 
The resulting simulated stress/strain curves can be compared with the experimental data in Fig. \ref{large1} for 5\% and 15\% strains. 
For the experimental data, the stress increases up to about 8\% strain, then it slowly decays until the unloading is performed. 
On the other hand, the stress decreases fastly after 5\% for the simulated results. 
As can be seen in Fig. \ref{quadraticdegradfunc}, the used degradation function $G_1$ decreases quickly as the process evolves. 
Once this function strongly influences the stress (see Eq. \eqref{final_second_piola}), it also decreases rapidly.
The geometric symbols presented in Figs. \ref{large1} and \ref{quadraticdegradfunc} correlate the strain and damage values for function $G_1$.

Based on these results, 
we see that the use of $G_1$ in the model does not give a correct damage behavior for viscoelastic materials in the case of large strain.
This is so because $G_1$ does not correctly describes the degradation mechanisms described in Remark \ref{Remark_viscoelastic_damage}.

To obtain the right degradation behavior,  
we used the degradation function $G_2$ of Eq. \eqref{functionG2} to perform a new fitting. 
Function $G_2$ depends on constants $a$, $b$ and $c$, which are included in the inverse parameter identification. 

The material parameters identified in this procedure were Young's modulus $E= 0.4\times 10^8\ Pa$; rate of damage propagation $\tilde{c}=0.115\times 10^{-2}\ m^2/Ns$; $p=0.67\times 10^9\ N/m^2s^\alpha$; order of  fractional derivative $\alpha=0.35$; and the degradation function parameters $a=3.8$, $b=1.5$ and $c=1.15$. 
Other  parameters are 
fracture toughness $f_t=0.89\times 10^6\ Pa.m^{\frac{1}{2}}$; length of the fracture layer width $\gamma= 0.006\ mm$; Poisson's ratio ${\nu}=0.45$; density $\rho=0.954\ g/m^3$; and $\zeta=1$. 
The same mesh, time step and tolerance for the Newton-Raphson of the preceding section were adopted in this case.

\begin{figure}[!h]
\centering
	\begin{subfigure}{.33\textwidth}
	\includegraphics[scale=0.42]{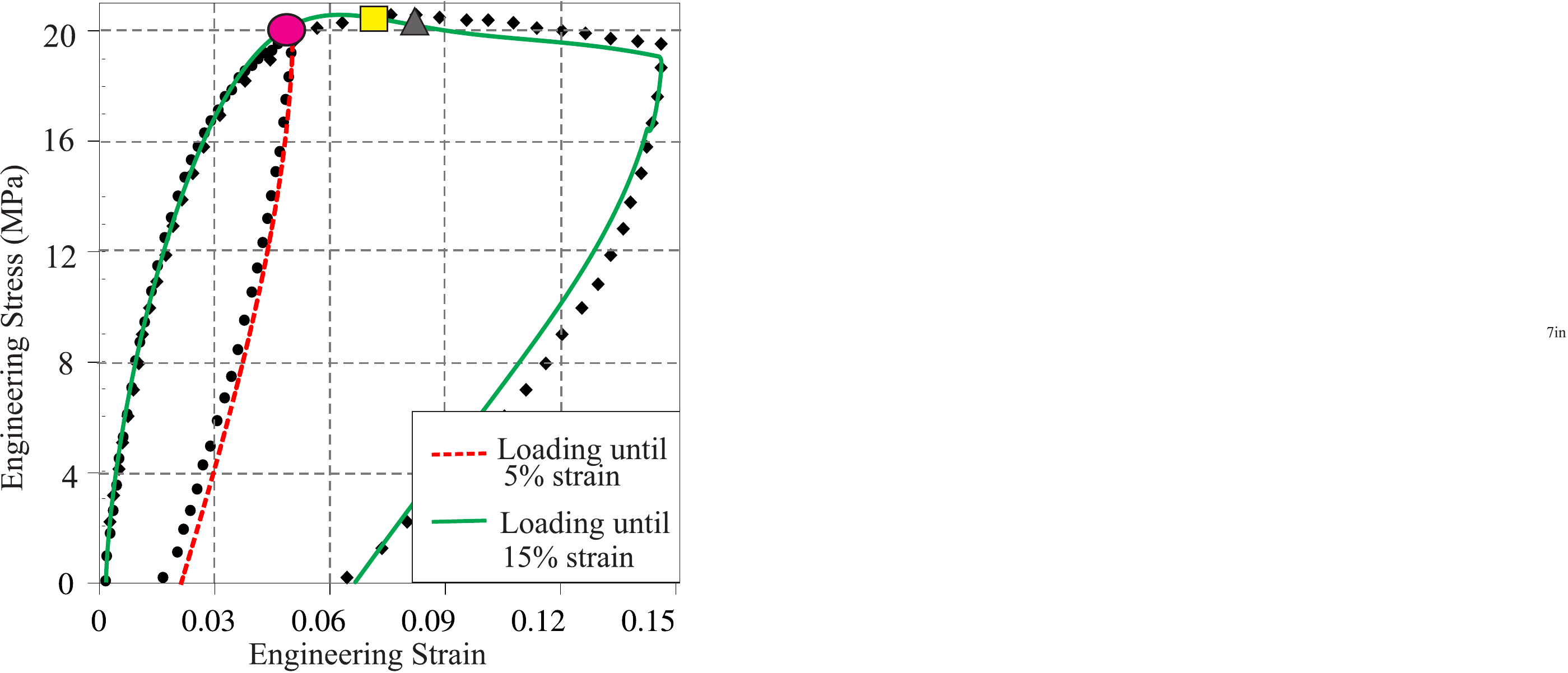} 
	\caption{\small{Stress/strain curves. Dashed line represents the loading/unloading process until 5\% strain. Solid line represents the loading until 15\% strain.}}\vspace{0.25cm}
	\label{grandes2}
	\end{subfigure}
	\begin{subfigure}{.33\textwidth}
	\includegraphics[scale=0.42]{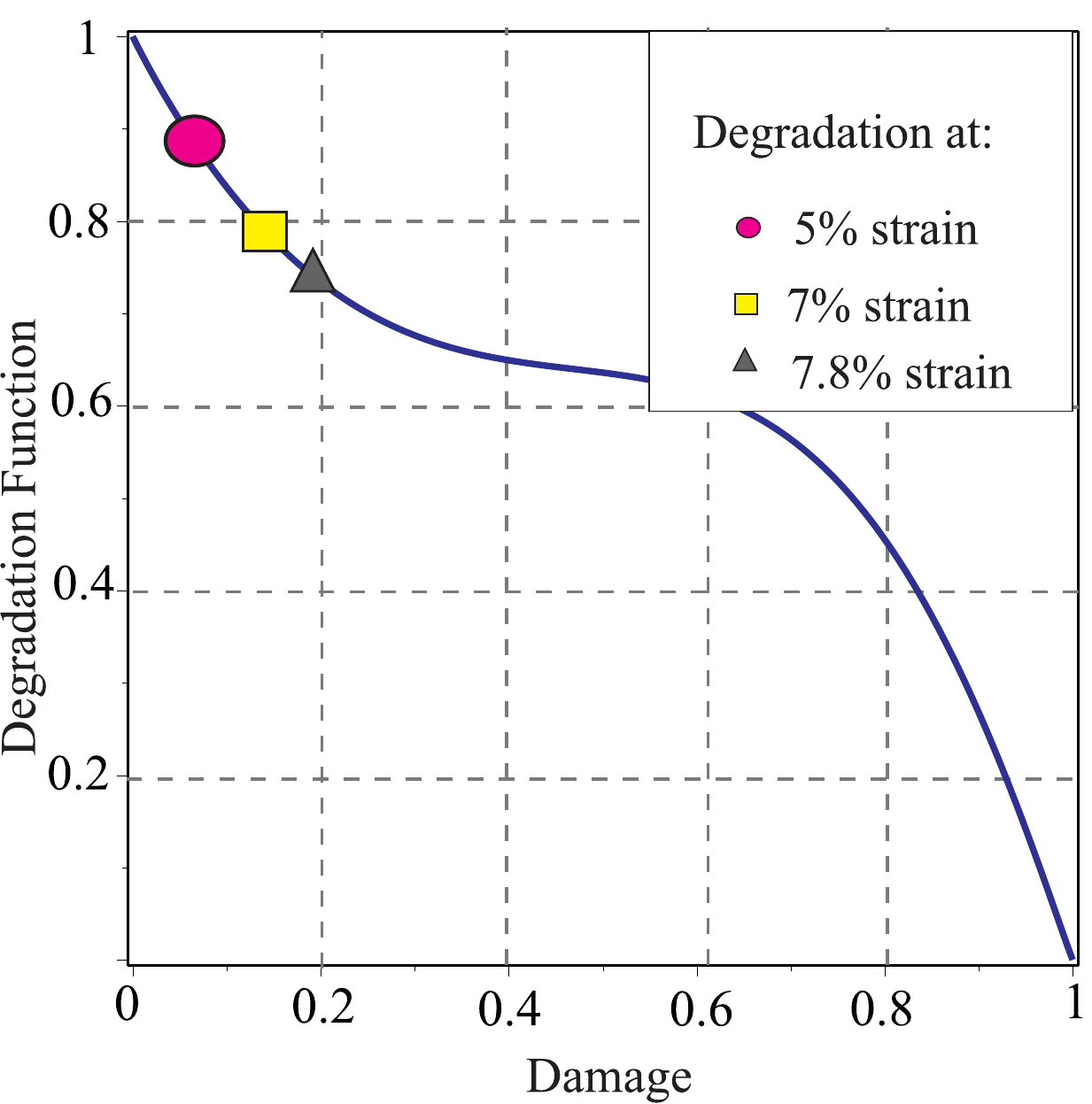} 
	\caption{\small{Degradation function $G_2$ for $a=3.8$, $b=1.5$ and $c=1.15$.}}
	\label{newdegradfunc}
	\end{subfigure}
	\caption{\small{Stress/strain relation in the tensile test for the HDPE \ref{large1} and the associated evolution of the degradation function $G_2$ \ref{newdegradfunc}. The geometric symbols correlate the percentage of strains with the damage values for the degradation function.}}
\end{figure}

\begin{figure}[h]
	\centering
	\includegraphics[scale=0.3]{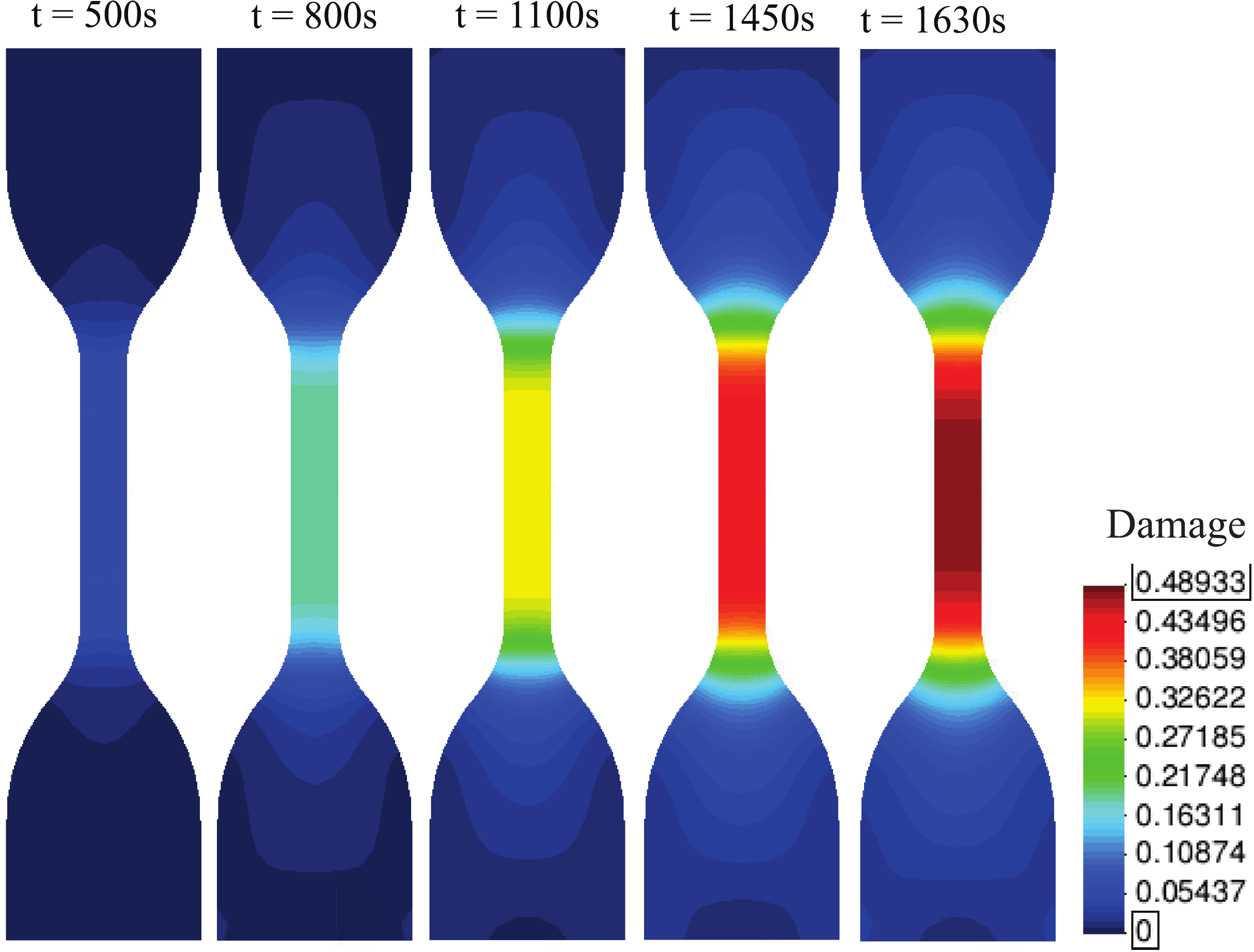} 
	\caption{\small{Damage distribution in the specimen until 15\% strain using degradation function $G_2$.}}
	\label{experimental_specimen_a}
\end{figure}

Figure \ref{grandes2} presents the new simulated stress/strain
curves.  
The new degradation function significantly improved the fitting, once the function $G_2$ was designed to describe the evolution of damage in the micro-structure (see Sec. \ref{degradation_function}). 
Figure \ref{newdegradfunc} shows the degradation function $G_2$, and the points which associate the damage values with the corresponding strain levels. 
Figure \ref{experimental_specimen_a} shows the damage evolution in the specimen for this case.	

The behavior of functions $G_1$ and $G_2$ are very similar until the damage achieves approximately $\varphi=0.16289$, which corresponds to 5\% strain.  
In fact, the fitting until $5\%$ of strain is not significantly affected by the choice of the degradation functions $G_1$ or $G_2$. 
In the unloading, both of these functions predict a rather correct level of degradation for the tested material. 
However, only the function $G_2$ is appropriate to predict results for situations where the material are subject to large strains.

In the previous simulations, we used the irreversible damage version of the model;
thus, damage does not decrease during unloading, as expected from Eq. (\ref{final_damage}).

\vspace{0.1cm}
To further check our model, we also consider another simulation, using again the specimen described in Fig. \ref{experimental_specimen}, but now keeping the tensile loading  until it breaks.

\begin{figure}[!h]
	\centering
	\includegraphics[scale=1.2]{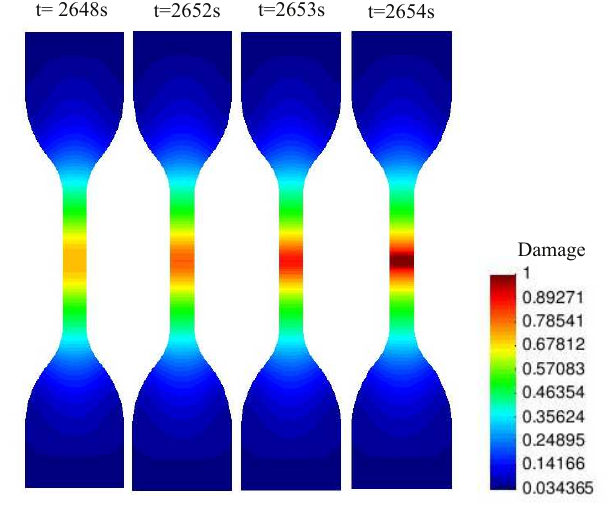} 
	\caption{\small{Damage distribution in the specimen until the fracture.}}
	\label{experimental_specimen_b}
\end{figure}

We assumed the same conditions and parameters mentioned above for the degradation function $G_2$. 
Figure \ref{experimental_specimen_b} shown the damage increasing and the localizated fracture in this case, where the strain achieves 23\% in the simulation time of 2654s.

\section{Conclusions}

We presented a general thermodynamically consistent phase-field model to describe damage {in viscoelastic materials}.
This model is constructed in Lagrangian configuration,
and damage is described by a dynamic phase-field variable.

Viscoelasticity is included in the model by using a suitable free-energy potential and a pseudo-potential of dissipation that lead to stress/strain constitutive relation in terms of fractional derivatives with finite strain
and ensure the validity of the second principle of thermodynamics 

We introduced a novel free-energy potential with memory effects and related to fractional derivatives;
this potential depends on function $ \mathcal{N}$ that can be chosen to represent different viscoelastic materials.
The free-energy potential allows damage evolution by including suitable degradation functions, which play an important role in modeling the change of stiffness between the undamaged and the fractured states. 
We proposed a new degradation function having suitable features to describe the viscoelastic behavior according to the evolution of damage in the micro-structure.

The development of the model results in a set of fractional order differential equations, which describe the evolution of motion, damage and temperature in a viscoelastic body. 
The numerical solutions of this system were obtained by using a semi-implicit/{explicit} method.

The behavior of our model was verified by numerical tests and comparisons with experimental data for the isothermal case.

To guarantee thermodynamical consistency, our model has, in addiction to an integral term related to the fractional derivatives, some extra singular integral terms.
To evaluate the importance of those extra-terms, a one-dimensional version of the model was considered to quantify their contribution to the stress evolution.
We compared the results of the simulations of the dynamic response of a viscoelastic rod subject to an external force
with and without the extra terms.

We observed that the extra terms do not significantly affect the stress response and can be disregarded in most cases. 
Afterwards, we used the same example to study the effect of varying the viscoelastic parameters on the displacement behavior.
We obtained displacement curves over time which agree qualitatively with the results of literature.
Simulated data presented for this case were used as a reference to obtain a suitable two-dimensional extension. 
Subsequently, two tests for an I-shaped specimen were performed in plane the stress state. 
One of these introduced the damage evolution, resulting in stress/strain curves for a tensile test.

Fitting of experimental data was performed to describe the viscoelastic response of HDPE samples in loading/unloading tensile tests for two cases: small strains (5\% strain) and large strains (15\% strain). 
Fittings were done just by using loading data; the predict unloading were then compared with the experimental unloading results.
We obtained results with the usual quadratic degradation function $G_1$ and also with function $G_2$ proposed in this work. 
The predicted results were both adequate for the case of small strains.
However, for large strain, $G_1$ gave incorrect results, while $G_2$ gave rather good results.

We observed numerical limitations related to the fractional derivative algorithm G1 when used with very small values for $\Delta t$
that made difficult the implementation of automatic optimization procedures for parameter identification.
Algorithm G1 is attractive due to its simple implementation, however, it requires a very high computational time.
Since the purpose of the present work is to verify the effectiveness of the proposed viscoelastic model, 
we did not investigate here the computational aspects of the fractional derivatives.
Future work will consider more economical algorithms for fractional derivatives.

Finally, the results presented indicate that the proposed model is successful in describing the response of viscoelastic materials under the conditions tested.
It is also an adequate thermodynamically consistent alternative to account for the viscoelastic behavior under damage.

\begin{acknowledgements}
The authors would like to thank the Coordination for the Improvement of Higher Education Personnel (CAPES) and the S\~ao Paulo Research Foundation (FAPESP), under grant 2015/20188-0, for their financial support.
\end{acknowledgements}

\appendix
\section{Conditions on $\mathcal{N}$ and Derivativation of $\dot{{\psi}}_m$}\label{AppendixA}

Consider  $\psi_m$ and  $\tilde{\psi}_m$ given respectively by Eqs. \eqref{PotentialMemoryEffectsIncludingDamage} and \eqref{FirstFormMemoryEffects}.
Function $\mathcal{N}( \bm{Z}_1,\bm{Z}_2 )$ of Eq. \eqref{FirstFormMemoryEffects} is a suitable {continuous} function of second-order symmetric tensors with the following properties:
\vspace{0.2cm}
\begin{enumerate}
	\item[(a)] $\mathcal{N} (\bm{Z}_1,\bm{Z}_2) \geq 0$, $\forall\quad \bm{Z}$;
	\item[(b)] $|\mathcal{N}(\bm{Z}_1, \bm{Z}_2)| \leq C(\bm{Z}_1,  \bm{Z}_2) \|\bm{Z}_1 - \bm{Z}_2 \|^\beta, $ with $\beta \geq 1+\alpha$ and
$C(\bm{Z}_1, \bm{Z}_2)$ bounded as  $\bm{Z}_1 - \bm{Z}_2 \rightarrow 0_+$;
	\item[(c)] $\| \partial_{\bm{Z}_1} \mathcal{N}( \bm{Z}_1, \bm{Z}_2 ) \| \leq  C_1(\bm{Z}_1,  \bm{Z}_2) \|\bm{Z}_1 - \bm{Z}_2 \|^{\beta_1},$ \mbox{with} $\beta_1 \geq \alpha $ \mbox{and} $C_1(\bm{Z}_1, \bm{Z}_2)$ \mbox{bounded as} $\bm{Z}_1 - \bm{Z}_2 \rightarrow 0_+$.
\end{enumerate}

From the previous properties and using the mean value theorem, we obtain
\begin{eqnarray}
	\label{Inequality_Mean_value_theorem}
	&&| \mathcal{N}( \bm{E}_t , \bm{E}_\tau)  | \leq C(\bm{E}_t , \bm{E}_\tau) \|  \bm{E}_t -\bm{E}_			\tau  \|^\beta \nonumber \\
	&&\leq C(\bm{E}_t , \bm{E}_\tau) \max\{  \|\dot{\bm{E}}_s \|, s \in [0,t) \}^\beta | t - \tau  |^			\beta ,
\end{eqnarray}
and
\begin{eqnarray}
	&&\| \partial_{Z_1} \mathcal{N}( \bm{E}_t , \bm{E}_\tau) \| 
	\leq 
	C_1(\bm{E}_t , \bm{E}_\tau) \|  \bm{E}_t -\bm{E}_\tau  \|^{\beta_1}  \nonumber\\
	&& \leq C_1(\bm{E}_t , \dot{\bm{E}}_\tau) \max\{  \|\bm{E}_s \|, s \in [0,t) \}^{\beta_1} | t - \tau  	|^{\beta_1}.
\end{eqnarray}
By using inequality \eqref{Inequality_Mean_value_theorem}, we have
\begin{eqnarray}
	\lim_{\tau \rightarrow t-} 
	\frac{ \mathcal{N}( \bm{E}(\bm{p}, t), \bm{E}(\bm{p}, \tau) )}{(t - \tau)^{1+\alpha}} = 0, \quad 0<	\alpha<1.
\end{eqnarray}
Moreover, we also consider that
$\displaystyle \frac{ \mathcal{N}( \bm{E}(\bm{p}, t),\bm{E}(\bm{p},\tau) )}{(t -\tau)^{2+\alpha}}$ and
$\displaystyle \frac{ \partial_{Z_1} \mathcal{N}( \bm{E}(\bm{p}, t), \bm{E}(\bm{p},\tau) )}{(t -\tau)^{1+\alpha}},$ 
are integrable in $[0,t)$ with respect to $\tau$.

Under the previous conditions for $\mathcal{N}$, 
the time derivative $\dot{{\psi}}_m$ for strains such are continuous at time $t = 0+$ and have bounded rates (i.e, $\| \dot{\bm{E}} (\bm{p}, t)\|$ bounded as $t \rightarrow 0_+)$, can be obtained as
\begin{eqnarray}
	\label{TimeDerivativeFormMemoryEffects}
	\dot{\psi}_{m}(\varphi, \Ha(\bm{E}))  
   	&=& \frac{G_{m}(\varphi )}{\rho}	\dot{\tilde{\psi}}_{m}(\Ha(\bm{E}))
  	+ \frac{G'_{m}(\varphi )}{\rho} \tilde{\psi}_{m}(\Ha(\bm{E})) \dot{\varphi}\nonumber\\
   	&=&  \bm{S}_m  : \dot{\bm{E}}    + \frac{G'_{m}(\varphi )}{\rho} \tilde{\psi}_{m}(\Ha(\bm{E})) 				\dot{\varphi} - R,	
\end{eqnarray}
where
\begin{eqnarray}
	\bm{S}_m 
	&=&  \frac{G_{m}(\varphi )}{\rho\Gamma(1-\alpha)}
	\bigg[\frac{  \partial_{\bm{E}_t} \mathcal{N}( \bm{E}_t, \bm{E}_0 ) }{t^{\alpha}} 
	+ \alpha \int_0^{t}\frac{ \partial_{\bm{E}_t} \mathcal{N}( \bm{E}_t, \bm{E}_\tau )}{(t -\tau)^{1+			\alpha}}\ \diff\tau\bigg],
\end{eqnarray}
and
\begin{eqnarray}
	R &=&  \frac{G_{m})\alpha}{\rho\Gamma(1-\alpha)} \bigg[   \frac{ \mathcal{N}(\bm{E}_t,\bm{E}_0 )}{t ^{1+\alpha}}   
  	+ (1+\alpha) \int_0^{t}\frac{ \mathcal{N}( \bm{E}_t,\bm{E}_\tau )}{(t -\tau)^{2+\alpha}}\ \diff			\tau \bigg].
\end{eqnarray}
Note that $R \geq 0$, due to the property (b) defined previously for $\mathcal{N}$.

\section{Examples for $\mathcal{N}$}
\label{Examples_for_N}

\vspace{0.2cm}
\noindent
{\bf Example 1:}

By considering $\mathcal{N} (Z_1, Z_2) = \frac{1}{2} (Z_1 - Z_2):\mathbfcal{{A}}:(Z_1 - Z_2)$,  with $\mathbfcal{{A}}$ a fourth order symmetric-positive definite constitutive tensor, it is easy to check that the conditions for $\mathcal{N}$ given in Appendix \ref{AppendixA} are satisfied. 
Replacing it in (\ref{FirstFormMemoryEffects}), we have
\begin{eqnarray}
	\label{pseudod2}
	\tilde{\psi}_{m} (\Ha(\bm{E}))(\bm{p},t)  
	&=&\frac{\kappa}{\rho}\bigg[\frac{\left[\bm{E}_t-\bm{E}_0\right)]:\mathbfcal{{A}} : \left[\bm{E}_t-			\bm{E}_0\right]}{t^{\alpha}}\nonumber\\
	&&+ \alpha \int_0^{t}\frac{\left[\bm{E}_t-\bm{E}_\tau\right]: \mathbfcal{{A}} : \left[\bm{E}_t-E_			\tau\right]}{(t - \tau)^{1+\alpha}}\ \diff s\bigg],
\end{eqnarray}
where $\kappa=1/2\Gamma(1-\alpha)$ and $\Gamma$ is the standard Gamma function \cite{artin2015gamma}.
In this case, we obtain
\begin{eqnarray}
	\bm{S}_m 
	= \frac{G_{m}(\varphi )}{\rho}\left( \mathbfcal{{A}} :{_0\mathrm{D}_t}^\alpha (\bm{E}_t)\right) 
	= \frac{G_{m}(\varphi )}{\rho} \left( {_0\mathrm{D}_t}^\alpha (\bm{E}_t ) : \mathbfcal{{A}})\right),
\end{eqnarray}
where ${_0\mathrm{D}_t}^\alpha \bm{E}$ is the Caputo fractional derivative of $\bm{E}$.

\vspace{0.2cm}
\noindent
{\bf Example 2:}

Now, we take $\mathcal{N} (Z_1, Z_2) = \frac{1}{2} ( \mathcal{N}_1 (Z_1) - \mathcal{N}_1 (Z_2) ):\mathbfcal{{A}}:(\mathcal{N}_1 (Z_1) -  \mathcal{N}_1 (Z_2))$ in (\ref{FirstFormMemoryEffects}),
where $\mathbfcal{{A}}$ is as in the previous example and $\mathcal{N}_1( Z ) $ is a suitable second-order-tensor valued function which satisfies the additional condition
\begin{itemize}
\item[(d)] $|\mathcal{N}_1(Z_1)  - \mathcal{N}_1( Z_2)| \leq C(Z_1, Z_2) \|Z_1 - Z_2 \|^{\tilde\beta}$, but now with $\tilde\beta \geq 1$ and $C(Z_1, Z_2)$ bounded as $Z_1, Z_2 \rightarrow 0+$.
\end{itemize}

Then,  (\ref{FirstFormMemoryEffects}),  becomes
\begin{eqnarray}
\label{SecondFormMemoryEffects}
	\tilde{\psi}_{m} 
	&=& \kappa
	\Bigg[ \frac{ [\mathcal{N}_1(\bm{E}_t)- \mathcal{N}_1(\bm{E}_0) ]: \mathbfcal{{A}} : [\mathcal{N}			_1(\bm{E}_t)- \mathcal{N}_1(\bm{E}_0) ]}{t^{\alpha}} \nonumber\\
	&&  + \alpha \int_0^{t}\frac{ [\mathcal{N}_1(\bm{E}_t)- \mathcal{N}_1(\bm{E}_\tau) ]: 						\mathbfcal{{A}} : [\mathcal{N}_1(\bm{E}_t)- \mathcal{N}_1(\bm{E}_\tau) ] )}{(t -\tau)^{1+\alpha}}\ 			\diff\tau \Bigg],\nonumber\\
\end{eqnarray}

For the special choice of $\mathcal{N}$ in this case and the previous property (d), one can easily either prove that the required properties (a), (b) and (c) are satisfied or observe directly that
\begin{eqnarray}
	&& \| \mathcal{N}_1( \bm{E}_t) -\mathcal{N}_1( \bm{E}_ \tau)  \| 
	 \leq C(\bm{E}_t) , \bm{E}_\tau) \|  \bm{E}_t -\bm{E}_\tau  \|^{\tilde\beta} \nonumber\\
	&& \leq C(\bm{E}_t , \bm{E}_\tau) \max\{  \|\bm{E}_s \|, s \in [0,t) \} | t - \tau  |^{\tilde\beta},
\end{eqnarray}
and so
\begin{eqnarray}
	&&\| [\mathcal{N}_1(\bm{E}_t)- \mathcal{N}_1(\bm{E}_t)) ]: \mathbfcal{{A}} : [\mathcal{N}_1(\bm{E}			_t)- \mathcal{N}_1(\bm{E}_t) ] )\| \nonumber\\
	&&\leq \| \mathbfcal{{A}} \| C^2(\bm{E}_t) , \bm{E}_\tau) \max\{  \|\bm{E}_s \|, s \in [0,t) \}^2 | 		t - \tau  |^{2 \tilde\beta} .
\end{eqnarray}
Thus,  we have that
\begin{eqnarray}
	\lim_{\tau \rightarrow t-}\frac{ [\mathcal{N}_1(\bm{E}_t)- \mathcal{N}_1(\bm{E}_\tau) ]: 					\mathbfcal{{A}} : [\mathcal{N}_1(\bm{E}_t)- \mathcal{N}_1(\bm{E}_\tau) ] }{(t - \tau)^{1+\alpha}} 
	= 0 .
\end{eqnarray}

Again we obtain relation (\ref{TimeDerivativeFormMemoryEffects}), 
but now with
\begin{eqnarray}
	\bm{S}_m    
	&=& \frac{G_{m}(\varphi )}{\rho} \mathbfcal{{A}} :\frac{1}{\Gamma(1-\alpha)}\left[ \frac{ 					[\mathcal{N}_1(\bm{E}_t- \mathcal{N}_1(\bm{E}_0 ]}{t^{\alpha}} \right.\nonumber\\
	&& \left. + \alpha \int_0^{t}\frac{ [\mathcal{N}_1(\bm{E}_t- \mathcal{N}_1(\bm{E}_s)])}{(t -				\tau)^{1+\alpha}}d\tau \right] : \partial_Z \mathcal{N}_1(\bm{E}_t\nonumber\\
 	&& = \displaystyle \frac{G_{m}(\varphi )}{\rho} \mathbfcal{{A}} :{_0\mathrm{D}_t}^\alpha ( 					\mathcal{N}_1(\bm{E}_t)): \partial_Z \mathcal{N}_1(\bm{E}_t) .
\end{eqnarray}

\vspace{0.2cm}
\noindent
{\bf Example 3:}

Now we take $\mathcal{N} (Z_1, Z_2) = \frac{1}{2} (Z_1 - Z_2):\mathbfcal{{A}}(Z_1):(Z_1 - Z_2)$
with $\mathbfcal{{A}} (Z_1)$ a fourth order symmetric-positive definite tensor continuously depending on $Z_1$.
Then (\ref{FirstFormMemoryEffects}) becomes
\begin{eqnarray}
	\label{pseudod2}
	\tilde{\psi}_{m}
	&=&\frac{\kappa}{\rho}\left[\frac{\left[\bm{E}_t-\bm{E}_0\right)]: \mathbfcal{{A}}: \left[\bm{E}_t-			\bm{E}_0\right]}{t^{\alpha}}\right. \nonumber\\
	&&\left. + \alpha \int_0^{t}\frac{\left[\bm{E}_t-\bm{E}_\tau\right]: \mathbfcal{{A}}: \left[\bm{E}_t-\bm{E}		_\tau\right]}{(t - \tau)^{1+\alpha}}\ \diff \tau\right].
\end{eqnarray}

\noindent
In this case, we obtain
\begin{eqnarray}
\label{Eq_Sm_appendix}
\bm{S}_m &= & \frac{G_{m}}{\rho} \left[  \mathbfcal{{A}} :{_0\mathrm{D}_t}^\alpha (\bm{E}_t) \right.
\nonumber\\
&&\left. + {\alpha\kappa} \int_0^{t}\frac{\left[\bm{E}_t-\bm{E}_\tau\right]: \partial_{\bm{E}} \mathbfcal{{A}}: \left[\bm{E}_t-E_\tau\right]}{(t - \tau)^{1+\alpha}}\ \diff \tau \right].
\end{eqnarray}
As before, ${_0\mathrm{D}_t}^\alpha \bm{E}$ is the Caputo fractional derivative of $\bm{E}$.

An interesting possibility is to take $ \mathbfcal{{A}} (\bm{E})  = \partial^2_{\bm{E}}\psi_e(\bm{E})$, where $\psi_e(\bm{E})$  is any standard elastic specific free-energy with continuous derivatives with respect to $\bm{E}$ up to order 3.
 
It is important to emphasize that Equation \eqref{pseudod2} is a modification of the free-energy potential proposed by 
Fabrizio \cite{fabrizio2014fractional}; in that work, 
the author shows that his proposal for the free-energy implies in a stress equation in terms of fractional derivatives.
However, the arguments presented in \cite{fabrizio2014fractional} do not make clear why the definition of fractional derivatives must appear.
In the present paper, we modified Fabrizio's suggestion including the first term of Eq. \eqref{pseudod2} to properly lead to the fractional derivative definition that appears in the associated stress $\mathbf{S}_M$ (see Eq. \eqref{Eq_Sm_appendix}).
We also extended his suggestion for the three-dimensional case, and added the possibility to consider $\mathbfcal{{A}}(E)$ nonlinear in relation to $E$.

\vspace{0.2cm}
\noindent
{\bf Example 4:}

Another possibility is to take $\mathcal{N}_1 (Z) = \psi_e(Z)$  in (\ref{FirstFormMemoryEffects}),
 where now $\psi_e(Z)$ is again a standard elastic specific free-energy but now normalized such that $\psi_e(Z) \geq 0$ for all $Z$, $\psi_e (0) = 0$ and with continuous derivatives with respect to $\bm{E}$ up to order 2.
Such conditions ensure that the required properties (a), (b) and (c) are satisfied.
In this case,  we then are left with
 \begin{eqnarray}
	\tilde{\psi}_{m}
	= \frac{1}{\Gamma(1-\alpha)}
	\left[
	\frac{  \psi_e ( \bm{E}_t -\bm{E}_0 ) }{t^{\alpha}} 
  	+\alpha \int_0^{t}\frac{ \psi_e( \bm{E}_t -\bm{E}_\tau )}{(t-\tau)^{1+\alpha}}\ \diff\tau
  	\right],
\end{eqnarray}

\noindent
and
\begin{eqnarray}
	\bm{S}_m 
	&=&  \frac{G_{m}}{\rho\Gamma(1-\alpha)} 
	\left[\frac{\partial_{\bm{E}_t} \psi_e( \bm{E}_t -\bm{E}_0 ) }{t^{\alpha}}\right. \nonumber\\
	&&\left. + \alpha \int_0^{t}\frac{ \partial_{\bm{E}_t} \psi_e ( \bm{E}_t -\bm{E}_\tau )}{(t -				\tau)^{1+\alpha}}d\tau\right]. 
\end{eqnarray}

\section{Alternative Expression for the Caputo Fractional Derivative}

Caputo \cite{caputo1966} proposed a fractional derivative definition for a function $f(t)\in C[a,b]$ and $a<t<b$ given by
\begin{eqnarray}
	_{a}{\mathrm{D}}_{t}^\alpha f(t) 
	= \frac{1}{\Gamma(1-\alpha)}\int_a^t\frac{{f^{(m)}}(\tau)}{(t-\tau)^{\alpha}}\ d \tau,
\end{eqnarray}
where $m=\lceil \alpha \rceil$ (ceiling function) such that $\alpha \in \mathbb{R}$ and $\Gamma$ is the usual Gamma function defined by
\begin{eqnarray}
\label{gammadefinition1}
	\Gamma(c)=\int_0^\infty e^{-\tau}\tau^{c-1}\ \diff \tau,
\end{eqnarray} 
with $c\in \mathbb{R}$. 
If $\alpha \in [0,1]$, then the particular definition is obtained
\begin{eqnarray}
	\label{caputofinal}
	_{a}{\mathrm{D}}_{t}^\alpha f(t) 
	= \frac{1}{\Gamma(1-\alpha)}\int_a^t\frac{{f'}(\tau)}{(t-\tau){^{\alpha}}}\ \diff \tau.
\end{eqnarray}

Equation \eqref{caputofinal} may be rewritten as
\begin{eqnarray}
	\label{caputo_new}
	_{a}{\mathrm{D}}_{t}^\alpha f(t) 
	= \frac{1}{\Gamma(1-\alpha)}\lim_{\epsilon \rightarrow t^-}	\underbrace{{\int_a^\epsilon\frac{{f'}			(\tau)}{(t-\tau)^{\alpha}}}\ \diff \tau}_{I}.
\end{eqnarray}
Applying integration by parts, we obtain
\begin{eqnarray}
	I&=&\frac{f(\tau)}{(t-\tau)^\alpha}\bigg|_a^\epsilon - \alpha \int_a^\epsilon\frac{f(\tau)}{(t-					\tau)^{\alpha+1}}\ \diff \tau.\nonumber\\
	&=&\frac{f(\epsilon)}{(t-\epsilon)^\alpha}-\frac{f(a)}{(t-a)^\alpha}- \alpha \int_a^\epsilon					\frac{f(\tau)}{(t-\tau)^{\alpha+1}}\ \diff \tau \nonumber\\
	&&+ \alpha \int_a^\epsilon\frac{f(\epsilon)}{(t-\tau)^{\alpha+1}}\ \diff \tau - \alpha \int_a^\epsilon\frac{f(\epsilon)}{(t-\tau)^{\alpha+1}}\ \diff \tau\nonumber\\
	&=&\frac{f(\epsilon)}{(t-\epsilon)^\alpha}-\frac{f(a)}{(t-a)^\alpha}+ \alpha \int_a^\epsilon					\frac{f(\epsilon)-f(\tau)}{(t-\tau)^{\alpha+1}}\ \diff \tau \nonumber\\
	&&- \alpha f(\epsilon)\int_a^\epsilon\frac{1}{(t-\tau)^{\alpha+1}}\ \diff \tau\nonumber\\
	&=&\frac{f(\epsilon)}{(t-\epsilon)^\alpha}-\frac{f(a)}{(t-a)^\alpha}+ \alpha \int_a^\epsilon					\frac{f(\epsilon)-f(\tau)}{(t-\tau)^{\alpha+1}}\ \diff \tau \nonumber\\
	&&-\frac{f(\epsilon)}{(t-\epsilon)^\alpha}+\frac{f(\epsilon)}{(t-a)^\alpha}\nonumber\\
	&=&-\frac{f(a)}{(t-a)^\alpha}+ \alpha \int_a^\epsilon\frac{f(\epsilon)-f(\tau)}{(t-\tau)^{\alpha+1}}		\ \diff \tau +\frac{f(\epsilon)}{(t-a)^\alpha}.
\end{eqnarray}
By replacing the above expression in Eq. \eqref{caputo_new}, then
\begin{eqnarray}
	\label{caputo_alternative}
	_{a}{\mathrm{D}}_{t}^\alpha f(t) 
	= \frac{1}{\Gamma(1-\alpha)}\left(\frac{f(t)-f(a)}{(t-a)^\alpha}+\alpha	\int_a^t \frac{f(t)-f(\tau)}		{(t-\tau)^{\alpha+1}}\ \diff \tau\right).
\end{eqnarray}

\noindent
This alternative expression plays a central role in this work for the definition of the viscoelastic pseudo-potential of dissipation of Eq. \eqref{pseudod2}.

\section{Matrices for the Numerical Evaluation of the Motion Equation}

In the two-dimensional case, the matrices $\overline{\bm{S}}$, $\overline{\bm{F}}$ and $\bm{s}$, that appear in Eqs. \eqref{Eq_Residum_motion} and \eqref{Eq_jac_motion}, are given respectively by
\begin{equation*}
	\refstepcounter{equation} \latexlabel{h}
	\refstepcounter{equation} \latexlabel{Sbar}
	\overline{\bm{S}}=\begin{bmatrix}
	S_{11} & S_{12}   &  &\\  
	S_{12} & S_{22} &  &  \\  
	 & & S_{11} & S_{12}\\  
	 & & S_{12} & S_{22}
	 \end{bmatrix}, \quad
	 \overline{\bm{F}} = \begin{bmatrix}
	F_{11} & 0 & F_{21} &0 \\ 
	0 & F_{12} & 0 & F_{22} \\ 
	F_{12} & F_{11} & F_{22} & F_{21}
	\end{bmatrix},
	\tag{\ref{h}-\ref{Sbar}}
\end{equation*}
and
\begin{eqnarray}
	\label{s}
	\bm{s}=\left[S_{11}\quad S_{22}\quad S_{12}\right].
\end{eqnarray}


\bibliographystyle{spbasic}
\bibliography{mybibfile}

\end{document}